\documentclass{amsart}
\usepackage{amsmath}
\usepackage{amsmath}
\usepackage{amsfonts}
\usepackage{amsthm}
\usepackage{amssymb}
\usepackage{esint}
   \usepackage{times}
\usepackage{accents}
\usepackage{tikz-cd}
\usepackage{eucal}

\newtheorem{theorem}{Theorem}
\newtheorem*{theorem*}{Theorem}
\numberwithin{equation}{section}
\numberwithin{theorem}{section}

\newtheorem*{acknowledgement*}{Acknowledgement}
\newtheorem*{definition*}{Definition}

\newtheorem{definition}[theorem]{Definition}

\newtheorem{lemma}[theorem]{Lemma}

\newtheorem{proposition}[theorem]{Proposition}
\newtheorem{remark}[theorem]{Remark}

\newtheorem*{question*}{Question}

\newcommand{\RR}[0]{\mathbb{R}}

\newcommand{\ZZ}[0]{\mathbb{Z}}

\newcommand{\pd}[2]{\frac{\partial #1}{\partial#2}}

\newcommand{\pdt}[0]{\frac{\partial}{\partial t}}

\newcommand{\gt}[0]{\tilde{g}}
  \newcommand{\ft}[0]{\tilde{f}}
\newcommand{\gb}[0]{\overline{g}}

\newcommand{\delb}[0]{\overline{\nabla}}
\newcommand{\delt}[0]{\widetilde{\nabla}}
\newcommand{\delh}[0]{\hat{\nabla}}

\newcommand{\ve}[1]{\mathbf{#1}}
\newcommand{\grad}[0]{\operatorname{grad}}
\newcommand{\Rc}[0]{\operatorname{Rc}}
\newcommand{\Rm}[0]{\operatorname{Rm}}

\newcommand{\dfn}[0]{\doteqdot}
\newcommand{\Xb}[0]{\mathbf{X}}
\newcommand{\Yb}[0]{\mathbf{Y}}
\newcommand{\Xc}[0]{\mathcal{X}}
\newcommand{\Yc}[0]{\mathcal{Y}}

\newcommand{\gh}[0]{\hat{g}}
\newcommand{\pdtau}[0]{\pd{}{\tau}}
\newcommand{\srad}[0]{\mathcal{S}}
\newcommand{\rad}{\mathcal{R}}

\newcommand{\Rt}{\tilde{R}}
\newcommand{\Deltat}{\widetilde{\Delta}}
\newcommand{\Deltah}{\hat{\Delta}}
\newcommand{\Lc}[0]{\mathcal{L}}

\newcommand{\abs}[1]{\left\vert#1\right\vert}

\newcommand{\ra}{\rangle}
\newcommand{\la}{\langle}

\newcommand{\Id}[0]{\operatorname{Id}}

\newcommand{\Rmt}[0]{\widetilde{\Rm}}

\newcommand{\End}[0]{\operatorname{End}}

\newcommand{\Ac}[0]{\mathcal{A}}

\newcommand{\Dc}[0]{\mathcal{D}}

\newcommand{\Ec}[0]{\mathcal{E}}

\newcommand{\Tc}[0]{\mathcal{T}}
\newcommand{\Rmc}[0]{\mathcal{R}}

\newcommand{\Cc}[0]{\mathcal{C}}

\newcommand{\ch}[1]{\breve{#1}}

\newcommand{\Rh}[0]{\hat{R}}

\newcommand{\supp}[0]{\operatorname{supp}}

\newcommand{\Zc}[0]{\mathcal{Z}}
\newcommand{\trace}[0]{\operatorname{tr}}

\title[Terminally conical Ricci flows]{Ricci flows which terminate in cones}

 \author{Brett Kotschwar}
 \email{kotschwar@asu.edu}
 \address{School of Mathematical and Statistical Sciences,
 	Arizona State University, Tempe, AZ 85287, USA}


  \thanks{The author was partially supported by Simons Foundation grant \#709004.}

\begin{document}
\begin{abstract}
We prove that a complete solution to the Ricci flow on $M\times [-T, 0)$ which has quadratic curvature decay on some end of $M$
and converges locally smoothly to the end of a cone on that neighborhood as $t\nearrow 0$ must be a gradient shrinking soliton.
\end{abstract}
\maketitle

\section{Introduction}\label{sec:intro}

 The asymptotically conical shrinking Ricci solitons constructed in \cite{AngenentKnopf}, \cite{FeldmanIlmanenKnopf}, \cite{DancerWang}, and \cite{Yang}, all give rise to examples of complete ancient solutions to the Ricci flow
\begin{equation}\label{eq:rf}
 \partial_t g = -2\Rc(g)
\end{equation}
which flow smoothly into a cone as $t\nearrow 0$  off of a small singular set. The purpose of this paper is to determine more generally what kinds of solutions may flow into a cone at on some neighborhood of infinity on some end.

Intuitively, it would seem that a solution to a parabolic equation such as \eqref{eq:rf} could not simply acquire a property as rigid as a cone's scaling invariance in finite time without possessing some trace of this invariance at earlier times. However, one would also expect such solutions (if complete for $t < 0$) to become singular at points off of the end as $t\nearrow 0$, so the analytic basis for this intuition is less certain. It turns out, however, that the behavior of the solution off of the end and as it approaches the terminal time does not matter: as long as the solution flows into the cone on the end in a reasonably controlled way,
the entire solution must be shrinking self-similar at all previous times.

\subsection{Statement of results}

We will need to introduce some notation in order to state our main theorem.  Throughout this paper,
$(\Sigma, g_{\Sigma})$ will denote a compact $(n-1)$-dimensional smooth manifold and
\[
\Cc_a \dfn \Cc_a^{\Sigma} \dfn (a, \infty)\times \Sigma
\]
will denote the generalized cylinder over $\Sigma$.
We will also write
\[
   \hat{g} \dfn dr^2 + r^2g_{\Sigma}
\]
for the regular conical metric on $\Cc_0$.  By an \emph{end} (or \emph{end neighborhood}), we mean an unbounded connected component of the complement of a compact set.

\begin{theorem}\label{thm:terminalcone} Let $g(t)$ be a smooth complete solution to the Ricci flow on $M\times [-1, 0)$. Assume that there is $\rad> 0$, an end $V\subset M$, and a diffeomorphism $F:\Cc_{\rad} \longrightarrow V$
such that $\gt(t) = F^*g(t)$ converges locally smoothly to $\hat{g}$ on $\Cc_{\rad}$ as $t\longrightarrow 0$ and
satisfies
\[
    |\Rmt|(x, t) \leq K r^{-2}(x)
\]
on $\Cc_{\rad}\times [-1, 0)$ for some constant $K$.
Then there is  $f\in C^{\infty}(M)$ such that $g = g(-1)$ satisfies
\[
    \Rc(g) + \nabla\nabla f = \frac{g}{2},
\]
and $g(t) = -t\Phi_{t}^*g$ on $M\times [-1, 0)$, where
$\Phi_t$ is the  family of diffeomorphisms defined by
\[
    \pdt \Phi_t =  -\frac{1}{t}(\operatorname{grad}_{g})f\circ\Phi_t, \quad \Phi_{-1} = \operatorname{Id}.
\]
Moreover,
\[
  (F^{-1})_*\operatorname{grad}_{g}f = \frac{r}{2}\pd{}{r},
\]
on $\Cc_{\rad}$ and
$F^*(f \circ \Phi_t)$ converges locally smoothly to $\frac{r^2}{4}$ as $t\nearrow 0$.
\end{theorem}

In Theorem \ref{thm:terminalcone}, we do not assume a priori that there exists a shrinking soliton asymptotic to the cone
$(\Cc_0, \hat{g})$, and neither do we require that the curvature of the solution is bounded off of the end $V$ (either uniformly on $M\times [-1, 0)$ or on individual time-slices $M\times \{t\}$). In particular, $g(t)$ may (and, in general, will) become singular off of $V$ as $t\nearrow 0$.

One side consequence of Theorem \ref{thm:terminalcone} is that a complete shrinking soliton with an asymptotically conical end must be \emph{gradient}. That is, if $(M, g)$ is asymptotic to a cone on some end and $X$ is a vector field on $M$ such that
\begin{equation}\label{eq:rs}
 \Rc(g) + \frac{1}{2}\mathcal{L}_Xg = \frac{g}{2},
\end{equation}
then $X = \nabla f + K$ for $f\in C^{\infty}(M)$ and some Killing vector $K$.
Naber \cite{Naber4D} has shown that every complete shrinker with bounded curvature tensor is gradient.

Without assuming the completeness of the slices $(M, g(t))$, we may still obtain a corresponding statement for the restriction of $g(t)$ to the end $V$; in fact, we will
prove Theorem \ref{thm:terminalcone} as a consequence of this local statement.

\begin{theorem}\label{thm:terminalcone2} Let $\rad > 0$ and assume that $g(t)$ is a smooth solution to the Ricci flow on $\Cc_{\rad}\times [-1, 0]$ satisfying $g(0) = \hat{g}$ on $\Cc_{\rad}$ and
\[
   |\Rm|(x, t) \leq K r^{-2}(x)
\]
on $\Cc_{\rad}\times [-1, 0]$ for some constant $K$.
Then there is $f\in C^{\infty}(\Cc_{\rad})$ such that
\[
   \Rc(g) + \nabla\nabla f = \frac{g}{2} \quad\mbox{ and } \quad
   \operatorname{grad}_{g} f = \frac{r}{2}\pd{}{r}
\]
on $\Cc_{\rad}$, where $g = g(-1)$. Moreover, $g(t) = -t\Phi_{t}^*g$ on $\Cc_{\rad}\times [-1, 0)$, where $\Phi_t:\Cc_{\rad}\longrightarrow\Cc_{\rad}$ is the map
$\Phi_t(r, \sigma) = (r/\sqrt{-t}, \sigma)$.
\end{theorem}

As we mentioned above, the self-similar solutions associated to asymptotically conical shrinking solitons
give rise to solutions to \eqref{eq:rf} meeting the hypotheses of Theorems \ref{thm:terminalcone} and \ref{thm:terminalcone2}.
It is shown in Proposition 2.1 of \cite{KotschwarWangCylindrical} (cf. Proposition 2.1 in \cite{KotschwarWangIsometries}) that if $(M, g, f)$
is a normalized gradient shrinking soliton on $\Cc_{\rad}$ satisfying
\[
  \Rc(g) + \nabla\nabla f = \frac{g}{2} \quad \mbox{ and } \quad  R + |\nabla f|^2 = f
\]
and
which is $C^2$-asymptotic to $\hat{g}$ along $\Cc_{\rad}$ in the sense that $\lambda^{-2}\rho_{\lambda}^*g \longrightarrow \hat{g}$ in $C^2_{\mathrm{loc}}(\Cc_{\rad}, \hat{g})$ as $\lambda \longrightarrow \infty$,
where $\rho_{\lambda}:\Cc_0\longrightarrow \Cc_0$ is the dilation map
\begin{equation}\label{eq:rhodef}
 \rho_{\lambda}(r, \sigma) = (\lambda r, \sigma),
\end{equation}
then, after modifying by an appropriate diffeomorphism and (possibly increasing $\rad$), one can arrange for the following to be true:
\[
   \grad_g f = \frac{r}{2}\pd{}{r}
\]
on $\Cc_{\rad}$ and the family $g(-t) = -t \Phi_{t}^*g$, where $\Phi_{t}(r, \sigma) = (r/\sqrt{-t}, \sigma)$,
solves \eqref{eq:rf} on $\Cc_{\rad}\times [-1, 0)$ and satisfies that $g(t) \longrightarrow \hat{g}$ smoothly locally
as $t\nearrow 0$.
Thus, defining $g(0) = \hat{g}$, one obtains a smooth solution to \eqref{eq:rf} on
$\Cc_{\rad}\times [0, 1]$ which interpolates between the shrinker $g$ at $t=-1$ and the cone $\hat{g}$ at $t=0$
and which has uniform quadratic curvature decay.  Theorem \ref{thm:terminalcone2} implies that all solutions to \eqref{eq:rf} which flow smoothly into a cone with quadratic curvature decay arise in this way.

Currently, there are very few examples of \emph{complete} nontrivial asymptotically conical solitons which have been rigorously verified to exist; in fact, to our knowledge, the families constructed in the papers Angenent-Knopf \cite{AngenentKnopf}, Dancer-Wang \cite{DancerWang}, Feldman-Ilmanen-Knopf \cite{FeldmanIlmanenKnopf}, and Yang \cite{Yang} exhaust the known examples. Incomplete examples are more plentiful, however, (see, e.g., \cite{FeldmanIlmanenKnopf} and Appendix B of \cite{KotschwarWangConical}),
and such solitons give rise to smooth solutions satisfying the hypotheses of Theorem \ref{thm:terminalcone2}.
Stolarski \cite{Stolarski} has recently shown that every complete asymptotically conical shrinking soliton arises as a finite-time singularity model of a compact Ricci flow.

\subsection{Outline of the proof}

We first prove the local version of the statement in Theorem \ref{thm:terminalcone2}, and then apply a general continuation principle for soliton structures (Theorem \ref{thm:solext}) to help deduce the global statement in Theorem 1.1. The proof of Theorem \ref{thm:terminalcone2} consists of several steps which we outline here in their logical order (though not quite the order that they appear in the paper).
\subsubsection{Reframing the question as one of backward uniqueness}
The first step in the proof of Theorem \ref{thm:terminalcone2} is to distill the statement into a clean
expression of backward uniqueness. Recall that the conical metric $\hat{g} = dr^2+r^2g_{\Sigma}$ is characterized by the scaling law
\begin{equation}\label{eq:conechar}
    \rho_{\lambda}^*\hat{g} = \lambda^2\hat{g}, \quad \lambda > 0,
\end{equation}
it satisfies under pullback by the dilation map $\rho_{\lambda}:\Cc_{\rad}\longrightarrow\Cc_{\rad}$ defined above in \eqref{eq:rhodef}. At the same time, homothetical scaling and pull-back by diffeomorphism are also symmetries
of the Ricci flow equation. With the statement of Theorem \ref{thm:terminalcone2} in mind, observe that if
$g(t)$ is a solution to the Ricci flow on $\Cc_{\rad}\times [-1, 0]$ with $g(0) = \hat{g}$, then, for all $\lambda \geq 1$,
the family of metrics
\begin{equation}\label{eq:glambdadef}
    g_{\lambda}(t) = \lambda^{-2}\rho_{\lambda}^*g(\lambda^2t)
\end{equation}
solves the Ricci flow on $\Cc_{\rad}\times [-\lambda^2, 0]$ and, per \eqref{eq:conechar},
satisfies
\[
    g_{\lambda}(0) = \lambda^{-2}\rho_{\lambda}^*\hat{g} = \hat{g}.
\]
That is, the solutions $\{g_{\lambda}(t)\}_{\lambda \geq 1}$ share the common terminal value $\hat{g}$.

We thus seek an appropriate statement of backward uniqueness  to prove that any two of these solutions coincide where they are both defined. This would allow us to conclude that
\[
 g(t) \equiv g_{\lambda}(t) =\lambda^{-2}\rho_{\lambda}^*g(\lambda^2t)
\]
on $\Cc_{\rad}\times[-\lambda^{-2}, 0]$ for each $\lambda\geq 1$. We prove such a statement in Theorem \ref{thm:rfbu}.

\subsubsection{Embedding the problem into one for a larger system}
Before we continue with the overview of the proof of Theorem \ref{thm:terminalcone2}, let us explain how we prove Theorem \ref{thm:rfbu}.
In \cite{KotschwarRFBU, KotschwarRFBU2}, we have previously established the backward uniqueness of complete solutions to the Ricci flow \eqref{eq:rf} with uniformly bounded curvature. We cannot apply those results here, on account of the incompleteness of the solutions $(\Cc_{\rad}, g_{\lambda}(t))$, but we may still borrow from that work a the device to circumvent the gauge-degeneracy of the equation.

Since the Ricci flow is not a strictly parabolic system, the backward uniqueness of it solutions is not simply a consequence of the standard theory for parabolic equations. Neither can we apply the method of De Turck to transform our problem
into one for a strictly parabolic system as one can for the correspond problem of forward uniqueness (e.g., \cite{ChenZhu, HamiltonSingularities}). (Indeed, as explained in \cite{KotschwarRFBU}, the solutions to the harmonic map heat flow
required to transform the system are now solutions to ill-posed \emph{terminal-value} parabolic problems.)

Instead, following \cite{KotschwarRFBU}, we will embed the backward uniqueness problem at the heart
of Theorem \ref{thm:rfbu} into one for a larger, prolonged ``PDE-ODE'' system.
Here we follow the method in \cite{KotschwarRFBU} albeit with a system which is somewhat simpler than the one in that reference. This is carried out in Section \ref{sec:rfpdeode}, with some of the calculations deferred
to Appendix \ref{app:pdeode}.

\subsubsection{A general backward uniqueness theorem}
We have thus recast Theorem \ref{thm:rfbu} in terms of the vanishing of solutions to an ``approximately parabolic'' PDE-ODE system. That these solutions vanish under our hypotheses
is a consequence of Theorem \ref{thm:xysysbu}, a general backward uniqueness principle for solutions to such PDE-ODE systems valued in vector bundles over an evolving asymptotically conical end. Our proof uses the method of Carleman estimates, using generalizations of the inequalities in \cite{KotschwarWangConical} (and thereby, ultimately, of the estimates for linear parabolic equations in \cite{EscauriazaSereginSverak}).

There are some key differences in our situation, however. In \cite{KotschwarWangConical}, where the aim was to establish the uniqueness of shrinking solitons asymptotic to a given cone, the \emph{existence} of a shrinking soliton asymptotic to the cone was presumed explicitly throughout the work, and the analysis was carried out
with respect to an asymptotically conical \emph{self-similar} shrinking background. For Theorem \ref{thm:terminalcone2}, we do not wish to assume a priori that there is any shrinking soliton asymptotic to $(\Cc_0, \hat{g})$, and, consequently, we cannot simply use the estimates and analysis from  \cite{KotschwarWangConical} out of the box.

A review of the arguments in \cite{KotschwarWangConical} reveals that the self-similarity of the background solution is never used   in an critical way
(in contrast, e.g., to the arguments needed for the uniqueness result proven in \cite{KotschwarWangCylindrical} for the asymptotically cylindrical case). Nevertheless, the various weights in the Carleman inequalities are defined in terms of the soliton potential function $f$,
and many of the estimates of subsequent error quantities appeal to identities for $f$ and $g$ which are particular to shrinking solitons.
Thus in order to adapt these arguments, we must replace the weights and detangle the rest of the arguments
from the underlying soliton structure.

In the end, what this detangling reveals is that all that the arguments in \cite{KotschwarWangConical} really required of $f$ was that it was a good proxy for the squared radial distance function on the cone, and that all they required of the family of the background metrics $g(\tau)$ were that that they approximated the conical metric suitably well at infinity and as $t\nearrow 0$.   Thus we are able to prove estimates analogous to those in \cite{KotschwarWangConical} using the weights based on the radial distance
and using our terminally-conical solution $g(t)$ on $\Cc_{\rad}$ as family of background metrics.  The resulting statement in Theorem \ref{thm:xysysbu} generalizes the backward uniqueness principle in \cite{KotschwarWangConical};
see \cite{KotschwarWarpedProduct} for one application of this generalization.

We carry out the proof of of Theorem \ref{thm:xysysbu} over several sections. Our primary case of interest is when the asymptotically conical end is evolving by the Ricci flow, but the argument does not really require that the background metric be an exact solution to \eqref{eq:rf}, and, later, for our application, we will need to allow for slightly more general families of metrics in this role. Thus in Section \ref{sec:acmetrics} we consider \emph{initially conical} families
of metrics: family of metrics $g = g(\tau)$ which emanate from a cone at $\tau =0$ in an appropriately controlled way.  Here we consider $\tau$ to be a ``backward'' time parameter, using which allows
us to avoid a proliferation of negative signs.  Thus, rather than solutions to \eqref{eq:rf}
we consider solutions $g(\tau)$ to the \emph{backward Ricci flow}
\begin{equation}\label{eq:brf}
 \pdtau g = 2\Rc(g)
\end{equation}
which we obtain from the former by the change of variables $\tau = T-t$. In Section \ref{sec:acmetrics} we derive some elementary estimates on the decay of these
initially conical metrics and on certain derivative quantities involving the radial distance.

With these decay estimates in hand, we prove in Section \ref{sec:carleman} two sets of Carleman estimates for sections of vector bundles over an evolving initially conical end. These estimates, modeled on those in \cite{EscauriazaSereginSverak}, are analogs of those in \cite{KotschwarWangConical}, but use weights constructed from the conical radial distance. The first of these sets of estimates prove that sections
which satisfy a suitable system of mixed differential inequalities and vanish initially must vanish identically provided that they also decay extremely rapidly in space.  The second set of estimates we prove allow us to establish that our sections of interest do indeed vanish sufficiently rapidly.

Finally, in Section \ref{sec:bu}, we combine these estimates and establish the general backward uniqueness result in Theorem \ref{thm:xysysbu}.
We note that Haskins-Khan-Payne \cite{HaskinsKhanPayne} have recently adapted the backward uniqueness principle in \cite{KotschwarWangConical}
in a different direction, for solutions to PDE-ODE systems coupled to the metric evolution associated to an
asymptotically conical shrinking self-similar solution to the $G_2$-Laplacian flow.
\subsubsection{A backward uniqueness theorem for terminally conical Ricci flows}
Now, we return to the outline of the proof of Theorem \ref{thm:rfbu}. Given two solutions $g(\tau)$, $\gt(\tau)$ to the backward Ricci flow on $\Cc_{\rad}\times [0, T]$
with quadratic curvature decay and $g(0) =\hat{g} = \tilde{g}(0)$, we build from them a prolonged system
as discussed above which satisfies a closed system of differential inequalities
suitable for application to Theorem \ref{thm:xysysbu}.
From this theorem we deduce that $g\equiv \gt$
on $\Cc_{\rad^{\prime}}\times [0, T^{\prime}]$ for some $\rad^{\prime} \geq \rad$ and $0 < T^{\prime} < T$.

For our application to the family $\{g_{\lambda}(\tau)\}_{\lambda\geq 1}$ discussed above, it is important that we know that
the two solutions $g$ and $\gt$ agree on \emph{all} of $\Cc_{\rad}\times [0, T]$. For $\tau\in [0, T)$, we may improve $\rad^{\prime}$ to $\rad$ via the slight improvement to the standard statement (\cite{Bando}, \cite{KotschwarRFAnalyticity}) of instantaneous spatial real-analyticity for the Ricci flow observed in \cite{KotschwarWangIsometries}, which asserts that a smooth solution to \eqref{eq:rf} on $M\times [0, T]$ is real-analytic with respect to a \emph{time-independent} atlas for all $t\in (0, T]$. Thus in our setting,
it follows that $g(\tau) \equiv \gt(\tau)$ on $\Cc_{\rad}\times [0, T^{\prime}]$.

To see that $g$ and $\gt$ must agree on $\Cc_{\rad}\times [0, T]$ requires a bit more work. (If the solutions were complete on $\Cc_{\rad}$, we could simply invoke the backward uniqueness theorem in \cite{KotschwarRFBU}.) We would like to simply translate the two solutions in time and invoke Theorem \ref{thm:xysysbu} iteratively to cover the entire time interval. However,
the background metric is not time-independent, and if we translate it, it will no longer emanate from the cone at the new initial time. Then our Carleman estimates (specifically Theorem \ref{prop:pdecarleman2}, by which we deduce the rapid decay of the system) will no longer be applicable. To get around this, we alter our background metric slightly in the interval $[0, T^{\prime}]$ in which we have already established that our $g$ and $\gt$ agree. The new background metric will no longer solve the backward Ricci flow on the entire interval, but will be initially and approximately conical in a sense (defined in Section \ref{sec:acmetrics}) which is sufficient for Theorem \ref{thm:xysysbu}. With these modifications, we can iterate to conclude that $g(\tau)\equiv \gt(\tau)$ on $\Cc_{\rad}\times [0, T]$. This gives the final, global statement in Theorem \ref{thm:rfbu}.

\subsubsection{Self-similarity}
Now at last we return to Theorem \ref{thm:terminalcone2}. Applying Theorem \ref{thm:rfbu}
to $g(t)$ and $g_{\lambda}(t)$,
we conclude that
\[
      g(t) = g_{\lambda}(t) = \lambda^{-2}\rho_{\lambda}^*g(\lambda^{2}t)
\]
on $\Cc_{\rad}\times [-\lambda^{-2}, 0]$ for all $\lambda \geq 1$. In particular, we see that
\begin{equation}\label{eq:selfsim}
        g(t) = -t \rho_{\frac{1}{\sqrt{-t}}}^*g(-1)
\end{equation}
on $\Cc_{\rad}\times [-1, 0)$. This says that $g(t)$ is a shrinking self-similar solution. Linearizing, we find that $g=g(-1)$ satisfies
\[
   \Rc(g) + \frac{1}{2}\mathcal{L}_{X}g = \frac{1}{2}g
\]
at $t=-1$ with $X = \frac{r}{2}\pd{}{r}$.

To complete the proof of Theorem \ref{thm:terminalcone2}, we only need to show that $X$ is gradient relative to $g(-1)$. For this, we note that the family of one-forms $A \dfn dX^{\flat}$
where $X^{\flat} = g(t)(X, \cdot)$ satisfies
\[
   (\partial_{t} - \Delta_{g(t)})A = \Rm \ast A
\]
and vanishes at $t=0$ (when $X$ is the gradient of $r^2/4$ relative to $\hat{g} = g(0)$).
Thus we can also attack this problem as one of backward uniqueness. Applying Theorem \ref{thm:xysysbu}
to the latter equation for $A$, we show in Theorem \ref{thm:xexact} that $X^{\flat}$
is closed for all $t\in [-1, 0]$. Then, using an explicit representation formula, we show
that it is exact on all of $\Cc_{\rad}$ relative to some potential $f = f(-1)$. As part of this argument we see that the pulled-back family of functions $f(t) = \Phi_{t}^*f$ satisfy that $-tf(t) \longrightarrow r^2/4$ locally smoothly on $\Cc_{\rad}$ as $t\nearrow 0$.

\subsubsection{Extending the soliton structure from $V$ to $M$}
Finally, to obtain Theorem \ref{thm:terminalcone} from Theorem \ref{thm:terminalcone2}, we need
to prove that the gradient shrinking soliton structure defined on the end $V$ of the manifold extends to the entire manifold. As Ricci solitons are real-analytic \cite{Ivey}, this can be approached as a problem of analytic continuation. In Theorem \ref{thm:solext}, we prove more generally that if $(M, g)$
is a connected, simply connected real-analytic manifold which satisfies
\[
   \Rc(g) + \frac{1}{2}\mathcal{L}_Xg = \frac{\lambda}{2}g
\]
for some vector field $X$ and constant  $\lambda$ on a connected open set $U$, then $X$ extends uniquely
to a vector field on $M$ relative to which $g$ is globally a Ricci soliton. Moreover, if $X$ is gradient on $U$, then so is its extension.
The proof is a modification of the classical extension principle
for Killing vector fields due to Nomizu \cite{NomizuKilling}.

Using Theorem \ref{thm:solext}, we prove first that the lift $\gt(t)$ of $g(t)$ to the universal cover
$\tilde{M}$ of $M$ is a shrinking gradient Ricci soliton for all $t\in [-1, 0)$. Finally we argue that the soliton potential $\ft$ on $\tilde{M}$ must descend to a function $f$ on $M$ which coincides with the potential already defined on $V$.

\begin{acknowledgement*} We would like to acknowledge the substantial contribution
of Lu Wang to this work.  A number of the results in this paper, particularly the Carleman estimates in Section \ref{sec:carleman}
and the general backward uniqueness theorem in Section \ref{sec:bu}, are generalizations of results from our earlier joint paper \cite{KotschwarWangConical}, and we borrow elsewhere from from the methods and and philosophy of that project.
\end{acknowledgement*}

\section{Asymptotically conical families of metrics}\label{sec:acmetrics}
Here, and in the rest of the paper, we will continue to use
$(\Sigma, g_{\Sigma})$ to denote a compact $(n-1)$-dimensional Riemannian manifold and write $\hat{g} = dr^2+ r^2g_{\Sigma}$ for the conical metric
on $\Cc_0$ where
\[
\Cc_a \dfn \Cc_a^{\Sigma} \dfn (0, \infty)\times \Sigma.
\]
We will also write
\[
    \Cc_{a, T} \dfn \Cc_a\times [0, T],
\]
and use $r$ to denote both the global radial coordinate on $\Cc_0$ and the radial distance function $r:\Cc_0\longrightarrow (0, \infty)$ which returns this coordinate, i.e., the projection onto the first factor.

To obtain the backward uniqueness statement in Theorem \ref{thm:rfbu}, we need to consider
a somewhat more general family of background metrics than solutions to \eqref{eq:brf}. It will be convenient to use a backward time parameter when working with these solutions. Thus, in this section, and the next, we will consider a smooth family of metrics $g = g(\tau)$ on $\Cc_{\rad, T}$ which satisfies
\begin{equation}\label{eq:genmetev}
 \partial_{\tau} g_{ij} = 2\beta_{ij}
\end{equation}
for some smooth family of two-tensors $\beta = \beta(\tau)$.

\begin{definition}\label{def:emsmooth}
We will say that a family of metrics $g(\tau)$ on $\Cc_{\rad, T}$ satisfying \eqref{eq:genmetev} \emph{emanates smoothly from $\hat{g}$} if
\begin{equation}\label{eq:geninitial}
  g(x, 0) \equiv \hat{g}(x)\quad \mbox{and} \quad \beta(x, 0)(\nabla r, \cdot) \equiv 0,
\end{equation}
on $\Cc_{\rad}$, and there is some constant $K > 0$ such that
\begin{equation}\label{eq:genquaddecay}
  |\nabla^{(m)}\partial^{(l)}_{\tau}\beta|(x, \tau) \leq \frac{K}{r^2(x)}
\end{equation}
on $\Cc_{\rad, T}$ for $0 \leq m \leq 3$ and $0 \leq l \leq 2$.
\end{definition}
Here and below, $|\cdot|= |\cdot|_{g(\tau)}$ and $\nabla = \nabla_{g(\tau)}$. We will also use the notation
\[
    x_* \dfn \min\{x, 1\}, \quad x^* \dfn \max\{x, 1\},
\]
and use $C = C(n)$ to denote a series of constants depending at most on $n$, and $N = N(K, \rad^*, T_*)$ a series
of positive constants potentially depending also on $K$, $\rad_*$, and $T^*$.

We will not require the full strength of the decay assumption \eqref{eq:genquaddecay} for our purposes later  (for example, we will not need to assume bounds on the third derivatives of $\partial_{\tau}\beta$ and $\partial^{(2)}_{\tau}\beta$); we have written it in this way for simplicity.  In our case of primary interest --- when $g(\tau)$ is a solution to the backward Ricci flow which emanates from a cone and has quadratic curvature decay --- we will have quadratic (and higher) bounds on the space-time derivatives of the speed $\beta = \Rc$ for all orders. See Proposition \ref{prop:rfemsmooth} and the discussion in Section \ref{sec:shrinkercase}.

\subsection{Basic decay estimates on the metric}
We now derive a few elementary consequences of the definition.
\begin{proposition}\label{prop:genmetricdecay}
Suppose that $g = g(\tau)$ emanates smoothly from $\hat{g}$ on $\Cc_{\rad, T}$.
Then there is a constant $N= N(n, K, \rad_*, T^*)$ such that
\begin{equation}\label{eq:genhckest}
   |\hat{\nabla}^{(m)}(g-\hat{g})|_{\hat{g}} \leq \frac{N\tau}{r^{2}}
\end{equation}
for $m=0, 1, 2, 3$.
In particular, the metrics $g(\tau)$ are uniformly equivalent to $\hat{g}$.
\end{proposition}
\begin{proof} As in Lemma 14.2 in \cite{Hamilton3D}, the
quadratic decay of the speed $\beta$ implies
 that the family of metrics $g(\tau)$ (including $\hat{g} = g(0)$) are uniformly equivalent on $\Cc_{\rad, T}$; in fact,
 \begin{equation}\
         e^{-\frac{2K\tau}{r^2}}\hat{g}(x) \leq g(x, \tau) \leq e^{\frac{2K\tau}{r^{2}}}\hat{g}(x).
 \end{equation}
The case $m=0$ of \eqref{eq:genhckest} follows then, as $|e^x-1|\leq C(a)|x|$ for $x\in [0, a]$. In particular, the uniform bounds
in \eqref{eq:genquaddecay} hold also for the fixed norms $|\cdot|_{\hat{g}}$ if the constant
is enlarged by a factor depending on $K$, $R_*$, and $T^*$.

Now write $h = g- \hat{g}$ and consider the fiberwise identity
\begin{equation}\label{eq:genhident}
      h = g(\tau) - g(0) = 2\int_0^{\tau}\beta(g(s))\,ds.
\end{equation}
Differentiating \eqref{eq:genhident} with respect to the (time-independent) connection $\delh$,
and using that
\[
  \nabla - \hat{\nabla} = g^{-1}\ast\hat{\nabla}h,
\]
where, in this proof, $A\ast B$ denotes a linear combination of contractions of $A\otimes B$ with respect to the conical metric $\hat{g}$,
it follows from the uniform equivalence proven above that
\begin{align*}
   |\hat{\nabla}h|_{\hat{g}} &\leq C\int_0^{\tau}\left(|g^{-1}|_{\hat{g}}|\hat{\nabla}h|_{\hat{g}}|\beta|_{\hat{g}} + |\nabla \beta|_{\hat{g}}\right)\,ds\\
                   &\leq \frac{N}{r^{2}}\left(\tau + \int_0^{\tau}|\hat{\nabla}h|_{\hat{g}}\, ds\right).
\end{align*}
Then the case $m=1$ of \eqref{eq:genhckest} then follows from Gronwall's inequality.

With these two cases and the identity
\[
 \delh\delh \beta = \nabla\nabla \beta
 + g^{-1}\ast\delh h \ast \nabla\beta+ g^{-2}\ast (\delh h)^{2} \ast \beta + g^{-1}\ast \delh\delh h \ast \beta ,
\]
we see that, similarly,
\begin{align*}
  |\hat{\nabla}\hat{\nabla}h|_{\hat{g}} &\leq C\int_0^{\tau}|\hat{\nabla}\hat{\nabla}\beta(s)|_{\hat{g}}\,ds \leq \frac{N}{r^{2}}\left(\tau + \int_0^{\tau}|\hat{\nabla}\hat{\nabla}h|_{\hat{g}}\,ds\right),
\end{align*}
from which the case $m=2$ follows. The case $m=3$ is proven in the same way.
\end{proof}

\begin{remark}
 It follows from Proposition \ref{prop:genmetricdecay} that when $g$ emanates smoothly from $\hat{g}$
we also have the bounds
\begin{equation}\label{eq:genhckest2}
|\nabla^{(m)}(g-\hat{g})| \leq Nr^{-2},
\end{equation}
and
 \begin{equation}\label{eq:genquaddecay2}
    |\delh^{(m)}\partial^{(l)}_{\tau}\beta|_{\hat{g}}(x, \tau) \leq Nr^{-2},
\end{equation}
expressed in terms of the connections $\nabla$ and $\delh$
for $l$, $m=0, 1, 2, 3$. We will use these alternative forms of the bounds \eqref{eq:genquaddecay}
and \eqref{eq:genhckest} freely below.
\end{remark}

\subsection{Decay estimates on the conical distance function}

With respect to the conical metric $\hat{g}$, the radial distance function $r$ satisfies the identities
\begin{equation}
   |\delh r|_{\gh} = 1, \quad \delh\delh r^2 = 2\hat{g}, \quad \Deltah r = \frac{n-1}{r},  \quad \mbox{and } \ \Rh_{ijkl}\delh^lr = 0,
\end{equation}
on $\Cc_0$. These identities will be satisfied with a small error for small $\tau$ and large $r$
on a family $g(\tau)$ which emanates smoothly from $\hat{g}$,
We will need to control this error.
To this end, we define
\begin{equation}\label{eq:eta1def}
 \eta_1 \dfn |\nabla r|^2 -1,
\end{equation}
\begin{equation}\label{eq:eta2def}
 \quad \eta_2 \dfn \nabla\nabla\frac{r^2}{2} - g,
\end{equation}
\begin{equation}\label{eq:eta3def}
  \eta_3 \dfn \Delta r - \frac{n-1}{r} = \frac{\operatorname{tr}_g{\eta_2} - \eta_1}{r},
\end{equation}
and
\begin{equation}\label{eq:eta4def}
 (\eta_4)_{ijk} \dfn R_{ijk}^l\nabla_l r,
\end{equation}
on $\Cc_{\rad, T}$ relative to $g = g(\tau)$.

\begin{proposition}\label{prop:genetaests} Suppose $g = g(\tau)$ emanates smoothly from
$\hat{g}$ on $\Cc_{\rad, T}$. There are constants $N = N(n, K, R_*, T^*)$
such that the bounds
\begin{align}
 \label{eq:eta1est}
     |\nabla^{(m)}\eta_1| + \tau |\partial_{\tau} \eta_1| &\leq \frac{N\tau^2}{r^2},\\
\label{eq:eta2est}
      |\nabla^{(m)}\eta_2| + \tau |\partial_{\tau} \eta_2| &\leq \frac{N\tau}{r},\\
\label{eq:eta3est}
  |\nabla^{(m)}\eta_3| + \tau |\partial_{\tau} \eta_3| &\leq \frac{N\tau}{r^2},
\end{align}
hold on $\Cc_{\rad, T}$
for $m=0, 1, 2$.
Additionally, we have
\begin{equation}
 \label{eq:eta4est}
  |\eta_4| \leq \frac{N\tau}{r^2},
\end{equation}
and, for $m= 2, 3, 4$, the bounds
\begin{equation}\label{eq:nablakr}
  |\nabla^{(m)} r| \leq \frac{N}{r},
\end{equation}
 on  $\Cc_{\rad, T}$ for $N$ with the same dependencies.
\end{proposition}

\begin{proof}
We start with  \eqref{eq:eta1est}, computing that
\[
  \partial_{\tau}\eta_1 = -2\beta(\nabla r , \nabla r),
  \quad \partial^{(2)}_{\tau} \eta_1 = -2\partial_{\tau}\beta(\nabla r, \nabla r) + 8\beta^2(\nabla r, \nabla  r),
\]
and
\[
 \partial^{(3)}_{\tau}\eta_1 = -2\partial_{\tau}\beta(\nabla r, \nabla r) + 40\partial_{\tau}\beta(\beta(\nabla r), \nabla r)
 - 48\beta^3(\nabla r, \nabla r).
\]
Using \eqref{eq:geninitial} we have
\[
  \eta_1(x, 0) = |\delh r|_{\hat{g}}^2(x) -1 = 0, \quad \partial_{\tau}\eta_1(x, 0) = -2\beta(x, 0)(\delh r,\delh r) = 0,
\]
and
\[
  \partial^{(2)}_{\tau}\eta_1(x, 0) = -2(\partial_{\tau}\beta)(x, 0)(\delh r, \delh r),
\]
so that, by Taylor's theorem,
\begin{equation}\label{eq:eta1ident}
  \eta_1(x, \tau) = -\tau^2\partial_{\tau}\beta(x, 0)(\delh r, \delh r) + \frac{1}{2}\int_0^{\tau}Q(x, s)^{ij}\delh_i r \delh_j r(\tau-s)^2\,ds,
\end{equation}
where
\[
Q^{ij} =
-2(\partial_{\tau}\beta)^{ij}+ 40 g_{ab}(\partial_{\tau}\beta^{ia}\beta^{bj} + \partial_{\tau}\beta^{jb}\beta^{ia})- 48(\beta^3)^{ij}.
\]

Now, by \eqref{eq:genquaddecay2} we have
\[
  |\delh^{(m)}\partial_{\tau}\beta| \leq \frac{K}{r^2}, \quad \mbox{and} \quad |\delh^{(m)}Q(\delh r, \delh r)| \leq \frac{N}{r^2},
\]
for $m=0, 1, 2, 3$, and given that
\[
    |\delh^{(m)}r|_{\hat{g}} \leq Nr^{-(m-1)}
\]
for all $m$, the latter implies that
\[
  |\delh^{(m)} (Q(\delh r, \delh r))| \leq N  |\delh^{(m)} (Q(\delh r, \delh r))|_{\hat{g}} \leq \frac{N}{r^2}.
\]

Thus we can differentiate \eqref{eq:eta1ident} with respect to the (time-independent) connection $\delh$ and estimate away the integral to obtain that
\[
   |\delh^{(m)}\eta_1| \leq \frac{N\tau^2}{r^2}
\]
for $m=0, 1, 2$,
for some constant $N$ on $\Cc_{\rad, T}$. Then the spatial estimates in \eqref{eq:eta1est} with respect to the connection $\nabla$ follow since
\[
|\nabla \eta_1| \leq N|\delh \eta_1|_{\hat{g}}, \quad \mbox{and} \quad
  |\nabla\nabla \eta_1| \leq N(|\delh h|_{\hat{g}}|\nabla \eta_1|_{\hat{g}} + |\delh \delh \eta_1|_{\hat{g}}),
\]
where $h = g -\hat{g}$,
in view of \eqref{eq:genhckest}.

For the estimate on the temporal derivative in \eqref{eq:eta1est}, we use \eqref{eq:eta1ident} to see that
\[
 \partial_{\tau}\eta_1(x, \tau) = -2\tau (\partial_{\tau}\beta)(x, 0)(\delh r, \delh r) + \int_0^{\tau}Q(x, s)^{ij}\hat{\nabla}_ir\delh_jr
 (\tau-s)\,ds,
\]
and then apply \eqref{eq:genquaddecay} and \eqref{eq:genhckest}.

For \eqref{eq:eta2est}, we first note that
\begin{align}
\begin{split}\label{eq:eta2ident}
 (\partial_{\tau} \eta_2)_{ij} &= -r\partial_{\tau}\Gamma_{ij}^k\nabla_k r - 2\beta_{ij}\\
 &=r(\nabla_m\beta_{ij} - \nabla_i \beta_{jm} - \nabla_j \beta_{im})g^{km}\nabla_kr - 2\beta_{ij}\\
 &\dfn r S_1(\nabla r)_{ij} + (S_2)_{ij},
\end{split}
 \end{align}
and
\begin{align}
\begin{split}\label{eq:eta2ident2}
 (\partial^{(2)}_{\tau} \eta_2)_{ij} &= -r\partial^{(2)}_{\tau}\Gamma_{ij}^k\nabla_k r - 2\partial_{\tau}\beta_{ij}\\
 &=r Q_{1}(\nabla r)_{ij} + (Q_2)_{ij},
\end{split}
 \end{align}
where
\[
   Q_1 = g^{-2}\ast \beta \ast \nabla \beta + g^{-1}\ast \nabla \partial_{\tau}\beta, \quad
   Q_2 = -2\partial_{\tau}\beta.
\]
Thus
\[
 \eta_{2}(x, \tau) = \tau (rS_1(x, 0)(\delh r) + S_2(x, 0)) + \int_0^{\tau}\left(r Q_1(x, s)^k\delh_k r + Q_{2}(x, s)\right)(\tau-s)\,ds.
\]
Differentiating this identity with respect to $\delh$ and using \eqref{eq:genhckest} and \eqref{eq:genquaddecay2}
we see that
\[
     |\delh^{(m)}\eta_2|_{\hat{g}} \leq \frac{N\tau}{r}, \quad m = 0, 1, 2,
\]
for some $N = N(n, K, R_*, T^*)$. With \eqref{eq:genhckest}, we obtain analogous
estimates on $\nabla \eta_2$ and $\nabla\nabla \eta_2$. For the estimate on the temporal derivative in \eqref{eq:eta2est}, one may apply \eqref{eq:genquaddecay} directly to \eqref{eq:eta2ident}.

Next, differentiating the identity
\[
 \nabla\nabla r = \frac{g - \nabla r \otimes\nabla r + \eta_2}{r},
\]
we obtain the estimates in \eqref{eq:nablakr} inductively from \eqref{eq:eta2est}.  Then using \eqref{eq:eta1est} and \eqref{eq:eta2est} together with the identity
\[
  \eta_3 = \frac{\operatorname{tr}_g \eta_2 - \eta_1}{r},
\]
we obtain the estimates in \eqref{eq:eta3est}.

Finally, for \eqref{eq:eta4est}, note that
\[
    \nabla_{i}(\eta_2)_{jk} - \nabla_{j}(\eta_2)_{ik} = 2rR_{ijk}^l\nabla_lr =2r(\eta_4)_{ijk}
\]
so, by \eqref{eq:eta2est}, we see that
\[
    |\eta_4| \leq \frac{1}{r}|\nabla \eta_2| \leq \frac{N\tau}{r^2}
\]
as claimed.
\end{proof}

\section{Carleman estimates on ends smoothly emanating from a cone}\label{sec:carleman}

In this section, we will consider a family $g= g(\tau)$ of metrics satisfying \eqref{eq:genmetev}
on $\Cc_{\rad, T}$ for $R \geq 1$
and $0 < T \leq 1$ which emanates smoothly from $\hat{g}$. Our goal is to prove the two sets of Carleman estimates
for families of sections of vector bundles over $\Cc_{\rad, T}$ that will provide the basis of
the proof of our general backward uniqueness theorem, Theorem \ref{thm:xysysbu}, in the next section.

The arguments broadly follow those in Sections 4 and 5 of \cite{KotschwarWangConical}, with adaptations to the weight functions and the error estimates since we are no longer assuming the background metric is a self-similar shrinking solution to the backward Ricci flow. The basic idea is that (with the bounds on $g(\tau)$ derived in Section \ref{sec:acmetrics}) we can substitute $r$ for the function $h = 2\sqrt{\tau f}$
defined in \cite{KotschwarWangConical} in terms of the shrinker potential $f$.

\subsection{Uhlenbeck's trick and the operator $D_{\tau}$} It will simplify some of our calculations
if we perform them with respect to evolving $g(\tau)$-orthonormal frames. For this, we will use the formalism of Uhlenbeck's trick as it appears in \cite{BamlerBrendle}.  Let
$\pi:\Cc_{\rad, T}\longrightarrow \Cc_{\rad}$ denote the projection on to the first factor. The family of metrics $g = g(\tau)$
define a metric on the pull-back bundle $\pi^*(T\Cc_{\rad})$ which we will continue to denote by $g$. The sections of $\pi^*(T\Cc_{\rad})$ (``spatial vector fields'') correspond to families $V= V(\tau)$
of vector fields on $\Cc_{\rad}$ defined for $\tau\in [0, T]$.  The family of Levi-Civita connections of $g$
defines a connection $D$ on $\pi^*(T\Cc_{\rad})$ characterized by the conditions
\[
D_{\pdtau} V = \partial_{\tau}V + \Rc(V), \quad\mbox{ and } \quad D_X V = \nabla_X V,
\]
for all $X\in T\Cc_{\rad}$ and sections $V$ of $T\Cc_{\rad}$.  This connection is metric since
\[
  \pdtau \left(g(V, W)\right) = \langle D_{\pdtau}V, W\rangle + \langle V, D_{\pdtau}W\rangle.
\]
For simplicity, we will henceforth write
\[
   D_{\tau} V \dfn D_{\pdtau}V
\]
and use $\nabla_X V$ in place of $D_{X}V$ for differentiation with respect to tangent vectors $X\in T\Cc_{\rad}$.
The metric and connection on $\pi^*(T\Cc_{\rad})$ induce  metrics and metric connections on all of the pull-back bundles $\pi^*(T^k(T^*\Cc_{\rad}))\cong \bigotimes^k (\pi^*(T\Cc_{\rad}))^*$.

One can also regard $D_{\tau}$ as an operator on families of $k$-tensor fields on $\Cc_{\rad}$
which acts by
\[
   D_{\tau} V_{\alpha_1 \cdots \alpha_{k}} = \partial_{\tau} V_{\alpha_1 \cdots \alpha_{k}}
    - R_{\alpha_1}^{\delta}V_{\delta \alpha_2 \cdots \alpha_{k}} - \cdots - R_{\alpha_{k}}^{\delta}
    V_{\alpha_1 \cdots \alpha_{k -1} \delta}
\]
This coincides with the ``total $\tau$-derivative'' taken relative to evolving $g(\tau)$-orthonormal frames.
In particular, $D_{\tau} g = 0$.

\subsection{A general divergence identity}
Now let $\Zc$ denote a bundle of the form
\[
\Zc = \pi^*(T^k(T^*\Cc_{\rad}))
\]
over $\Cc_{\rad, T}$.
We will prove two sets of Carleman inequalities for sections of $\Zc$, which both, in effect, boil down to a weighted $L^2$-estimate on the commutator of the operators
\begin{equation}\label{eq:asdef}
    \Ac \dfn D_{\tau} - \nabla_{\nabla\phi} + \frac{F}{2}\Id,
  \quad\mbox{and}\quad \mathcal{S}\dfn \Delta + \nabla_{\nabla\phi} - \frac{F}{2}\Id,
\end{equation}
defined in terms of a choice of function $F\in C^{\infty}(\Cc_{\rad, T})$,
acting on sections $Z$ of $\Zc$.

The starting point for both estimates is the following divergence identity, which can be verified
by a long but straightforward computation.  It is effectively the  same identity as that in Lemma 5.1 of \cite{KotschwarWangConical} where $\beta = \Rc$, but for the use of $D_{\tau}$ instead of $\pd{}{\tau}$, which also eliminates certain terms involving $\beta$. Here and below, we write $d\mu = d\mu_{g(\tau)}$ and
 \begin{equation}\label{eq:bdef}
B \dfn \trace_{g}\beta,
 \end{equation}
so that, in particular, we have
\[
     \partial_{\tau} d\mu = B\, d\mu.
\]

\begin{lemma}\label{lem:dividentity2}
Let $G$ and $H$ be smooth positive functions on $\Cc_{\rad, T}$, and $\theta$ and $\sigma$ smooth positive functions on $[0, T]$.
Write $\phi = \log G$ as before.
Then for all smooth sections $Z$ of $\Zc$,
we have the identity
\begin{align}\label{eq:dividentity2}
  \begin{split}
  & \theta\sigma^{-\alpha}\nabla_i\bigg\{2\left\langle D_{\tau} Z, \nabla_i Z\right\rangle G
    + |\nabla Z|^2\nabla_i G -2\langle \nabla_{\nabla G} Z,  \nabla_i Z\rangle
    +\langle \nabla_i Z, Z\rangle GH\\
    &\quad+\frac{1}{2}|Z|^2(H\nabla_i G - G \nabla_i H) \bigg\}\;d\mu - \pdtau\left\{\left(|\nabla Z|^2 + \frac{1}{2}|Z|^2H\right)\theta\sigma^{-\alpha}G\;d\mu\right\}\\
    &= \bigg\{2\big\langle (D_\tau + \Delta)Z, \mathcal{A}Z\big\rangle  - 2|\mathcal{A}Z|^2
    -\bigg(2\nabla_i\nabla_j \phi + \frac{\dot{\theta}}{\theta}g_{ij}-2\beta_{ij}\bigg)  \langle \nabla_i Z, \nabla_j Z\rangle \\
   &\quad + \bigg(H - G^{-1}\left(\partial_{\tau} G - \Delta G\right) - B + \alpha\frac{\dot{\sigma}}{\sigma}\bigg)\left(|\nabla Z|^2 + \frac{H}{2}|Z|^2\right)\\
    &\quad -\frac{1}{2}\left(\partial_{\tau} H + \Delta H + \frac{\dot{\theta}}{\theta}H\right)|Z|^2
    + E_{\phi}(Z, \nabla Z)
    \bigg\} \theta \sigma^{-\alpha} G\;d\mu,
  \end{split}
\end{align}
where $\phi = \log G$ and
\begin{equation}\label{eq:edef}
 E_{\phi}(Z, \nabla Z) \dfn 2\big\langle[\nabla_i, D_{\tau}] Z -\beta_{ip}\nabla_p Z, \nabla_i Z\big\rangle + 2\nabla_j \phi \big\langle [\nabla_i, \nabla_j]Z, \nabla_i Z\big\rangle.
\end{equation}
\end{lemma}

Let us put the commutator terms in $E_{\phi}$ into a more explicit form. We have (compare, e.g., Appendix F of \cite{RFV2P2} for the case $\beta = \Rc$)
\begin{equation}\label{eq:dtaudericomm}
  [\nabla_i, D_{\tau}]Z = \left(\nabla_b \beta_{ip} - \nabla_p \beta_{ib}\right)\Lambda^b_p Z + \beta_{ip}\nabla_p Z
\end{equation}
and
\begin{equation}\label{eq:didjcomm}
  \nabla_j\phi [\nabla_i, \nabla_j] Z = \nabla_j\phi R_{jibp} \Lambda^b_p Z,
\end{equation}
where
$\Lambda^p_q$ is the operator
\begin{align*}
    \Lambda^p_q(Z_{\alpha}^{\beta}) &= \delta^p_{\alpha_1} Z_{q\alpha_2\cdots\alpha_\nu}^{\beta_1\beta_2\cdots\beta_\kappa}
	+ \delta^p_{\alpha_2}Z_{\alpha_1q\alpha_{3}\cdots\alpha_\nu}^{\beta_1\beta_2\cdots\beta_\kappa} + \cdots +
      \delta^p_{\alpha_\nu}Z_{\alpha_1\alpha_2\cdots q}^{\beta_1\beta_2\cdots\beta_\kappa}\\
      &\phantom{=}-\delta_q^{\beta_1} Z_{\alpha_1\alpha_2\cdots\alpha_\nu}^{p\beta_2\cdots\beta_\kappa}
	- \delta_q^{\beta_2} Z_{\alpha_1\alpha_2\cdots\alpha_\nu}^{\beta_1 p \cdots\beta_\kappa} - \cdots -
      \delta_q^{\beta_\kappa} Z_{\alpha_1\alpha_2\cdots\alpha_\nu}^{\beta_1\beta_2\cdots p}.
\end{align*}
In our applications, $\phi = \log G$ will have the form
\[
\phi(x, \tau) = \varphi(r(x), \tau)
\]
for some smooth $\varphi(r, \tau)$. In this case,
\begin{align}\label{eq:comm1}
 \begin{split}
   |E_{\phi}(Z, \nabla Z)| &\leq C\left(|\nabla\beta| + |\partial_r\varphi||\eta_4|\right)|\nabla Z| |Z|
   + C|\beta||\nabla Z|^2
 \end{split}
\end{align}
for some $C = C(n)$. Here $\eta_4 = \Rm(\nabla r, \cdot, \cdot, \cdot)$ is as defined in \eqref{eq:eta4def}.

\subsection{Carleman estimates to imply backward uniqueness}

For our first pair of Carleman estimates, we will take $\sigma \equiv 1$ and $\theta\equiv 1$.
Taking $G$ to be some fixed but as-yet-unspecified smooth positive function,
we set
\[
  H = F\dfn G^{-1}\left(\partial_{\tau} G - \Delta G + BG\right),
\]
so that the penultimate line of \eqref{eq:dividentity2} vanishes, leaving us with
\begin{align}\label{eq:dividentity2a}
  \begin{split}
  &\nabla_i\bigg\{2\left\langle D_{\tau} Z, \nabla_i Z\right\rangle G
    + |\nabla Z|^2\nabla_i G -2\langle \nabla_{\nabla G} Z,  \nabla_i Z\rangle
    +\langle \nabla_i Z, Z\rangle GF\\
    &\quad+\frac{1}{2}|Z|^2(F\nabla_i G - G \nabla_i F) \bigg\}\;d\mu - \pdtau\left\{\left(|\nabla Z|^2 + \frac{1}{2}|Z|^2F\right)G\;d\mu\right\}\\
    &= \bigg\{2\big\langle (D_\tau + \Delta)Z, \mathcal{A}Z\big\rangle  - 2|\mathcal{A}Z|^2
    -2\bigg(\nabla_i\nabla_j \phi -\beta_{ij}\bigg)  \langle \nabla_i Z, \nabla_j Z\rangle \\
    &\quad -\frac{1}{2}\left(\partial_{\tau} F + \Delta F\right)|Z|^2
    + E_{\phi}(Z, \nabla Z)
    \bigg\}G\;d\mu.
  \end{split}
\end{align}
Applying this identity to a smooth section $Z$ of $\mathcal{Z}$ which has compact support in space and vanishes
at $\tau =0$, we obtain the following consequence. (Cf. Lemma 4.3 of \cite{KotschwarWangConical}.)

\begin{proposition}\label{prop:pdecarlemanineq1} Let $G$ be any smooth positive function on $\Cc_{\rad, T}$. For any
 compactly supported section $Z$ of $\mathcal{Z}$ with $Z(\cdot, 0) \equiv 0$, we have the inequality
\begin{align}\label{eq:pdecarlemanineq1}
\begin{split}
&\frac{1}{2}\int_{\Cc_{\rad, T}} \!\left|D_{\tau}Z + \Delta Z\right|^2 G\, d\mu\, d\tau
    +\int_{\Cc_{\rad}\times\{T\}}\!\!\left(|\nabla Z|^2 + \frac{F}{2}|Z|^2\right)\!G\, d\mu\\
&\qquad \geq\int_{\Cc_{\rad, T}} \!\!\big(Q_1(\nabla Z,\nabla Z) + Q_2(Z, Z) + E_{\phi}( Z, \nabla Z)\big) G\,d\mu\,d\tau,
\end{split}
\end{align}
where
\[
  \phi = \log G, \quad H = G^{-1}\left(\partial_{\tau} G - \Delta G + BG\right),
\]
 and
\begin{align}
\begin{split}\label{eq:q1q2def}
  Q_1(\nabla Z, \nabla Z) &= 2(\nabla_i\nabla_j\phi -\beta_{ij})\langle\nabla_i Z, \nabla_j Z\rangle,\\
  Q_2(Z, Z) &= \frac{1}{2}\left(\partial_{\tau}F + \Delta F\right)|Z|^2,
\end{split}
\end{align}
where $E_{\phi}$ is as defined in \eqref{eq:edef}.
\end{proposition}
\begin{proof}
 Given our assumptions on $Z$, integrating the left-hand side of \eqref{eq:dividentity2a} over $\Cc_{\rad, T}$ leaves just one term, namely,
 \[
     -\int_{\Cc_{\rad}\times\{T\}}\!\!\left(|\nabla Z|^2 + \frac{F}{2}|Z|^2\right)\!G\, d\mu,
 \]
On the other hand, integrating the right-hand side yields
\begin{align*}
 & 2\int_{\Cc_{\rad, T}}\left(\big\langle D_\tau Z + \Delta Z, \mathcal{A}Z\big\rangle  - |\mathcal{A}Z|^2\right)\,d\mu\,d\tau\\
&\qquad   -\int_{\Cc_{\rad, T}}\left(Q_{1}(\nabla Z, \nabla Z) + Q_{2}(Z, Z) + E_{\phi}(Z, \nabla Z)\right)\,d\mu\,d\tau.
\end{align*}
and the inequality \eqref{eq:pdecarlemanineq1} then follows from Cauchy-Schwarz.
\end{proof}

We will also need the following elementary inequality (cf. Lemma 4.4 of \cite{KotschwarWangConical}).
\begin{proposition}
\label{prop:odecarlemanineq1} Let $G$ be any smooth positive function on $\Cc_{\rad, T}$.
 There exists a constant $N=N(n, K)$ such that if $Z$ is a compactly supported section of $\Cc_{\rad, T}$ with $Z(\cdot, 0)\equiv 0$, then
\begin{equation}
\label{eq:odecarlemanineq}
-\int_{\Cc_{\rad, T}} \left(N+\partial_{\tau}\phi\right)|Z|^2G\, d\mu\, d\tau\le
\int_{\Cc_{\rad, T}} \left|D_{\tau} Z\right|^2 G\, d\mu\, d\tau,
\end{equation}
where $\phi = \log G$.
\end{proposition}

\begin{proof}
Note that
\begin{equation}
\label{eq:tweightnorm}
\frac{\partial}{\partial \tau}\left(\abs{Z}^2\,d\mu\right)=\left(2\left\la D_{\tau} Z, Z\right\ra +\partial_{\tau} \phi|Z|^2+B|Z|^2\right)G\,d\mu.
\end{equation}
Since $|B| \leq C(n)K/r^2$ by \eqref{eq:genquaddecay}, the inequality (\ref{eq:odecarlemanineq}) follows by integrating \eqref{eq:tweightnorm} over $\Cc_{\rad, T}$ and applying the Cauchy-Schwarz inequality.
\end{proof}

\subsection{The weight function}

We now specify the weight function we will use in our first Carleman inequality. For fixed $\alpha > 0$, $0 < \delta < 1$, and $\tau_0\in (0, T)$, we define
\begin{equation}
\label{eq:g1def}
G_1 \dfn G_{1; \alpha, \delta, \tau_0} \dfn \exp(\alpha(\tau_0 - \tau)r^{2-\delta} + r^2),
\end{equation}
and
\begin{equation}\label{eq:phi1def}
\phi_1 \dfn \phi_{1; \alpha, \delta, \tau_0} \dfn \log G_{1; \alpha, \delta, \tau_0}.
\end{equation}

\begin{proposition}\label{prop:phi1prop}
 The function $\phi_1 = \log G_1$ satisfies
\begin{align}
\label{eq:phi1dtau}
 \partial_{\tau}\phi_1&= -\alpha r^{2-\delta},\\
\label{eq:phi1grad}
 \nabla \phi_1 &= \left(\alpha(2-\delta)(\tau_0 - \tau)r^{1-\delta} + 2r\right)\nabla r,\\
\begin{split}
 \label{eq:phigradsq}
 |\nabla \phi_1|^2 &= \big(4r^2 + 4\alpha(2-\delta)(\tau_0 -\tau) r^{2-\delta} \\
 &\phantom{=}\qquad+ \alpha^2(2-\delta)^2(\tau_0 -\tau)^2r^{2-2\delta}\big)(1+ \eta_1),
\end{split}\\
\begin{split}\label{eq:phi1hess}
 \nabla\nabla \phi_1 &= 2(g +\eta_2) +\frac{\alpha(2-\delta)(\tau_0 - \tau)}{r^\delta}\left(g -\delta\nabla r\otimes \nabla r + \eta_2\right),
\end{split}\\
 \begin{split}\label{eq:phi1lap}
  \Delta \phi_1&= 2(n + \eta_1 +  r\eta_3) + \frac{\alpha(2-\delta)(\tau_0 - \tau)}{r^{\delta}}\left(n - \delta + (1-\delta)\eta_1 + r \eta_3\right),
 \end{split}
\end{align}
where the quantities $\eta_i$ are as defined in \eqref{eq:eta1def} - \eqref{eq:eta3def}.
In particular, there is $\srad \geq \rad$ depending on  $n$, $\delta$, and $K$ such that
\begin{equation}\label{eq:phi1hesslowerbound}
  \nabla\nabla \phi_1 \geq (1+\delta)g
\end{equation}
on $\Cc_{\srad, T}$.
\end{proposition}

\begin{proof}
 The identities \eqref{eq:phi1dtau} and \eqref{eq:phi1grad} follow directly from the definition,
 and \eqref{eq:phigradsq} follows immediately from \eqref{eq:phi1grad}.

 For \eqref{eq:phi1hess}, we compute
 from \eqref{eq:phi1grad} that
\begin{align}\label{eq:hessphi1alt}
 \nabla\nabla\phi_1 &= \alpha(2-\delta)(\tau_0 - \tau)r^{-\delta}\left((1-\delta)\nabla r\otimes\nabla r + r\nabla\nabla r\right) + \nabla\nabla r^2\\
 \nonumber &= \alpha(2-\delta)(\tau_0 - \tau)r^{-\delta}\left(g + \eta_2 -\delta\nabla r\otimes\nabla r \right) + 2(g + \eta_2),
\end{align}
whereas for \eqref{eq:phi1lap}, we take the trace of \eqref{eq:hessphi1alt} and compute that
\begin{align*}
 \Delta \phi_1 &= \alpha(2-\delta)(\tau_0 - \tau)r^{-\delta}\left((1-\delta)|\nabla r|^2 + r \Delta r \right) + 2(r\Delta r + |\nabla r|^2)\\
\begin{split}
               &= \alpha(2-\delta)(\tau_0 - \tau)r^{-\delta}\left(n - \delta + (1-\delta)\eta_1 + r \eta_3\right)
                + 2(n + \eta_1 +  r\eta_3).
\end{split}
\end{align*}

For the lower bound \eqref{eq:phi1hesslowerbound}, we use \eqref{eq:eta1est} and \eqref{eq:eta2est} to choose $\srad$ sufficiently large to ensure that $|\eta_1|$, $|\eta_2| < (1-\delta)/2$ on $\Cc_{\srad, T}$. Then, the first term on the right of \eqref{eq:phi1hess} will be greater than or equal to $g$, while the second term will be nonnegative definite. Indeed, if $V = \nabla r/|\nabla r|$, then
\[
 (g -\delta\nabla r\otimes \nabla r + \eta_2)(V, V) \geq \frac{1}{2}(1+\delta) - \delta(1 + |\eta_1|) \geq 0,
\]
while, if $V$ is perpendicular to $\nabla r$, then
\[
  (g -\delta\nabla r\otimes \nabla r + \eta_2)(V, V) = (g +\eta_2)(V, V) \geq \frac{1}{2}(1+\delta)|V|^2\geq 0.
\]
\end{proof}

Next, we estimate the function
\begin{equation}\label{eq:f1def}
  F_1 \dfn G_1^{-1}(\partial_{\tau} - \Delta)G_1 + B.
\end{equation}
With Proposition \ref{prop:phi1prop}, the following bounds replace those in Lemmas 4.7 and 4.8 of \cite{KotschwarWangConical}
for the functions $G_1$ and $F_1$ defined in terms of the soliton potential of the background solution.
\begin{proposition}\label{prop:f1est}
There exist constants $\srad \geq \rad$ and $N > 0$ depending on $n$ and  $K$ such that, for all $\alpha > 0$, $0 < \delta < 1$, and $0 < \tau_0 \leq T$, we have

\begin{equation}\label{eq:f1c0est}
 0 \geq F_1 \geq -Nr^2\left(1 + \alpha r^{-\delta} + \alpha^2(\tau_0 - \tau)^2 r^{-2\delta}\right),
\end{equation}
and
\begin{equation}\label{eq:f1heatest}
  \pd{F_1}{\tau} + \Delta F_1 \geq 3\alpha r^{2-\delta} + \alpha^2(\tau_0 - \tau)r^{2-2\delta},
\end{equation}
 on $\Cc_{\srad, \tau_0}$.

\end{proposition}
\begin{proof} We will assume initially that $\srad \geq \rad$ and continue to increase $\srad$ as needed throughout the proof. For \eqref{eq:f1c0est}, we start from the identity
\[
  F_1 = \partial_{\tau} \phi_1  - \Delta\phi_1 - |\nabla \phi_1|^2 + B.
\]
By \eqref{eq:phi1dtau} the first term is nonpositive, and so too is the third.
From \eqref{eq:phi1lap} (or from \eqref{eq:phi1hesslowerbound}), we have $\Delta \phi_1 \geq n$
on $\Cc_{\srad, T}$ for $\srad$ large enough. Since $B \leq Nr^{-2}$ on $\Cc_{\rad, T}$, we see that
\[
 F_1 \leq -n + Nr^{-2} \leq -n/2,
\]
on $\Cc_{\srad, T}$ provided $\srad$ is taken larger still. This gives the upper bound in \eqref{eq:f1c0est}.

For the lower bound, we group terms according to like powers of $\alpha$, and use
\eqref{eq:phi1dtau}, \eqref{eq:phi1grad}, and \eqref{eq:phi1lap} to write
\[
 F_1 = H_0 + \alpha H_1 + \alpha^2 H_2,
\]
where
\begin{align*}
 H_0 &=  - 2n - 4r^2  -2(1+2r^2)\eta_1 - 2r\eta_3 + B,\\
H_1 &= -(1+4(2-\delta)(\tau_0 - \tau)(1+\eta_1))r^{2-\delta}\\
&\phantom{=-}
-r^{-\delta}(2-\delta)(\tau_0 - \tau)(n - \delta + (1-\delta) \eta_1  +r\eta_3),\\
 H_2 &= -(2-\delta)^2(\tau_0 - \tau)^2 r^{2-2\delta}(1+\eta_1).
\end{align*}
According to the bounds in  \eqref{eq:genquaddecay}, \eqref{eq:eta1est}, and \eqref{eq:eta3est} we have that
\begin{align*}
 H_0 \geq -Nr^2, \quad H_1 \geq -N r^{2-\delta}, \quad \ \mbox{and} \ \quad H_2 \geq -N(\tau_0 - \tau)^2 r^{2-2\delta},
\end{align*}
on $\Cc_{\srad, \tau_0}$ for some $N = N(n, K)$, provided $\srad$ is taken large enough. Combining these yields the lower bound in
\eqref{eq:f1c0est}.

For \eqref{eq:f1heatest}, we define
\[
  J_i = \partial_{\tau}H_i + \Delta H_i, \quad i=0, 1, 2,
\]
and estimate each $J_i$ separately. First, we compute that, for any $\beta > 0$,
\begin{align*}
 \Delta r^{\beta} &= \beta r^{\beta -1}\Delta r + \beta(\beta - 1)r^{\beta-2}|\nabla r|^2\\
 &=\beta r^{\beta-2}(n + \beta - 2) + \beta r^{\beta-2}\left((\beta -1) \eta_1 + r\eta_3\right)\\
 &= \beta r^{\beta-2}(n + \beta -2) + O(r^{\beta-3})
 \end{align*}
 in view of \eqref{eq:eta1est} and \eqref{eq:eta3est}. Here we use the big-O notation $O(r^{a})$ to denote a term $Q$ for which
 $|Q(r, \sigma, \tau)|\leq Nr^{a}$ on $\Cc_{\srad, \tau_0}$ as $r\longrightarrow \infty$.

In particular, $\Delta r^2 = 2n + O(r^{-1})$ and $\Delta r = O(r^{-1})$. Since the quantities
$|\partial_{\tau}\eta_i|$, $|\Delta\eta_i|$,
and $|\nabla\eta_i|$ for $i=1$ and $i=3$ are at least $O(r^{-2})$ and
\begin{equation}\label{eq:bev}
|(\partial_{\tau} + \Delta) B| \leq C(|\partial_{\tau}\beta| + |\beta|^2) = O\left(r^{-2}\right),
\end{equation}
applying the backward heat operator $\partial_{\tau}+ \Delta$ to the expression above for $H_0$,
we see that all terms in $J_0$ decay at least at the rate $O(r^{-1})$ with the exception of
\[
-4\Delta r^{2} = O(1), \quad \mbox{and} \quad -(1+4r^2)(\partial_{\tau} + \Delta)\eta_1 = O(1).
\]
Hence
\begin{align*}
J_0 & \geq -N
\end{align*}
on $\Cc_{\srad, \tau_0}$.

Computing similarly, we see that
\[
  J_1 \geq 4(2-\delta)r^{2-\delta} -Nr^{-\delta}
\]
and
\[
 J_2 \geq 2(2-\delta)^2(\tau_0 - \tau)r^{2-2\delta} - N(\tau_0 -\tau)r^{-2\delta}.
\]
Thus, since $\alpha \geq 1$ and $\delta \in (0, 1)$, we can choose $\srad$
sufficiently large to ensure that
\begin{align*}
 (\partial_{\tau} + \Delta)F_1 &= J_0 + \alpha J_1 + \alpha^2 J_2 \\
 &\geq 3\alpha(2-\delta) r^{2-\delta} + \alpha^2(2-\delta)^2(\tau_0-\tau)r^{2-2\delta}
\end{align*}
on $\Cc_{\srad, \tau_0}$.
\end{proof}

\subsection{Carleman estimates to imply backward uniqueness}

We now derive  a  Carleman estimate from the inequality in Proposition \ref{prop:pdecarlemanineq1}, using
 Proposition \ref{prop:f1est} to bound from below the quadratic forms $Q_1$ and $Q_2$ which appear in that inequality.
We will use the estimate in the proof of Theorem \ref{thm:xysysbu} to verify the vanishing of the elements
of the system.

Here and below, we will write $\|\cdot\|_{\Omega}$ for the space-time (respectively, spatial) $L^2$-norm induced by the time-dependent
family of measures $d\mu = d\mu_{g(\tau)}$ on $\Omega \subset \Cc_{\rad, T}$ (respectively, $\Omega \subset \Cc_{\rad}\times\{\tau\}$). This estimate is analogous to Proposition 4.9 in \cite{KotschwarWangConical}.

\begin{proposition}\label{prop:carlemanest1}
 Let $0< \delta< 1$, $\rad\geq 1$, and $0 < T \leq 1$. There exists $\srad \geq \rad$ depending on $n$, $\delta$, and $K$, such that, for all $\alpha \geq 1$ and $0 < \tau_0 < T$,
 and all smooth sections $Z$ of $\mathcal{Z}$ with compact support in $\Cc_{\srad, T}$ which vanish on $\Cc_{\srad}\times\{0\}$, we have the estimates
\begin{align}
\label{eq:pdecarlemanest1}
\begin{split}
&\alpha\|ZG_1^{\frac12}\|^2_{\Cc_{\srad, \tau_0}}+ \|\nabla Z G_1^{\frac12}\|_{\Cc_{\srad, \tau_0}}^2
\\
  &\qquad\le \frac{1}{2}\|(D_{\tau} +\Delta) ZG_1^{\frac12}\|_{\Cc_{\srad, \tau_0}}^2
 + \|\nabla Z G_1^{\frac12}\|^2_{\Cc_{\srad}\times\{\tau_0\}}
\end{split}
\end{align}
and
\begin{align}
\label{eq:odecarlemanest1}
&\alpha\|Z G_1^{\frac12}\|_{\Cc_{\srad, \tau_0}}^2 \le \|D_{\tau}ZG_1^{\frac12}\|_{\Cc_{\srad, \tau_0}}^2,
\end{align}
where $G_1 = G_{1; \alpha, \tau_0}$.
\end{proposition}
\begin{proof}
Consider the right hand side of the inequality in \eqref{eq:pdecarlemanineq1}. Using \eqref{eq:phi1hesslowerbound}
and \eqref{eq:genquaddecay}, we have
\begin{align*}
  Q_1(\nabla Z, \nabla Z) &= 2(\nabla_i\nabla_j \phi_1 - \beta_{ij})\langle \nabla_i Z, \nabla_j Z\rangle\\
  &\geq \left(2(1+\delta) - \frac{N}{r^2}\right)|\nabla Z|^2 \geq 2|\nabla Z|^2
\end{align*}
on $\Cc_{\srad, T}$ for sufficiently large $\srad$. Combining this with Proposition \ref{prop:f1est},
it follows that by increasing $\srad$ further, we can achieve that
\begin{align*}
 Q_1(\nabla Z, \nabla Z) + Q_2(Z, Z) &\geq 2|\nabla Z|^2 + \left(\frac{3\alpha}{2} + \alpha^2(\tau_0 - \tau)r^{2-2\delta}\right)|Z|^2
\end{align*}
on $\Cc_{\srad, T}$.

It remains to estimate the the commutator quantity $E_{\phi_1}$. For this, recall that, by \eqref{eq:eest1}, we have
\begin{align*}
 |E_{\phi_1}(Z, \nabla Z)| &\leq C\left(|\nabla \beta| + |\partial_r\Phi_1||\eta_4|\right)|\nabla Z| |Z|  + C|\beta||\nabla Z|^2
\end{align*}
where $\Phi_1(r(x), \tau ) = \phi_1(x, \tau)$,
for some $C = C(n)$ .
Since
\[
  |\partial_{r} \Phi_1| = (\alpha(2-\delta)(\tau_0-\tau)r^{1-\delta} + 2r) \leq Nr(1+ \alpha(\tau_0-\tau)r^{-\delta}),
\]
it follows from \eqref{eq:genquaddecay} and \eqref{eq:eta4est} that
\begin{align*}
|E_{\phi_1}(Z, \nabla Z)| &\leq \frac{N\alpha^2(\tau_0-\tau)}{r^{1+2\delta}}|Z|^2 + \frac{N}{r}|\nabla Z|^2\\
                 &\leq \alpha^2(\tau_0 - \tau)r^{2-2\delta}|Z|^2 + |\nabla Z|^2
\end{align*}
on $\Cc_{\srad, \tau_0}$ provided $\srad$ is taken large enough.

Putting things together, we have
\[
Q_1(\nabla Z, \nabla Z) + Q_2(Z, Z) + E_{\phi_1}(\nabla Z, Z) \geq |\nabla Z|^2 + \alpha|Z|^2
\]
on $\Cc_{\srad, \tau_0}$, for suitably large $\srad$. The estimate \eqref{eq:pdecarlemanest1}
then follows from  \eqref{eq:pdecarlemanineq1}, using that $F_1 \leq 0$.

For \eqref{eq:odecarlemanest1}, note that by \eqref{eq:phi1dtau}, we have
\[
 \partial_{\tau}\phi_1 = -\alpha r^{2-\delta},
\]
so that the left-hand side of \eqref{eq:odecarlemanineq} is greater than $\alpha\|ZG^{1/2}\|^2_{\Cc_{\srad, \tau_0}}$ provided $\srad$ is taken large enough.
\end{proof}

\subsection{Carleman estimates to imply rapid decay}
The weight $G_1$ in the estimates in Proposition \ref{prop:carlemanest1} grows at an exponential-quadratic rate in the radial distance $r$. We can therefore apply them only to sections we know a priori to vanish decay even more rapidly. As in \cite{KotschwarWangConical}, we will use estimates modeled on (1.4) in \cite{EscauriazaSereginSverak} to establish this decay for sections satisfying the hypotheses of Theorem \ref{thm:xysysbu}.

We begin, as in the proof of the last estimate, with the divergence identity \eqref{eq:dividentity2}. Here,
however, we will need to make nontrivial choices of $\sigma$ and $\theta$.  Let us assume for the moment that $\sigma= \sigma(\tau)$ is smooth, positive, and increasing, and
$G = G(x, \tau)$ is smooth and positive. As in \cite{EscauriazaSereginSverak}, we set $\theta = \sigma/\dot{\sigma}$,
a choice which gives
\[
   \frac{\dot{\theta}}{\theta} \cdot \frac{\dot{\sigma}}{\sigma}  = - \ddot{\widehat{\log \sigma}}.
\]
Next, we take $H$ in the form $H = F - \alpha\dot{\sigma}/\sigma$ where (as before)
\[
   F =  G^{-1}\left(\partial_{\tau} G - \Delta G\right) +   B.
\]
With these choices, the quantity in the coefficient of $|Z|^2$ in the last line of \eqref{eq:dividentity2}
satisfies
\begin{align*}
\partial_{\tau}H + \Delta H + \frac{\dot{\theta}}{\theta} H &= \partial_{\tau}H
+ \Delta F +\frac{\dot{\theta}}{\theta}\left(F - \alpha \frac{\dot{\sigma}}{\sigma}\right) -\alpha \ddot{\widehat{\log \sigma}} \\
&= \partial_{\tau} F + \Delta F + \frac{\dot{\theta}}{\theta} F,
\end{align*}
and the penultimate line vanishes entirely, so that
we have
\begin{align}\label{eq:dividentity3}
  \begin{split}
  & \theta\sigma^{-\alpha}\nabla_i\bigg\{2\left\langle D_{\tau} Z, \nabla_i Z\right\rangle G
    + |\nabla Z|^2\nabla_i G -2\langle \nabla_{\nabla G} Z,  \nabla_i Z\rangle
    +\langle \nabla_i Z, Z\rangle GH\\
    &\quad+\frac{1}{2}|Z|^2(H\nabla_i G - G \nabla_i H) \bigg\}\;d\mu - \pdtau\left\{\left(|\nabla Z|^2 + \frac{1}{2}|Z|^2H\right)\theta\sigma^{-\alpha}G\;d\mu\right\}\\
    &= \bigg\{2\big\langle (D_\tau + \Delta)Z, \mathcal{A}Z\big\rangle  - 2|\mathcal{A}Z|^2
    -\bigg(2\nabla_i\nabla_j \phi + \frac{\dot{\theta}}{\theta}g_{ij}-2\beta_{ij}\bigg)  \langle \nabla_i Z, \nabla_j Z\rangle \\
    &\quad -\frac{1}{2}\left(\partial_{\tau} F + \Delta F + \frac{\dot{\theta}}{\theta}F\right)|Z|^2
    + E_{\phi}(Z, \nabla Z)
    \bigg\} \theta \sigma^{-\alpha} G\;d\mu,
  \end{split}
\end{align}
for any smooth section $Z$ of $\Zc$. Here $\phi = \log G$ as before.

Applying \eqref{eq:dividentity3} to a section $Z$ which vanishes identically at $\tau =0$
and has compact support in $\Cc_{\rad}\times [0, T)$, and integrating over $\Cc_{\rad, T}$, we may use
the Cauchy-Schwarz as in the proof of Proposition \ref{prop:pdecarlemanineq1} to obtain the following inequality. (Cf. Lemma 5.2 of \cite{KotschwarWangConical}.)

\begin{proposition}\label{prop:pdeintineq2}
Let $\sigma(\tau)$ be a smooth positive increasing function on $[0, T]$ and $G$ a smooth positive function on $\Cc_{\rad, T}$.
For any $\alpha > 0$ and any smooth section $Z$ of $\Zc$ which vanishes for $\tau=0$ and has compact support in $\Cc_{\rad} \times[0, T)$, we have the inequality
 \begin{align}\label{eq:pdeintident2}
 \begin{split}
 &\int_{\Cc_{\rad, T}}\frac{\sigma^{1-\alpha}}{\dot{\sigma}}
 	\left|D_{\tau}Z + \Delta Z\right|^2G\,d\mu\,d\tau\\
 &\qquad\geq \int_{\Cc_{\rad, T}}\,\frac{\sigma^{1-\alpha}}{\dot{\sigma}}\big(Q_3(\nabla Z, \nabla Z) + Q_4(Z, Z)
 + E_{\phi}(Z, \nabla Z)\big)G\,d\mu\,d\tau,
 \end{split}
 \end{align}
 where $\theta = \sigma/\dot{\sigma}$, $F = G^{-1}(\partial_{\tau} G - \Delta G) + B$,
 \begin{align}
 \begin{split}\label{eq:q3def}
   Q_3(\nabla Z, \nabla Z)  &= \left(2\nabla_i\nabla_j\phi
     +\frac{\dot{\theta}}{\theta} g_{ij}-2\beta_{ij}\right)\langle \nabla_i Z, \nabla_j Z\rangle,
\end{split}\\
\begin{split}\label{eq:q4def}
  Q_4(Z, Z) &= \frac{1}{2}\left(\partial_{\tau} F + \Delta F +\frac{\dot{\theta}}{\theta}F\right)|Z|^2,
\end{split}
\end{align}
and $E_{\phi}$ is as defined in \eqref{eq:edef}.
\end{proposition}

We will use the above inequality in combination with the following auxiliary ``absorbing'' inequality. (Cf. Lemma 5.3 in \cite{KotschwarWangConical}.)
\begin{proposition}\label{prop:pdeintineq3} For $\sigma$, $G$, and $Z$ as in Proposition \ref{eq:pdeintident2}, we have
\begin{align}
\begin{split}\label{eq:pdeintineq3}
  &\int_{\Cc_{\rad, T}}\,\sigma^{-\alpha}\left(\frac{\alpha\dot{\sigma}}{2\sigma}- F\right)|Z|^2G\,d\mu\,d\tau\\
  &\qquad \leq2\int_{\Cc_{\rad, T}}
    \,\sigma^{-\alpha}\left(|\nabla Z|^2+\frac{\sigma}{\alpha\dot{\sigma}}\left|D_{\tau}Z + \Delta Z\right|^2
	    \right)G\,d\mu\,d\tau
\end{split}
\end{align}
for all $\alpha > 0$.
\end{proposition}
\begin{proof}
The result follows from integrating the identity
\begin{align}\label{eq:dividentity4}
\begin{split}
 & \pdtau\bigg\{\sigma^{-\alpha}|Z|^2 G\,d\mu \bigg\} + \sigma^{-\alpha}\nabla_i\left(\nabla_i|Z|^2 G - |Z|^2\nabla_iG\right)\,d\mu\\
&\qquad=\sigma^{-\alpha}\bigg\{\left(\pdtau + \Delta\right)|Z|^2 + \left(F- \alpha\frac{\dot{\sigma}}{\sigma}\right)\bigg\} G\,d\mu
\end{split}
\end{align}
over $\Cc_{\rad, T}$
and using the estimate
\begin{align*}
    \left(\pdtau + \Delta\right) |Z|^2 &= 2 |\nabla Z|^2 + 2\langle (D_{\tau} + \Delta)Z, Z\rangle \\
          &\leq 2|\nabla Z|^2 +  \frac{\alpha\dot{\sigma}}{2\sigma}|Z|^2
          + \frac{2\sigma}{\alpha\dot{\sigma}}|D_{\tau}Z + \Delta Z|^2
\end{align*}
on the first term on the righthand side.
\end{proof}

\subsection{The second Carleman estimate for the PDE component}
Now we specialize these inequalities to our case of interest.  As in \cite{EscauriazaSereginSverak}, for fixed $a\in (0, 1)$ and $\rho \geq 1$, we define the functions
\begin{equation}\label{eq:sigmadef}
\sigma_a(\tau) \dfn (\tau+ a) e^{-\frac{\tau+a}{3}},
\end{equation}
and
\begin{equation}\label{eq:g2def}
 G_2(x, \tau) \dfn G_{2; a, \rho}(x, \tau) \dfn (\tau+a)^{-n/2}\exp{\left(-\frac{(r(x) - \rho)^2}{4(\tau+a)}\right)}.
\end{equation}
The function $\sigma_a$ is positive, increasing, and approximately linear for small $\tau$. We note for later that it satisfies the bounds
\begin{equation}\label{eq:sigmaest}
\frac{\tau + a}{3} \leq \sigma_a(\tau)\leq \tau+a, \quad\mbox{ and }\quad \frac{1}{10} \leq \dot{\sigma}_a(\tau) \leq 1,
\end{equation}
for $a$, $\tau\in [0, 1]$. Also, recalling \eqref{eq:pdeintident2}, with this choice for $\sigma$, we have the identity
\begin{equation}
\label{eq:sigmaident}
  \frac{\dot\theta}{\theta} = -\frac{\sigma_{a}}{\dot{\sigma_a}}\cdot \ddot{\widehat{\log \sigma_a}} = \frac{1}{\tau+a}\left(\frac{1}{1-(\tau+a)/3}\right),
\end{equation}
where $\theta = \sigma_{a}/\dot{\sigma_a}$ as above.

The function $G_2$ is an exact solution to the heat equation relative to the conical metric $\hat{g}$. On $\Cc_{\rad}$, relative to $g(\tau)$ it will be an approximate solution. To estimate
the error of this approximation, define
\[
\phi_2 \dfn \phi_{2; a, \rho} \dfn \log G_{2; a, \rho} = -\frac{(r- \rho)^2}{4(\tau+a)} - \frac{n}{2}\log(\tau+a),
\]
and define
\[
  F_2 = G_2^{-1}(\partial_{\tau} G_2 - \Delta G_2) + B.
\]

\begin{proposition}\label{prop:phi2ident} Let $0 < a < 1$ and $\rho > 0$.
  The function $\phi_2 = \phi_{2; a, \rho}$ satisfies
\begin{align}\label{eq:dtauphi2}
   \partial_{\tau} \phi_2 &= \frac{(r-\rho)^2}{4(\tau+a)^2} - \frac{n}{2(\tau+a)},\\
\label{eq:gradphi2}
   \nabla \phi_2 &= - \left(\frac{r-\rho}{2(\tau+a)}\right)\nabla r,\\
\label{eq:normgradphi2}
|\nabla \phi_2|^2 & = \frac{(r-\rho)^2}{4(\tau+a)^2}(1 + \eta_1),\\
\begin{split}\label{eq:hesslogG2ident}
   \nabla \nabla \phi_2 &= -\frac{1}{2(\tau + a)}\bigg\{\left(1- \frac{\rho}{r}\right)g  + \frac{\rho}{r}\nabla r \otimes \nabla r
   +\left(1 -\frac{\rho}{r}\right)\eta_2\bigg\},
\end{split}\\
\label{eq:lapphi2}
   \Delta \phi_2 &= -\frac{1}{2(\tau + a)}\bigg\{n- \frac{(n-1)\rho}{r} + \eta_1 + (r-\rho)\eta_3\bigg\}.
\end{align}

In particular, there is a constant $N$ depending on $n$ and $K$ (but independent of $\rho$ and $a$) such that
\begin{equation}\label{eq:hesslogG2est}
  \nabla\nabla \phi_2 \geq -\frac{g}{2(\tau+a)} - \frac{N}{r}g
\end{equation}
on $\Cc_{\rad, T}\cap \Cc_{\rho/16, T}$.
\end{proposition}
\begin{proof}
The identities \eqref{eq:dtauphi2} - \eqref{eq:lapphi2} follow by direct computation. For \eqref{eq:hesslogG2est}, we argue as for \eqref{eq:phi1hesslowerbound}. First consider a unit vector $V$ orthogonal to $\nabla r$. According to \eqref{eq:hesslogG2ident} and \eqref{eq:eta2est}, we have
\begin{align*}
 \begin{split}
   \nabla \nabla \phi_2(V, V) &= -\frac{1}{2(\tau + a)}\bigg(\left(1- \frac{\rho}{r}\right)
   +\left(1 -\frac{\rho}{r}\right)\eta_2(V, V)\bigg)\\
   &\geq -\frac{1}{2(\tau+a)} -\frac{N}{r}
\end{split}
\end{align*}
for $r > \rho/16$.
On the other hand, if $V = \nabla r/|\nabla r|$, then
\begin{align*}
\begin{split}\label{eq:hesslogG2ident}
   \nabla \nabla \phi_2(V, V) &=
   -\frac{1}{2(\tau + a)}\bigg\{ 1 + \frac{\rho}{r}\eta_1
   +\left(1 -\frac{\rho}{r}\right)\eta_2(V, V)\bigg\}\\
   &\geq -\frac{1}{2(\tau+a)} - \frac{N}{r},
\end{split}
\end{align*}
by \eqref{eq:eta1est} and \eqref{eq:eta2est}. Combining the two cases yields \eqref{eq:hesslogG2est}.

\end{proof}

Next, we estimate the function $F_2$. Together with those in Proposition \ref{prop:phi2ident}, the bounds below replace
those in Lemmas 5.3 - 5.4 of \cite{KotschwarWangConical} for the functions $G_2$ and $F_2$ defined in terms of the soliton potential.
\begin{proposition}\label{prop:phif2est}
Let $0 < a < 1$ and $\rho > 0$. There is a constant $N= N(n, K)$ such that
$F_2 = F_{2; a, \rho}$
satisfies the bounds
\begin{equation}\label{eq:f2est}
 \left|F_2 + \frac{(n-1)\rho}{2r(\tau+a)}\right| \leq N,
\end{equation}
and
 \begin{equation}\label{eq:f2bheat}
  \left|(\partial_{\tau} + \Delta)F_2 - \frac{(n-1)\rho}{2r(\tau+a)^2}\right| \leq \frac{N}{(\tau+a)},
\end{equation}
on the set $\Cc_{\rho/16, T}\cap\Cc_{\rad, T}$.
\end{proposition}
\begin{proof}

From Proposition \ref{prop:phi2ident}, we see that
\begin{align}
\begin{split}\label{eq:f2ident}
 F_2 &= \partial_{\tau}\phi_2 - \Delta \phi_2 - |\nabla \phi_2|^2 + B\\
 &= -\frac{(n-1)\rho}{2(\tau+a)r} + \left(\frac{1}{2(\tau+a)} - \frac{(r-\rho)^2}{4(\tau+a)^2}\right)\eta_1
 +\left(\frac{r-\rho}{2(\tau+a)}\right)\eta_3 + B.
\end{split}
\end{align}
We then obtain \eqref{eq:f2est} follows directly from the decay bounds in
\eqref{eq:genquaddecay}, \eqref{eq:eta1est}, and \eqref{eq:eta3est} on $B$, $\eta_1$, and $\eta_3$.

For \eqref{eq:f2bheat}, we argue as for \eqref{eq:f1heatest}. First of all, from \eqref{eq:f2ident},
we have
\begin{equation}
 F_2 = A_0 + A_1 \eta_1 + A_2 \eta_3  + B,
\end{equation}
where
\[
   A_0 = -\frac{(n-1)\rho}{2(\tau+a)r}, \quad A_1 = \frac{1}{2(\tau+a)} - \frac{(r-\rho)^2}{4(\tau+a)^2}, \quad \mbox{and} \quad
  A_2 = \frac{r-\rho}{2(\tau+a)}.
\]
Since (by Proposition \ref{prop:genetaests}, as before)
\begin{align}
\begin{split}
\label{eq:deltarbeta} \Delta r^{\beta} &= \frac{\beta(n+\beta -2)}{r^{2-\beta}} + O\left(\frac{\tau}{r^{3-\beta}}\right),
\end{split}
\end{align}
for any $\beta\neq 0$, we see that
\begin{align}\label{eq:a0ev}
 \begin{split}
    (\partial_{\tau} + \Delta)A_0 &= \frac{(n-1)\rho}{2(\tau+a)r}\left(\frac{1}{\tau+a} + \frac{n-3}{r^2}
    + O\left(\frac{\tau}{r^3}\right)\right)\\
    &=  \frac{(n-1)\rho}{2r(\tau+a)^2} + O\left(\frac{1}{(\tau+a)r^2}\right),
 \end{split}
\end{align}
for $r > \rho/16$.

Next, again using \eqref{eq:deltarbeta}, we compute that
\begin{align*}
 \Delta (r-\rho)^2 &= 2\left(n - \frac{(n-1)\rho}{r}+O\left(\frac{\tau}{r}\right)\right)= O(1)
\end{align*}
for $r > \rho/16$. Then
\begin{align*}
\begin{split}
 (\partial_{\tau} + \Delta)A_1 &= O\left(\frac{r^2}{(\tau+a)^3}\right),\quad \mbox{and}\quad |\nabla A_1| = O\left(\frac{r}{(\tau+ a)^2}\right).
\end{split}
 \end{align*}
Since
\[
\eta_1 = O\left(\frac{\tau^2}{r^2}\right), \quad|\nabla \eta_1| = O\left(\frac{\tau^2}{r^2}\right), \quad \mbox{and}\quad (\partial_{\tau} + \Delta)\eta_1 = O\left(\frac{\tau}{r^2}\right),
\]
 by \eqref{eq:eta1est}, we find that
\begin{align}
\begin{split}\label{eq:a1termev}
  (\partial_{\tau} + \Delta)(A_1 \eta_1)
  &=  \eta_1(\partial_{\tau} + \Delta)A_1
  + A_1(\partial_{\tau} + \Delta)\eta_1
  + 2\langle \nabla A_1, \nabla \eta_1\rangle\\
  &= O\left(\frac{1}{\tau+a}\right).
\end{split}
\end{align}

Similarly, we find that
\begin{equation*}
 (\partial_{\tau} + \Delta)A_2 = O\left(\frac{r}{(\tau+a)^2}\right), \quad \mbox{ and }\quad |\nabla A_2|= O\left(\frac{1}{\tau+a}\right),
\end{equation*}
which, since
\begin{equation*}
\eta_3 = O\left(\frac{\tau}{r^2}\right),
\quad |\nabla \eta_3| = O\left(\frac{\tau}{r^2}\right), \quad \mbox{and} \quad (\partial_{\tau} + \Delta)\eta_3 = O\left(\frac{1}{r^2}\right),
\end{equation*}
implies that
\begin{equation}\label{eq:a2termev}
(\partial_{\tau} + \Delta)(A_2\eta_3) = O\left(\frac{1}{r(\tau+a)}\right).
\end{equation}

Combining \eqref{eq:a0ev} with \eqref{eq:a1termev}, \eqref{eq:a2termev}, and \eqref{eq:bev},
we see that
\begin{align}
\begin{split}\nonumber
  (\partial_{\tau} + \Delta)F_2 &=
\frac{(n-1)\rho}{2r(\tau+a)^2} + O\left(\frac{1}{\tau+a}\right),
\end{split}
\end{align}
and \eqref{eq:f2bheat} follows.
\end{proof}

Returning now to Proposition \ref{prop:pdeintineq2} with the estimates in Propositions \ref{prop:phif2est}
we can now prove the key inequality in the second set of Carleman estimates. This estimate is the analog
of Proposition 5.8 in \cite{KotschwarWangConical}. Note that, while the estimate there is stated in terms of the radial distance function, it is proven via a weight defined in terms of the soliton potential (and in terms of a self-similar shrinking background).

\begin{proposition}\label{prop:pdecarleman2} Suppose $\rad \geq 1$ and $0 < T\leq 1$. There is a constant $C = C(n) > 0$ and constants
$k_0\geq 1$ and $\srad \geq \rad$ depending on $n$ and $K$
such that, for any $\rho \geq 1$ and any section $Z$ of $\mathcal{Z}$ which vanishes on $\Cc_{\srad}\times\{0\}$ and is compactly supported in $(\Cc_{\srad}\cap \Cc_{\rho/16})\times [0, T)$, we have the estimate
\begin{align}
\begin{split}
\label{eq:pdecarleman2}
 &\sqrt{k}\|\sigma_a^{-k-1/2}Z G_{2}^{1/2}\|_{\Cc_{\srad, T}}
  +\|\sigma_a^{-k}\nabla ZG_{2}^{1/2}\|_{\Cc_{\srad, T}}\\
  &\qquad\qquad\leq C\|\sigma_a^{-k}(D_{\tau}Z+\Delta Z)G_{2}^{1/2}\|_{\Cc_{\srad, T}}
\end{split}
\end{align}
for all $a \in (0, 1)$ and $k \geq k_0$.
\end{proposition}
 Here $\sigma_a$ and $G_2 = G_{2; a, \rho}$ are as defined in \eqref{eq:sigmadef} and \eqref{eq:g2def}.
\begin{proof} Fix $0 < a < 1$ and $\rho \geq 1$.
We will specify $k_0$ and $\srad$ over the course of the proof, assuming to begin with that $k_0 \geq 1$
and $\srad \geq \rad$. Below, we will use $C$ to denote a series of positive constants depending at most on the dimension $n$, and $N$ a series of constants depending potentially also on $K$.

We begin by inspecting the right-hand side of the inequality \eqref{eq:pdeintident2}. Using \eqref{eq:genquaddecay}, \eqref{eq:hesslogG2est}, and the identity \eqref{eq:sigmaident} for $\theta = \sigma_a/\dot{\sigma}_a$, we have
\begin{align}\label{eq:q3est}
 \begin{split}
   Q_3(\nabla Z, \nabla Z)  &= \left(2\nabla_i\nabla_j\phi
     +\frac{\dot{\theta}}{\theta} g_{ij}-2\beta_{ij}\right)\langle \nabla_i Z, \nabla_j Z\rangle\\
     &\geq \left(\frac{1}{\tau+a}\left(\frac{1}{1-(\tau+a)/3} - 1\right) - \frac{N}{r}\right)|\nabla Z|^2\\
     &\geq\left(\frac{1}{3} -\frac{N}{r}\right)|\nabla Z|^2\\
     &\geq \frac{7}{24}|\nabla Z|^2,
\end{split}
\end{align}
for all $r \geq \srad$, provided $\srad$ is chosen large enough.

Next note that
\begin{align*}
\begin{split}
  Q_4(Z, Z) &= \frac{1}{2}\left(\partial_{\tau}F_2 + \Delta F_2 +\frac{\dot{\theta}}{\theta}F_2\right)|Z|^2\\
  &\geq \frac{1}{2}\left( \frac{(n-1)\rho}{2r(\tau+a)^2} - \frac{N}{(\tau+a)}\right)|Z|^2\\
  &\phantom{\geq}-
  \frac{1}{2}\left(\frac{1}{1-(\tau+a)/3}\right)\left(\frac{(n-1)\rho}{2r(\tau+a)^2} + \frac{N}{\tau+a}\right)|Z|^2\\
&\geq -\frac{N}{(\tau+a)}|Z|^2,
\end{split}
\end{align*}
so that
\begin{equation}\label{eq:q4est}
Q_4(Z, Z) \geq -N\sigma_{a}^{-1}|Z|^2.
\end{equation}

Finally, we estimate the derivative commutator term $E_{\phi}$. By \eqref{eq:comm1},
we have
\begin{align}
|E_{\phi}(Z, \nabla Z)|
&= C\left(|\nabla \beta| + \left(\frac{|r-\rho|}{\tau+a}\right)|\eta_4|\right)|\nabla Z||Z| + C|\beta||\nabla Z|^2.
\end{align}
where $\eta_4$ is defined by \eqref{eq:eta4def}.
Since we have previously estimated
$|\eta_4| \leq Nr^{-2}\tau$
in \eqref{eq:eta4est}, we have
\begin{equation}\label{eq:eest1}
 |E_{\phi}(Z, \nabla Z)| \leq \frac{N}{r}|Z||\nabla Z| + \frac{N}{r^2}|\nabla Z|^2 \leq N|Z|^2 + \frac{1}{24}|\nabla Z|^2
\end{equation}
for $r\geq S$ if $S$ is sufficiently large.

Putting \eqref{eq:q3est}, \eqref{eq:q4est}, and \eqref{eq:eest1} together, we obtain that
\begin{align}\label{eq:quadtermest}
Q_3(Z, Z) + Q_4(\nabla Z, \nabla Z) + E_{\phi}(Z, \nabla Z)
\geq \frac{1}{4}|\nabla Z|^2 - N\sigma_a^{-1}|Z|^2
\end{align}
and hence, via \eqref{eq:pdeintident2}, that
\begin{align}\label{eq:carleman2est1}
\begin{split}
&\frac{1}{4}\int_{\Cc_{\srad, T}}\,\frac{\sigma_a^{1-\alpha}}{\dot{\sigma_a}}|\nabla Z|^2G\,d\mu\,d\tau\\
 &\quad\quad
 \leq \int_{\Cc_{\srad, T}}
 	\left( \frac{N\sigma_a^{-\alpha}}{\dot{\sigma}}|Z|^2+\frac{\sigma_a^{1-\alpha}}{\dot{\sigma}_a}\left|D_{\tau}Z + \Delta Z\right|^2
 	\right)G\,d\mu\,d\tau.
 \end{split}
\end{align}

Now define $k = (\alpha-1)/2$. Using \eqref{eq:sigmaest} to bound $\dot{\sigma}$ from above and below, we obtain from \eqref{eq:carleman2est1} that, for all $k > 0$,
\begin{align}\label{eq:carleman2est2}
\begin{split}
& 2\|\sigma_a^{-k}\nabla Z G_2^{1/2}\|^2_{\Cc_{\srad, T}}\\
&\qquad \leq N_0\|\sigma_a^{-k-1/2}ZG_2^{1/2}\|^2_{\Cc_{\srad, T}} + N_0\|\sigma_a^{-k}(\partial_{\tau}Z+\Delta Z)G_2^{1/2}\|^2_{\Cc_{\srad, T}},
 \end{split}
\end{align}
for some $N_0 = N_0(n, K)$.

On the other hand, applying Proposition \ref{prop:pdeintineq3} with $\alpha = 2k$, we know
\begin{align}
\begin{split}
  &\int_{\Cc_{\srad, T}}\,\sigma^{-2k}\left(k\frac{\dot{\sigma}_a}{\sigma_a}- F_2\right)|Z|^2G\,d\mu\,d\tau\\
  &\qquad \leq \int_{\Cc_{\srad, T}}
    \,\sigma^{-2k}\left(2|\nabla Z|^2+\frac{\sigma_a}{k\dot{\sigma}_a}\left|D_{\tau}Z + \Delta Z\right|^2
	    \right)G\,d\mu\,d\tau.
\end{split}
\end{align}
By \eqref{eq:f2est}, we have $F_2 \leq N$, hence $F_2\leq \sigma_{a}^{-1}N$,
on $\Cc_{\srad, T}$. As $\dot{\sigma}_a\geq 1/10$ on $[0, 1]$,
we see that if $k_0$ is chosen large enough (depending on $N = N(n, K)$), then, for all $k\geq k_0$,
we have
\begin{align}\label{eq:carleman2est3}
\begin{split}
& k\|\sigma_a^{-k-1/2}ZG_2^{1/2}\|^2_{\Cc_{\srad, T}} \\
&\qquad \leq C_0\|\sigma_a^{-k}\nabla Z G_{2}^{1/2}\|^2_{\Cc_{\srad, T}} + C_0\|\sigma_a^{-k}(D_{\tau}Z + \Delta Z)G_2^{1/2}\|^2_{\Cc_{\srad, T}},
\end{split}
\end{align}
for some constant $C_0 = C_0(n)$.

Combining \eqref{eq:carleman2est2} and  \eqref{eq:carleman2est3}, we obtain
\begin{align}\label{eq:carleman2est4}
\begin{split}
 2\|\sigma_a^{-k}\nabla Z G_2^{1/2}\|^2_{\Cc_{\srad, T}}
 &\leq \frac{C_0N_0}{k}\|\sigma_a^{-k}\nabla Z G_{2}^{1/2}\|^2_{\Cc_{\srad, T}}\\
 &\phantom{\leq}+ N_0\left(1+\frac{C_0}{k}\right)\|\sigma_a^{-k}(D_{\tau}Z+\Delta Z)G_2^{1/2}\|^2_{\Cc_{\srad, T}}
\end{split}
\end{align}
for all $k \geq k_0$. If we choose $k_0$ larger still so that $k_0 \geq C_0 N_0$, we can absorb the first term on the right
into the term on the left side and obtain that
\begin{align}\label{eq:carleman2est5}
\begin{split}
 &\|\sigma_a^{-k}\nabla Z G_2^{1/2}\|^2_{\Cc_{\srad, T}} \leq N_2\|\sigma_a^{-k}(D_{\tau}Z+\Delta Z)G_2^{1/2}\|^2_{\Cc_{\srad, T}}
\end{split}
\end{align}
for all $k \geq k_0$.  Finally, adding \eqref{eq:carleman2est3} to $C_0 +1$ times \eqref{eq:carleman2est5}
and using that $\sqrt{x} + \sqrt{y} \leq \sqrt{2(x+y)}$ for $x$, $y\geq 0$
yields
\begin{align*}
\begin{split}
 &\sqrt{k}\|\sigma_a^{-k-1/2}ZG_2^{1/2}\|_{\Cc_{\srad, T}} + \|\sigma_a^{-k}\nabla Z G_2^{1/2}\|_{\Cc_{\srad, T}}\leq N_3\|\sigma_a^{-k}(D_{\tau}Z+\Delta Z)G_2^{1/2}\|_{\Cc_{\srad, T}},
\end{split}
\end{align*}
for all $k\geq k_0$,
where $N_3 = \sqrt{2(C_0 + (C_0+1)N_2)}$.
\end{proof}

\subsection{The second Carleman inequality for the ODE component}

We will need a matching inequality with the weights $\sigma_a$ and $G_2$ to bound the norm of $D_{\tau}Z$.  The following estimate generalizes (and simplifies the proof of) the estimate in Proposition 5.9 in \cite{KotschwarWangConical}.

\begin{proposition}\label{prop:odecarleman2}  Assume that $\rho$, $\rad \geq 1$ and  $0 < T \leq 1/4$.
There are constants $\gamma =\gamma(n) > 1$, $c= c(n) > 0$, and $N = N(n, K) > 0$ such that, for any section $Z$ of $\Zc$ which vanishes at $\tau= 0$ and has compact support in $\Cc_{\rad}\times [0, T)$,
we have
\begin{align}
\begin{split}\label{eq:odecarleman2}
\|\sigma_{a}^{-k-\frac12}ZG_{2}^{\frac12}\|_{\Cc_{\rad, T}}
&\leq Nk^{-1}\|\sigma_{a}^{-k}D_{\tau}ZG_{2}^{\frac12}\|_{\Cc_{\rad, T}} + N(\gamma a)^{-(k+c)}\rho^{\frac{n-1}{2}}\|Z\|_{\infty}
\end{split}
\end{align}
for all $k > 0$ and all $0 < a < 1/8$.
\end{proposition}

\begin{proof}
Fix $0 < a < 1/8$ and $0 < T \leq 1/4$.
We start with the identity
\begin{align}
\begin{split}
\label{eq:odedivident2}
 &\sigma_a^{-2k}\pdtau|Z|^2G_2
   - \pdtau\left(\sigma_a^{-2k}|Z|^2G_2\right)\\
&\qquad\qquad= \sigma_a^{-2k}\left(\frac{2k+ n/2}{(\tau+a)} - \frac{2k}{3} - \frac{(r-\rho)^2}{4(\tau+a)^2}\right)|Z|^2G_2,
\end{split}
\end{align}
valid for any $k$, $\rho > 0$.
Since
\[
 \pdtau|Z|^2 =2\langle D_{\tau}Z, Z\rangle \leq \frac{3}{k}|D_{\tau}Z|^2 + \frac{k}{3}|Z|^2,
\]
we have
\begin{align}
\begin{split}
\label{eq:odedivineq2}
&\frac{3}{k}\sigma_a^{-2k}\left|D_{\tau}Z\right|^2G_2
   - \pdtau\left(\sigma_a^{-2k}|Z|^2G_2\right)\\
&\qquad\qquad\geq \sigma_a^{-2k}\left(\frac{k}{(\tau+a)} - \frac{(r-\rho)^2}{4(\tau+a)^2}\right)|Z|^2G_2
\end{split}
\end{align}
on $\Cc_{\rad, T}$, and consequently that
\begin{align}
 \label{eq:odeineq3}
 \frac{3}{k}\int_{\Cc_{\rad, T}}\sigma_a^{-2k}|D_{\tau} Z|^2G_2\,d\hat{\mu}\,d\tau \geq
 \int_{\Cc_{\rad, T}}\frac{\sigma_{a}^{-2k}}{\tau+a}\left(k - \frac{(r-\rho)^2}{4(\tau+a)}\right)|Z|^2G_2\,d\hat{\mu}\,d\tau,
\end{align}
where $d\hat{\mu} = d\mu_{\hat{g}}$ is the volume form for the (time-independent) conical metric.

Define the set
\begin{align}
\begin{split}
\label{eq:omegadef}
  \Omega = \Omega_{a, k, \rho} &\dfn \left\{\,(x, \tau)\,|\, (r(x)-\rho)^2 > 2k(\tau+a)\,\right\}.
\end{split}
\end{align}
Then, on $\Omega\cap\operatorname{supp} Z$, we have
\begin{align*}
   \frac{(r-\rho)^2}{\tau+a}\exp\left(-\frac{1}{4}\cdot\frac{(r-\rho)^2}{\tau+a}\right)
   &\leq C\cdot \exp\left(-\frac{7}{32}\cdot\frac{(r-\rho)^2}{\tau+a}\right)\\
   &\leq C\cdot \exp\left(\frac{-(r-\rho)^2}{6} - \frac{5k}{16}\right).
\end{align*}
Here we have used that $xe^{-x/32} \leq C$ and that $\exp(-(r-\rho)^2/(16(\tau+a)))\leq e^{-(r-\rho)^2/6}$ as $\tau + a \leq T + a< 3/8$.
Also, as
\[
\sigma_a(\tau) = (\tau+a)e^{-\frac{\tau+a}{3}} \geq ae^{-\frac{1}{8}},
\]
we have $\sigma_a(\tau)^{-2k} \leq a^{-2k}e^{\frac{k}{4}}$, so
\begin{align}
\begin{split}
\label{eq:odeineq4}
&\frac{1}{4}\int_{\Omega} \sigma_a^{-2k}\frac{(r-\rho)^2}{(\tau+a)^2}|Z|^2G_2\,d\hat{\mu}\,d\tau\\
&\qquad\leq Ce^{-\frac{5k}{16}}\int_{\Cc_{\rad, T}}\frac{\sigma_a^{-2k}}{(\tau+a)^{\frac{n}{2}+1}}|Z|^2e^{-\frac{(r-\rho)^2}{6}}\,d\hat{\mu}\, d\tau\\
&\qquad\leq  Ce^{-\frac{k}{16}} a^{-(2k +\frac{n}{2}+1)}\|Z\|^2_{\infty}\operatorname{vol}_{g_{\Sigma}}(\Sigma)\int_{0}^{\infty}e^{-\frac{(s-\rho)^2}{6}}s^{n-1}\,ds\\
&\qquad\leq N(ae^{\frac{1}{32}})^{-(2k+\frac{n}{2}+1)}\rho^{n-1}\|Z\|_{\infty}^2
\end{split}
\end{align}
for some $N$ depending on $\hat{g}$.

On the other hand, on $\Omega^c$, we have $(r-\rho)^2/(4(\tau+a)^2) \leq k/(2(\tau+a))$, so
\begin{align}
\begin{split}
\label{eq:odeineq5}
\frac{1}{4}\int_{\Omega^c} \sigma^{-2k}\frac{(r-\rho)^2}{(\tau+a)^2}|Z|^2G_2\,d\hat{\mu}\,d\tau
\leq \frac{k}{2} \int_{\Cc_{\rad, T}}\frac{\sigma^{-2k}}{\tau+a}|Z|^2G_2\,d\hat{\mu}\,d\tau.
\end{split}
\end{align}
Thus, taking \eqref{eq:odeineq4} and \eqref{eq:odeineq5} together and returning to \eqref{eq:odeineq3}, we have
\begin{align*}
\begin{split}
 &\frac{3}{k}\int_{\Cc_{\rad, T}}|D_{\tau} Z|^2\sigma_a^{-2k}G\,d\hat{\mu}\,d\tau \\
 &\qquad\geq
 \frac{k}{2}\int_{\Cc_{\rad, T}}\frac{\sigma_{a}^{-2k}}{\tau+a}|Z|^2G_2\,d\hat{\mu}\,d\tau
 - N(ae^{\frac{1}{32}})^{-(2k + \frac{n}{2}+ 1)}\rho^{n-1}\|Z\|^2_{\infty}
\end{split}\\
\begin{split}
 &\qquad\geq
 \frac{k}{6}\int_{\Cc_{\rad, T}}\sigma_{a}^{-2k-1}|Z|^2G_2\,d\hat{\mu}\,d\tau
 - N(ae^{\frac{1}{32}})^{-(2k + \frac{n}{2}+ 1)}\rho^{n-1}\|Z\|^2_{\infty}.
\end{split}
\end{align*}
Using the uniform equivalence of the metrics $g(\tau)$ and $\hat{g}$, we thus obtain the bound
\begin{align*}
\begin{split}
\|\sigma_{a}^{-k-\frac12}ZG_{2}^{\frac12}\|_{\Cc_{\rad, T}}
&\leq Nk^{-2}\|\sigma_{a}^{-k}D_{\tau}ZG_{2}^{\frac12}\|^2_{\Cc_{\rad, T}} +N(ae^{\frac{1}{32}})^{-(2k + \frac{n}{2}+ 1)}\rho^{n-1}\|Z\|^2_{\infty}
\end{split}
\end{align*}
in terms of the space time $L^2$-norms $\|\cdot\|$ induced by $d\mu_{g(\tau)}$. Equation \ref{eq:odecarleman2} then follows
with $\gamma \dfn e^{1/32}$ and $c \dfn n/4 + 1/2$.
\end{proof}

\section{A general backward uniqueness theorem}\label{sec:bu}
In this section, we will continue to assume that $\rad \geq 1$, $0 < T\leq 1$, and that $g = g(\tau)$ is a smooth family of metrics on $\Cc_{\rad, T}$ satisfying \eqref{eq:genmetev} and emanating smoothly from $\gh$ at $\tau = 0$.
We will take $\Xc = \Xc_{\rad, T}$ and $\Yc_{\rad, T}$ to be vector bundles of the form
$\bigoplus \pi^*(T^{k_i}(T^*\Cc_{\rad}))$ over $\Cc_{\rad, T}$ where $\pi:\Cc_{\rad, T}\longrightarrow \Cc_{\rad}$
is the projection as before.

The following theorem is a generalized and  slightly strengthened version of the backward uniqueness principle implicit in the proof of the main theorem in \cite{KotschwarWangConical}.
\begin{theorem}\label{thm:xysysbu}
Let $\Xb$ and $\Yb$ be smooth sections of $\Xc$ and $\Yc$ such that
\begin{equation}\label{eq:xybounds}
  \sup_{\Cc_{\rad, T}}\left\{|\Xb| + |\nabla \Xb| + |\Yb|\right\} \leq L,
\end{equation}
and
\begin{align}\label{eq:xysys2}
\begin{split}
    \left|D_{\tau}\ve{X} + \Delta\ve{X}\right| &\leq \varepsilon\left(|\ve{X}| + |\nabla \ve{X}|+ |\ve{Y}|\right)\\
    \left|D_{\tau}\ve{Y}\right| &\leq L\left(|\ve{X}| + |\nabla \ve{X}| + |\ve{Y}|\right),
\end{split}
\end{align}
for some constant $L > 0$ and some function $\varepsilon = \varepsilon(r) > 0$ with $\varepsilon(r)\longrightarrow 0$ as $r\longrightarrow\infty$. Then, if $\ve{X}(x, 0) \equiv 0$ and $\ve{Y}(x, 0) \equiv 0$ on $\Cc_{\rad}$, there are $\rad_0 = \rad_0(n, K, L, T) \geq \rad$ and $\lambda = \lambda(n) \in  (0, 1)$ such that
$\Xb(x, \tau) \equiv 0$ and $\Yb(x, \tau) \equiv 0$ on $\Cc_{\rad_0, T_0}$ where $T_0 = \lambda T$.
\end{theorem}

\subsection{The rapid decay of $\Xb$ and $\Yb$}
We first show that solutions $\ve{X}$, $\ve{Y}$ which satisfy the assumptions of Theorem \ref{thm:xysysbu}
must at least decay at an exponential-quadratic rate. We follow the outline of Claim 6.2 of \cite{KotschwarWangConical}, taking some care
to track the dependencies of the constant $A$.

\begin{proposition}\label{prop:rapiddecay} Let $g$, $\ve{X}$, and $\ve{Y}$ be as in Theorem \ref{thm:xysysbu}
on $\Cc_{\rad, T}$ where $\rad \geq 1$ and $0 < T\leq 1$. Then there are dimensional constants $0 < A \leq 1$
and $B > 0$, and further constants $N > 0$ and $\srad\geq \rad$ depending on $n$, $K$, and $L$ such that
\begin{equation}\label{eq:rapiddecay}
  \big\||\ve{X}| + |\nabla\ve{X}| + |\ve{Y}| \big\|_{\Dc_{\rho}(\sqrt{a})}
  \leq N e^{-\frac{B\rho^2}{a}}
\end{equation}
for all $\rho > 16\srad$ and $0 < a \leq A T$. Here $\Dc_{r}(s)$ denotes the annular parabolic cylinder
\[
  \Dc_{\rho}(s) = (\rho- s,\rho+s) \times\Sigma \times [0, s^2]
\]
\end{proposition}
\begin{proof}
We will specify the constants $A$, $B$, $N$, and $\srad$ over the course of the proof.
Assume initially that $\srad \geq 2\rad$ and fix an arbitrary $\rho > 16\srad$. As in previous arguments, we will use $C$ to denote a sequence of constants depending only on $n$
 and use $N$ to denote a sequence of constants which may depend in addition on the parameters $K$ and $L$.

We will need two cutoff functions, one temporal and one spatial. For the temporal cutoff, we select $\varphi\in C^{\infty}(\RR)$ satisfying
\[
  \varphi(\tau) = \left\{\begin{array}{rl}
                          1 & \tau \leq T/6,\\
                          0 & \tau \geq T/5,
                         \end{array}\right.  \quad \mbox{and} \quad  -CT^{-1} \leq \varphi^{\prime}(\tau) \leq 0.
\]
For the spatial cutoff, we select for each $\xi > 4\rho$ a function $\psi_{\rho, \xi}$ in the form
\[
\psi_{\rho, \xi}(x) = \Psi_{\rho, \xi}(r(x)),
\]
where
$\Psi_{\rho, \xi} \in C^{\infty}(\RR, [0, 1])$ satisfies
\[
  \Psi_{\rho, \xi}(r) = \left\{\begin{array}{rl}
                                1 & \rho/3 \leq r \leq 2\xi,\\
                                0 &  r < \rho/6 \,\, \mbox{or} \,\, r > 3\xi,
                               \end{array}\right. \quad \mbox{and}  \quad \rho|\Psi_{\rho, \xi}^{\prime}| + \rho^2|\Psi_{\rho, \xi}^{\prime\prime}| \leq N.
\]
(For the latter estimate, we have used that $\xi > 4\rho$ and $\rho > 16\srad$.) In view of Proposition \ref{prop:genetaests}, we will have the bounds
\begin{align}\label{eq:psiderivest}
\begin{split}
  |\nabla \psi_{\rho, \xi}| &\leq |\Psi_{\rho, \xi}^{\prime}|(1+ \eta_1)^{1/2} \leq N\rho^{-1},\\
  |\Delta \psi_{\rho, \xi}| &\leq |\Psi_{\rho, \xi}^{\prime}||(n-1)r^{-1} + \eta_3| + |\Psi_{\rho, \xi}^{\prime\prime}|(1+ |\eta_1|)
  \leq N\rho^{-2},
\end{split}
\end{align}
on the derivatives of $\psi_{\rho, \xi}$ by Proposition \ref{prop:genetaests}.

The sections
\[
   \ve{X}_{\rho, \xi} = \varphi \psi_{\rho, \xi}\ve{X} \quad \mbox{and}\quad \ve{Y}_{\rho, \xi} = \varphi \psi_{\rho, \xi} \ve{Y}
\]
will then have compact support in the annular cylinder $\Ac(\rho/6, 3\xi) \times [0, T/4)$ where
\[
  \Ac(r_1, r_2) = \{\,(r, \sigma)\in \Cc_0\,|\, r_1 < r < r_2\,\},
\]
and will vanish identically for $\tau=0$.

Applying Proposition \ref{prop:pdecarleman2} to the components
of $\ve{X}_{\rho, \xi}$, it follows that if $\srad$ is sufficiently
large, we have
\begin{align}\label{eq:xcarleman2}
\begin{split}
 &k^{\frac12}\|\sigma_a^{-k-\frac{1}{2}}\Xb_{\rho, \xi}G_2^{\frac{1}{2}}|\|_{\Ac(\frac{\rho}{6},3\xi)\times [0, \frac{T}{5}]}
  +\|\sigma_a^{-k}\nabla \Xb_{\rho, \xi}G^{\frac12}_2\|_{\Ac(\frac{\rho}{6},3\xi)\times [0, \frac{T}{5}]}\\
   &\leq N\|\sigma_a^{-k}(D_{\tau}+\Delta) \Xb_{\rho, \xi}G_2^{\frac12}\|_{\Ac(\frac{\rho}{6},3\xi)\times [0, \frac{T}{5}]}
\end{split}
\end{align}
for all $0 < a \leq 1$, $\rho > 16\srad$, and $k \geq k_0= k_0(n, K)$. Here and below, we use $\|\cdot\|$ to denote the either the space-time $L^2$-norm or spatial $L^2$-norm induced by the family of metrics $g(\tau)$.

On the other hand, by \eqref{eq:xysys2} and \eqref{eq:psiderivest}, we have
 \begin{align*}
  |D_{\tau}\Xb_{\xi, \rho} + \Delta \Xb_{\rho, \xi}| &\leq \varepsilon(|\Xb_{\rho, \xi}| + |\nabla \Xb_{\rho, \xi}|
  + |\Yb_{\rho, \xi}|) + \psi_{\rho, \xi}|\varphi^{\prime}||\Xb|\\
  &\phantom{\leq}+\varphi\left((|\Delta \psi_{\rho, \xi}| + \varepsilon |\nabla\psi_{\rho, \xi}|)|\Xb| + 2 |\nabla \psi_{\rho, \xi}||\nabla \ve{X}|\right)\\
  &\leq \varepsilon(|\Xb_{\rho, \xi}| + |\nabla \Xb_{\rho, \xi}|
  + |\Yb_{\rho, \xi}|) + CT^{-1}|\Xb|\chi_{\supp \varphi^{\prime}}\\
  &\phantom{\leq}+ N\rho^{-1}\left(|\Xb|+ |\nabla \Xb|\right)\chi_{\supp \nabla\psi_{\rho, \xi}}.
 \end{align*}
 Integrating this inequality over space and time, and combining it with \eqref{eq:xcarleman2}, we see that
  if $\srad$ is taken large enough, then
\begin{align}\label{eq:xcarleman2b}
\begin{split}
  &k^{\frac{1}{2}}\|\sigma_a^{-k-\frac{1}{2}}\Xb_{\rho, \xi}G_2^{\frac{1}{2}}||_{\Ac(\frac{\rho}{6},3\xi)\times [0, \frac{T}{5}]}
  +\|\sigma_a^{-k}\nabla \Xb_{\rho, \xi}G^{\frac{1}{2}}_2\|_{\Ac(\frac{\rho}{6},3\xi)\times [0, \frac{T}{5}]}\\
   &\qquad\leq \frac{1}{2}\|\sigma_a^{-k-\frac{1}{2}}\Yb_{\rho, \xi}G_2^{\frac{1}{2}}||_{\Ac(\frac{\rho}{6},3\xi)\times [0, \frac{T}{5}]}
  +CT^{-1} \|\sigma_a^{-k-\frac{1}{2}}\Xb G_2^{\frac{1}{2}}\|_{\Ac(\frac{\rho}{6},3\xi)\times [\frac{T}{6}, \frac{T}{5}]}\\
  &\qquad\phantom{\leq} + N\|\sigma_a^{-k}(|\Xb| +|\nabla \Xb|)G_2^{\frac{1}{2}}\|_{\Ac(\frac{\rho}{6},\frac{\rho}{3})\times [0, \frac{T}{5}]}\\
  &\qquad\phantom{\leq}
  + N\|\sigma_a^{-k}(|\Xb| +|\nabla \Xb|)G_2^{\frac{1}{2}}\|_{\Ac(2\xi, 3\xi)\times [0, \frac{T}{5}]},
\end{split}
\end{align}
using that $\varepsilon \longrightarrow 0$ as $\rho\longrightarrow\infty$, and that $\rho > 16\srad$.

On the other hand,  applying Proposition \ref{prop:odecarleman2} to $\ve{Y}_{\rho, \xi}$, we obtain
\begin{align}\label{eq:ycarleman2}
\begin{split}
\|\sigma_{a}^{-k-\frac12}\Yb_{\rho, \xi}G_{2}^{\frac12}\|_{\Ac(\frac{\rho}{6}, 3\xi)\times [0, \frac{T}{5}]}
&\leq Nk^{-1}\|\sigma_{a}^{-k}D_{\tau}\Yb_{\rho, \xi}G_{2}^{\frac12}\|_{\Ac(\frac{\rho}{6},3\xi)\times [0, \frac{T}{5}]}\\
&\phantom{\leq}+ N(\gamma a)^{-(k+c)}\rho^{\frac{n-1}{2}}\|\Yb\|_{\infty},
\end{split}
\end{align}
for some $\gamma > 1$ and $c > 0$ and all $0 < a < 1/8$ and $k > 0$.
Since
\begin{align*}
\begin{split}
\left|D_{\tau}\Yb_{\rho, \xi}\right| &\leq
  L\left(|\Xb_{\rho, \xi}| + |\nabla \Xb_{\rho, \xi}| + |\Yb_{\rho, \xi}|\right) + L\varphi|\nabla\psi_{\rho, \xi}||\Xb|+\psi_{\rho, \xi}|\varphi^{\prime}||\Yb|
\end{split}\\
\begin{split}
       &\leq N\left(|\Xb_{\rho, \xi}| + |\nabla \Xb_{\rho, \xi}|+ |\Yb_{\rho, \xi}|\right)
       + N\rho^{-1}|\Xb|\chi_{\supp \nabla\psi_{\rho, \xi}}\\
  &\phantom{\leq}
      +CT^{-1}\psi_{\rho, \xi}|\Yb|\chi_{\supp \varphi^{\prime}},
\end{split}
\end{align*}
returning to \eqref{eq:ycarleman2}, we see that, for  suitably large $k_1 = k_1(n, K, L) \geq k_0$
we have
\begin{align}\label{eq:ycarleman2b}
\begin{split}
&\|\sigma_{a}^{-k-\frac12}\Yb_{\rho, \xi}G_{2}^{\frac12}\|_{\Ac(\frac{\rho}{6}, 3\xi)\times [0, \frac{T}{5}]}\\
&\qquad
\leq \frac{1}{2}\left(\|\sigma_a^{-k-\frac{1}{2}}\Xb_{\rho, \xi}G_2^{\frac{1}{2}}||_{\Ac(\frac{\rho}{6},3\xi)\times [0, \frac{T}{5}]}
  +\|\sigma_a^{-k}\nabla \Xb_{\rho, \xi}G^{\frac12}_2\|_{\Ac(\frac{\rho}{6},3\xi)\times [0, \frac{T}{5}]}\right)\\
  &\qquad\phantom{\leq} + N\|\sigma_a^{-k}\Xb G_2^{\frac12}\|_{\Ac(\frac{\rho}{6},\frac{\rho}{3})\times [0, \frac{T}{5}]}
  + N\|\sigma_a^{-k}\Xb G_2^{\frac12}\|_{\Ac(2\xi, 3\xi)\times [0, \frac{T}{5}]}\\
  &\qquad\phantom{\leq} +CT^{-1} \|\sigma_a^{-k-\frac{1}{2}}\Yb G_2^{\frac{1}{2}}|\|_{\Ac(\frac{\rho}{6},3\xi)\times [\frac{T}{6}, \frac{T}{5}]} + N(\gamma a)^{-(k+c)}\rho^{\frac{n-1}{2}}
\end{split}
\end{align}
for all $k \geq k_1$.

Now we add \eqref{eq:ycarleman2b} to \eqref{eq:xcarleman2b} and use again \eqref{eq:ycarleman2b} to estimate the first term on the right of \eqref{eq:xcarleman2b}. The coefficients on of the terms on the right of this new inequality involving $L^2$-norms of $\Xb_{\rho, \xi}$, $\nabla \Xb_{\rho, \xi}$, and $\Yb_{\rho, \xi}$
over $\Ac(\rho/6, 3\xi)\times [0, \frac{T}{5}]$ are small enough that they may be absorbed into
their counterparts on the left. The result is the inequality
\begin{align}\label{eq:xycarleman2a}
\begin{split}
&\|\sigma_a^{-k-\frac{1}{2}}(|\Xb_{\rho, \xi}|+|\Yb_{\rho, \xi}|)G_2^{\frac{1}{2}}|\|_{\Ac(\frac{\rho}{6},3\xi)\times [0, \frac{T}{5}]}
  +\|\sigma_a^{-k}\nabla \Xb_{\rho, \xi}G_2^{\frac12}\|_{\Ac(\frac{\rho}{6},3\xi)\times [0, \frac{T}{5}]}\\
  &\qquad\leq CT^{-1} \|\sigma_a^{-k-\frac{1}{2}}(|\Xb|+|\Yb|)G_2^{\frac{1}{2}}|\|_{\Ac(\frac{\rho}{6},3\xi)\times [\frac{T}{6}, \frac{T}{5}]}\\
  &\qquad\phantom{\leq} + N\|\sigma_a^{-k}(|\Xb| +|\nabla \Xb| + |\Yb|)G_2^{\frac12}\|_{\Ac(\frac{\rho}{6},\frac{\rho}{3})\times [0, \frac{T}{5}]}\\
  &\qquad \phantom{\leq} + N\|\sigma_a^{-k}(|\Xb| +|\nabla \Xb| + |\Yb|)G_2^{\frac12}\|_{\Ac(2\xi, 3\xi)\times [0, \frac{T}{5}]} + N\rho^{\frac{n-1}{2}}(\gamma a)^{-(k+c)},
\end{split}
\end{align}
which is valid for all
for all $0 < a < 1/8$, $k \geq k_1$, and $\rho$ and $\xi$ satisfying $\xi \geq 4\rho$ and $\rho \geq 16 \srad$.

Since the metrics $g(\tau)$ are uniformly equivalent to the cone metric, we have
\[
     \operatorname{vol}_{g(\tau)}(\Ac(s_0, s)) \leq Ns^n
\]
for all $\rad < s_0< s$.
Thus, given \eqref{eq:xybounds} and the exponential-quadratic decay of $G_2$, if we send $\xi\longrightarrow \infty$, the limits of all terms in the inequality \eqref{eq:xycarleman2a}
exist and are finite, and the penultimate term in the last line vanishes. Shrinking the domain of the left side of the resulting inequality, and using \eqref{eq:sigmaest} to compare $\sigma_a(\tau)$ with $\tau+a$, we obtain
\begin{align}\label{eq:xycarleman2b}
\begin{split}
&\|(\tau+a)^{-k}(|\Xb|+|\nabla \Xb| + |\Yb|)G_2^{\frac{1}{2}}|\|_{\Cc_{\frac{\rho}{3}, \frac{T}{6}}}\\
  &\qquad\leq C^k T^{-1} \|(\tau+a)^{-k}(|\Xb|+|\Yb|)G_2^{\frac{1}{2}}|\|_{\Cc_{\frac{\rho}{6}}\times [\frac{T}{6}, \frac{T}{5}]}\\
  &\qquad\phantom{\leq} + C^kN\|(\tau+a)^{-k}(|\Xb| +|\nabla \Xb| + |\Yb|)G_2^{\frac12}\|_{\Ac(\frac{\rho}{6},\frac{\rho}{3})\times [0, \frac{T}{5}]}\\
  &\qquad\phantom{\leq} + N(\gamma a)^{-(k+c)}\rho^{\frac{n-1}{2}}\\
  &\qquad\dfn \mathrm{I} + \mathrm{II} + \mathrm{III}.
\end{split}
\end{align}

Consider the term $\mathrm{I}$. As $0 < a < 1/8$ and $0 < T\leq 1$, we have $1/6 < \tau + a < 1/3$ and
\begin{equation*}
   G_2(x, \tau) \leq Ce^{-\frac{3(r-\rho)^2}{4}}
\end{equation*}
on the domain $\Cc_{\rho/6}\times[T/6, T/5]$ of the integral.
Thus we can estimate
\begin{align}\label{eq:Iest}
 \begin{split}
    \mathrm{I}^2 &= C^k T^{-1} \|(\tau+a)^{-k}(|\Xb|+|\Yb|)G_2^{\frac{1}{2}}|\|^2_{\Cc_{\frac{\rho}{6}}\times [\frac{T}{5}, \frac{T}{6}]}\\
    &\leq C^kNT^{-1}\int_{\frac{T}{6}}^{\frac{T}{5}}\int_{\frac{\rho}{6}}^{\infty}e^{-\frac{3(s-\rho)^2}{4}}s^{n-1}\,ds \leq C^kN\rho^{n-1}.
 \end{split}
\end{align}

Now consider the term $\mathrm{II}$. On its domain of integration, we have
\[
  G^{\frac12}_2(x, \tau) \leq (\tau + a)^{-\frac{n}{4}} e^{-\frac{\rho^2}{18(\tau+a)}}.
\]
Now, Stirling's inequality implies that, for any integer $m > 0$,
\[
 \max_{s > 0} s^{-m} e^{-\frac{\rho^2}{18s}}
 = \rho^{-2m}(18m)^m e^{-m} \leq \rho^{-2m} C^m m!,
\]
so, taking $m = k +n/4$, we see that
\[
  (\tau+a)^{-k} G_2^{\frac12} \leq C^k\rho^{-2k-\frac{n}{2}}k!
\]
on that domain. Hence
\begin{align}\label{eq:IIest}
\begin{split}
 \mathrm{II} &=  C^kN\|(\tau+a)^{-k}(|\Xb| +|\nabla \Xb| + |\Yb|)G_2^{\frac12}\|_{\Ac(\frac{\rho}{6},\frac{\rho}{3})\times [0, \frac{T}{5}]}\\
 &\leq C^kN\rho^{-2k-\frac{n}{2}}k!\left(\int_0^{\frac{T}{5}}\int_{\Ac(\frac{\rho}{6}, \frac{\rho}{3})}\,d\mu_{g(\tau)}\,d\tau\right)^{\frac12}\\
 &\leq C^kN\rho^{-2k}k!.
\end{split}
\end{align}

Combining \eqref{eq:xycarleman2b} with \eqref{eq:Iest} and \eqref{eq:IIest},
we obtain that
\begin{align*}
  \begin{split}
     &\|(\tau+a)^{-k}(|\Xb|+|\nabla \Xb| + |\Yb|)G_2^{\frac{1}{2}}\|_{\Cc_{\frac{\rho}{3}, \frac{T}{6}}}\\
     &\qquad \leq N\rho^{\frac{n}{2}}\left(C^k(1 + \rho^{-2k}k!) + (\gamma a)^{-(k+c)}\right)
  \end{split}
 \end{align*}
 for all $k \geq k_1$, that is,
\begin{align}\label{eq:xycarleman2d}
 \begin{split}
    &\|(\tau+a)^{-(l + k_1)}(|\Xb|+|\nabla \Xb| + |\Yb|)G_2^{\frac{1}{2}}\|_{\Cc_{\frac{\rho}{3}}\times [0, \frac{T}{6}]}\\
    &\qquad \leq N\rho^{\frac{n}{2}}\left(C_1^{l}\left(1 + \rho^{-2l}(l+ k_1)!\right) + (\gamma a)^{-(l + k_1+c)}\right)
 \end{split}
\end{align}
for all $l \geq 0$ where $C_1$ is some universal constant. Now,
\[
(l + k_1)! \leq C^{k_1+l}(k_1!)(l!)\leq NC_2^l(l!)
\]
for all $l\geq 0$ for some other universal constant $C_2$. Thus, if we multiply both sides of \eqref{eq:xycarleman2d}
by $\rho^{2l}/(C^{l}l!)$ for some $C \geq 2\times\max\{C_1, C_2\}$  (and also discard the factors $(\tau+a)^{-(k_1 +n/4)} \geq 1$ on the left),
we obtain that
\begin{align*}
 \begin{split}
    &\left\|\frac{\rho^{2l}}{(C(\tau+a))^{l}l!}(|\Xb|+|\nabla \Xb| + |\Yb|) e^{-\frac{(r-\rho)^2}{8(\tau+a)}}\right\|_{\Cc_{\frac{\rho}{3}. \frac{T}{6}}}\\
    &\qquad \leq N\rho^{\frac{n}{2}}\left(\frac{\rho^{2l}}{2^ll!}+ \frac{1}{2^l} + a^{-(k_1+c)}\frac{\rho^{2l}}{(C\gamma a)^ll!}\right)
 \end{split}
\end{align*}
for all $l\geq 0$. Summing over $l$ then implies
\begin{align}\label{eq:xycarleman2f}
 \begin{split}
    &\left\|(|\Xb|+|\nabla \Xb| + |\Yb|)e^{\frac{\rho^2}{C(\tau+a)} - \frac{(r-\rho)^2}{8(\tau+a)}}\right\|_{\Cc_{\frac{\rho}{3}, \frac{T}{6}}}\\
    &\qquad \leq N\rho^{\frac{n}{2}}\left(1 + e^{\frac{\rho^2}{2}} + a^{-(k_1+c)}e^{\frac{\rho^2}{C\gamma a}}\right).
 \end{split}
\end{align}

Now we assume at least that $0 < A< 1/8$ and restrict to $0 < a \leq AT$.
Recall that the constant $\gamma$ from Proposition \ref{prop:odecarleman2} is given by $\gamma = e^{\frac{1}{32}}\in (1, 1.04)$.
Define $\delta > 0$ by $\gamma = 1+2\delta$. Then, if $\tau\in [0, \delta a]$ and $0 < a \leq AT$ we have $\tau \in [0, T/6]$
and $e^{\frac{\rho^2}{C(1+\delta)a}} \leq e^{\frac{\rho^2}{C(\tau+a)}}$. Multiplying
both sides of \eqref{eq:xycarleman2f} by $e^{-\frac{\rho^2}{C(1+\delta)a}}$, we find that
\begin{align}\label{eq:xycarleman2g}
 \begin{split}
    &\left\|(|\Xb|+|\nabla \Xb| + |\Yb|)e^{-\frac{(r-\rho)^2}{8(\tau+a)}}\right\|_{\Cc_{\frac{\rho}{3}}\times [0, \delta a]}\\
    &\qquad \leq N\rho^{\frac{n}{2}}\left(e^{\frac{-\rho^2}{C(1+\delta)a}}(1 + e^{\frac{\rho^2}{2}})
    + a^{-(k_1+c)}e^{\frac{-\delta\rho^2}{C(1+2\delta)(1+\delta) a}}\right).
 \end{split}
\end{align}
Now (given the dependencies of $k_1$) we can bound $a^{-(k_1+c)}e^{\frac{-\delta\rho^2}{2C(1+2\delta)(1+\delta) a}} \leq N$ for $N$ depending on $n$, $K$, and $L$ but not on $a$. Also, taking $A= A(n)$ smaller still if necessary, the expression
$e^{\frac{-\rho^2}{C(1+\delta)a}}(1 + e^{\frac{\rho^2}{2}})$ is bounded above independently of $\rho$ for all $0 < a \leq AT \leq A$.
Thus we obtain that
\begin{align}\label{eq:xycarleman2h}
 \begin{split}
    \left\|(|\Xb|+|\nabla \Xb| + |\Yb|)e^{-\frac{(r-\rho)^2}{8(\tau+a)}}\right\|_{\Cc_{\frac{\rho}{3}}\times [0, \delta a]}\leq Ne^{-\frac{\rho^2}{C_3a}},
 \end{split}
 \end{align}
for all $0 < a \leq AT$ and some $C_3 = C_3(n) > 0$.

On the other hand, $e^{-\frac{(r-\rho)^2}{8(\tau+a)}} \geq e^{-\frac{\delta}{8}}$ on
$\Ac(\rho - \sqrt{\delta a}, \rho + \sqrt{\delta a})\times [0, \delta a]$,
and given that $\rho  > 16 \srad > 16$, we have $\rho - \sqrt{\delta a} > \rho/3$. Thus,
shrinking the domain of integration in the term on the left of \eqref{eq:xycarleman2h} we find that
\begin{align}\label{eq:xycarleman2i}
 \begin{split}
    \left\||\Xb|+|\nabla \Xb| + |\Yb|\right\|_{\Dc_{\rho}(\sqrt{\delta a})}\leq Ne^{\frac{-\rho^2}{C_3a}}
 \end{split}
 \end{align}
 for all $0 < a \leq AT$.  Rewriting $\delta a$ as $a$ and $\delta A$ as $A$, the claim follows with $B = \delta/C_3$.
\end{proof}

\subsection{Proof of backward uniqueness}
With the decay estimates from Proposition \ref{prop:rapiddecay} in hand, we can now apply the Carleman estimates in Proposition \ref{prop:carlemanest1} to complete
the proof of Theorem \ref{thm:xysysbu}. We modify the argument of the corresponding statement (Claim 6.3) in \cite{KotschwarWangConical} so as to demonstrate that the length of the time-interval on which the sections are guaranteed to vanish is independent of the solutions.
(In \cite{KotschwarWangConical}, the length of this interval is non-explicit, and arises from an invocation of the mean-value theorem.)

\begin{proof}[Proof of Theorem \ref{thm:xysysbu}]
We will again use $C$ to denote a series of constants depending only on $n$, and $N$ to denote a series of  constants which may depend in addition on the parameters $K$ and $L$. Let $A_0$, $B_0$, and $S_0$ be those values
specified in Proposition \ref{prop:rapiddecay}. Define
\[
\lambda = \frac{1}{2}\min\{A_0, B_0/4\}
\]
and assume initially that $\srad \geq 1 + 16\srad_0$.
We will
continue to increase $\srad$ over the course of the proof.

From \eqref{eq:rapiddecay}, it follows from our choice of $\lambda$ that
then
\begin{equation}\label{eq:basicdecayest}
    \||\ve{X}|+ |\nabla \ve{X}| + |\Yb|\|_{\Dc_s(\sqrt{2\lambda T})} \leq N e^{-4s^2}
\end{equation}
for all $s \geq \srad$.  Expressing $\Cc_{\srad_0, 2\lambda T}$ as the union
$\Cc_{\srad, 2\lambda T} = \cup_{i=1}^{\infty}\Dc_{s_i}(\sqrt{2\lambda T})$, where the $s_i$ are suitably close to each other and tend to $+\infty$,
applying \eqref{eq:basicdecayest} to each set in the union, and summing the result, it follows
that
\begin{equation}\label{eq:growthestimate1}
   \|(|\ve{X}|+ |\nabla \ve{X}| + |\Yb|)e^{2r^2}\|_{\Cc_{\srad, 2\lambda T}} \leq N_0,
\end{equation}
for some $N_0 = N_0(n, K, L)$.

The inequality \eqref{eq:growthestimate1} implies, in particular, that, for each $m \geq 1$, there is $\tau_m$ satisfying
\[
 2\lambda T\left(1-\frac{1}{m}\right) < \tau_m < 2\lambda T,
\]
and
\begin{equation}\label{eq:timeslicedecay}
  \int_{\Cc_{\srad}\times\{\tau_m\}}(|\Xb| + |\nabla \Xb| + |\Yb|)^2\,e^{4r^2}\,d\mu \leq N_m,
\end{equation}
where $N_m$ depends on $m$ as well as $n$, $K$, and $L$ (in fact, $N_m$ will have
the form $N_m = N\cdot m$). Indeed, given \eqref{eq:growthestimate1}, the set
\[
  \left\{\,\tau\in [0, 2\lambda T]\,\middle|\, \left\|\left(|\Xb| + |\nabla\Xb|
  + |\Yb|\right)e^{2r^2}\right\|^2_{\Cc_{\srad}\times\{\tau\}} > m N_0^2/(\lambda T)\right\}
\]
must have measure less than $\lambda T/m$ by Chebychev's theorem.

Select any such sequence $\{\tau_m\}_{m=1}^{\infty}$ of times. We will show below that
$\ve{X}$ and $\ve{Y}$ vanish on $\Cc_{\srad}\times [0, \tau_m/2]$, for some $\srad$ independent of $m$. Sending
$m\longrightarrow \infty$, we will have $\tau_m/2\longrightarrow \lambda T$
and Theorem \ref{thm:xysysbu} will follow.

Let us fix some element $\tau_m$ of the sequence described above. Define
\[
   G_1(x, \tau) \dfn G_{1; \alpha, \tau_m}(x, \tau) = \exp\left(\alpha(\tau_m - \tau)r^{2-\delta}(x) + r^2(x)\right)
\]
on $\Cc_{\rad_0, \tau_m}$ for $\alpha \geq 0$. Note that we have $G_1(x, \tau_m) = e^{r^2}$ on $\Cc_{\rad_0}$,
and, in general, for any $\theta > 1$ we have $G_1 \leq e^{\theta r^{2}}$ for all $r$ sufficiently large.

Recalling \eqref{eq:xysys2}, we may choose $\srad$ so large (depending on $L$) such that
\begin{align}\label{eq:xysys3}
\begin{split}
  |D_{\tau}\ve{X} + \Delta\ve{X}| &\leq \frac{1}{100}(|\Xb| + |\nabla \Xb| + |\Yb|),\\
       |D_{\tau} \ve{Y}| &\leq N(|\ve{X}| + |\nabla \ve{X}| + |\ve{Y}|),
\end{split}
\end{align}
on $\Cc_{\srad, T}$.

Let $\rho > \srad$  and $\xi > 4\rho$. As in the proof of Proposition \ref{prop:rapiddecay}, we can construct a spatial cutoff function $\psi_{\rho, \xi}:\Cc_{\rad}\longrightarrow [0, 1]$
with
\[
    \psi_{\rho, \xi}(x)  = \left\{\begin{array}{rl}
                             1 &  2\rho < r(x) < \xi,\\
                             0 &  r(x) < \rho \ \mbox{ or } \ r(x) > 2\xi,
                         \end{array}\right.
\]
and
\[
|\nabla \psi_{\rho, \xi}| + |\Delta \psi_{\rho, \xi}| \leq N,
\]
on $\Cc_{\rad, T}$.  With these choices, the sections
\[
   \Xb_{\rho, \xi} \dfn \psi_{\rho, \xi} \Xb \quad \mbox{and} \quad \Yb_{\rho, \xi} \dfn \psi_{\rho, \xi}\Yb
\]
will be compactly supported in $\Cc_{\srad, \tau_m}$ and vanish identically at $\tau=0$.

Now we apply the estimates \eqref{eq:pdecarlemanest1} and \eqref{eq:odecarlemanest1}
of Proposition \ref{prop:carlemanest1} to $\Xb_{\rho, \xi}$ and $\Yb_{\rho, \xi}$, respectively, and sum
the resulting inequalities to obtain that
\begin{align}
\label{eq:xycarleman1a}
\begin{split}
 &\left\|\left(\sqrt{\alpha}|\Xb_{\rho, \xi}| + |\nabla \Xb_{\rho, \xi}| + |\Yb_{\rho, \xi}|\right)G_1^{\frac12}\right\|_{\Cc_{\srad, \tau_m}}
  \\
  &\qquad\le C\|(D_{\tau} +\Delta) \Xb_{\rho, \xi}G_1^{\frac12}\|_{\Cc_{\srad, \tau_m}} +
  \frac{C}{\sqrt{\alpha}}\|D_{\tau}\Yb_{\rho, \xi} G_1^{\frac12}\|_{\Cc_{\srad, \tau_m}}\\
  &\qquad\phantom{\le}
 + C\|\nabla \Xb_{\rho, \xi} G_1^{\frac12}\|_{\Cc_{\srad}\times \{\tau_m\}}
\end{split}
\end{align}
for all $\alpha \geq 1$.  Since \eqref{eq:xysys3} implies that
\begin{align*}
  \left|\left(D_{\tau} + \Delta\right)\Xb_{\rho, \xi}\right|
    &\leq \frac{1}{100} \left(|\Xb_{\rho, \xi}| + |\nabla \Xb_{\rho, \xi}|+ |\Yb_{\rho, \xi}|\right) + 2|\nabla\psi_{\rho, \xi}||\nabla\Xb|\\
      &\phantom{\leq} + \left(|\Delta \psi_{\rho, \xi}| + \frac{1}{100}|\nabla \psi_{\rho, \xi}|\right)|\Xb|,
\end{align*}
and
\begin{align*}
  \left|D_{\tau}\Yb_{\rho, \xi}\right| &\leq N(|\Xb_{\rho, \xi}| + |\nabla\Xb_{\rho, \xi}|+|\Yb_{\rho, \xi}|)
	  +N|\nabla\psi_{\rho, \xi}| |\Xb|,
\end{align*}
on $\Cc_{\srad, \tau_m}$,
using that
\[
  |\nabla \ve{X}_{\rho, \xi}| \leq |\nabla\psi_{\rho, \xi}||\ve{X}| + \psi_{\rho, \xi}|\nabla \ve{X}|,
\]
it follows that there exists $\alpha_0 = \alpha_0(n, K)$ such that
\begin{align}\label{eq:xycarleman1b}
\begin{split}
&\|(|\Xb| + |\Yb|) G_1^{\frac12}\|_{\Ac(2\rho, \xi)\times[0, \tau_m]}\\
  &\qquad\le N\|(|\Xb| + |\nabla \Xb|)G_1^{\frac12}\|_{\Ac(\rho, 2\rho)\times[0, \tau_m]}\\
  &\qquad\phantom{\le} + N\|(|\Xb| + |\nabla \Xb|)G_1^{\frac12}\|_{\Ac(\xi, 2\xi)\times[0, \tau_m]}\\
   &\qquad\phantom{\le} + N\|\Xb e^{\frac{r^2}{2}}\|_{\Ac(\rho, 2\rho)\times\{\tau_m\}}
     + N\|\Xb e^{\frac{r^2}{2}}\|_{\Ac(\xi, 2\xi)\times\{\tau_m\}}\\
   &\qquad\phantom{\le} +N\|\nabla \Xb e^{\frac{r^2}{2}}\|_{\Ac(\rho,2\xi)\times\{\tau_m\}}
\end{split}
\end{align}
for all $\alpha \geq \alpha_0$. Given the bounds \eqref{eq:growthestimate1} and \eqref{eq:timeslicedecay},
all of the integrals in \eqref{eq:xycarleman1b} tend to finite limits as $\xi\longrightarrow\infty$,
and the limits of the second and fourth terms on the left vanish. Thus, sending $\xi\longrightarrow\infty$,
we obtain that
\begin{align}
\begin{split}\nonumber
\|(|\Xb| + |\Yb|) G_1^{\frac12}\|_{\Cc_{4\rho, \frac{\tau_m}{2}}}
&\le N\|(|\Xb| + |\nabla \Xb|)G_1^{\frac12}\|_{\Ac(\rho, 2\rho)\times[0, \tau_m]}\\
&\phantom{\leq} + N\|\Xb e^{\frac{r^2}{2}}\|_{\Ac(\rho, 2\rho)\times\{\tau_m\}} +N\|\nabla \Xb e^{\frac{r^2}{2}}\|_{\Cc_{\rho}\times\{\tau_m\}}
\end{split}\\
\begin{split}\label{eq:xycarleman1c}
 &\le N\|(|\Xb| + |\nabla \Xb|)G_1^{\frac12}\|_{\Ac(\rho, 2\rho)\times[0, \tau_m]} + N,
\end{split}
\end{align}
after shrinking the domain of integration on the right.

Now, on $\Ac(\rho, 2\rho)\times [0, \tau_m]$, we have that
\[
  G_{1}(x, \tau) \leq \exp\left(4\rho^2(1+\alpha\tau_m \rho^{-\delta})\right),
\]
while, on $\Cc_{4\rho, \tau_m/2}$, that
\[
  G_1(x, \tau) \geq \exp\left(8\rho^2(1+\alpha\tau_m\rho^{-\delta})\right).
\]
So, \eqref{eq:xycarleman1c} implies that
 \begin{align}\label{eq:xycarleman1d}
\begin{split}
&\|(|\Xb| + |\Yb|)\|_{\Cc_{4\rho, \frac{\tau_m}{2}}} \le Ne^{-4\rho^2\left(1 + \alpha \tau_m \rho^{-\delta}\right)}
\end{split}
\end{align}
for all $\alpha > \alpha_0$, $\rho > \srad$, and $m > 0$. Sending first $\alpha\longrightarrow \infty$ and then $\rho\longrightarrow \srad$
in \eqref{eq:xycarleman1d}, we obtain that
gives that $\Xb \equiv 0$ and $\Yb\equiv 0$ on $\Cc_{4\srad, \tau_m/2}$ for all $m > 0$. Sending at last
$m\longrightarrow\infty$, we conclude that $\Xb \equiv 0$ and $\Yb \equiv 0$ on $\Cc_{4\srad, \lambda T}$.
\end{proof}

\section{Backward uniqueness for terminally conical  Ricci flows}\label{sec:rfbu}
We now return to the Ricci flow, though, for convenience, we will continue to work with the backward time parameter $\tau$ and consider solutions $g(\tau)$ to the backward Ricci flow \eqref{eq:brf} which flow out of a cone. We will say that a solution $g(\tau)$ to \eqref{eq:brf} on $\Cc_{\rad, T}$ \emph{emanates from $\hat{g}$ with quadratic curvature decay} if
\[
 g(x, 0) = \hat{g}(x)
\]
on $\Cc_{\rad}$ and
\begin{equation}\label{eq:rfquaddecay2}
      |\Rm|_g(x, \tau) \leq K r^{-2}(x)
 \end{equation}
 on $\Cc_{\rad, T}$ for some constant $K$.

In this section, we will use Theorem \ref{thm:xysysbu} to prove the the following uniqueness result.
\begin{theorem}\label{thm:rfbu}
 Suppose $g(\tau)$ and $\gt(\tau)$ are solutions to \eqref{eq:brf} on $\Cc_{\rad, T}$
 which emanate from $\hat{g}$ at $\tau=0$ with quadratic curvature decay.
Then $g(x, \tau) \equiv \tilde{g}(x, \tau)$ on $\Cc_{\rad, T}$.
\end{theorem}
This statement generalizes that of Theorem 2.2 in \cite{KotschwarWangConical}. Here, we do not assume that $g(\tau)$ and $\gt(\tau)$ are (a priori) asymptotically conical shrinking self-similar solutions, and we assert that the two solutions agree on all of $\Cc_{\rad, T}$ rather than on some smaller time interval on some more remote neighborhood of infinity. The latter is particularly important for our application to Theorem \ref{thm:terminalcone2}.

\subsection{Quadratic decay of derivatives of curvature}

We intend to use one of the solutions in Theorem \ref{thm:rfbu} as a background metric, and so we first verify that a solution to \eqref{eq:brf} which emanates from $\hat{g}$ with quadratic curvature decay emanates smoothly from $\hat{g}$ in the sense of Definition \ref{def:emsmooth}. Given our assumption of quadratic curvature decay, Shi's  local derivative estimates guarantee that we will have uniform quadratic decay of all derivatives of curvature
on any region $\Cc_{\rad^{\prime}, T^{\prime}}\subset \Cc_{\rad, T}$ with $T^{\prime}  < T$ and $R^{\prime} > R$.

\begin{proposition}\label{prop:rfemsmooth} Suppose $g(\tau)$ is a solution to \eqref{eq:brf} on $\Cc_{\rad, T}$ which emanates from $\hat{g}$ with quadratic curvature decay.
 Then, for any $\rad^{\prime} > \rad$ and $0 < T^{\prime} < T$, $g(\tau)$ emanates smoothly from $\hat{g}$ on $\Cc_{\rad^{\prime}, T^{\prime}}$ in the sense of Definition \ref{def:emsmooth}. Moreover, for any $m$, there is a constant $N_m$
 such that
\[
   |\nabla^{(m)}\Rm|(x, \tau) \leq N_mr^{-2}(x)
\]
on $\Cc_{\rad^{\prime}, T^{\prime}}$.
\end{proposition}
\begin{proof} The uniform curvature bound \eqref{eq:rfquaddecay2} together with Shi's local estimates (interpreted
for the backward Ricci flow) imply that for each $\rad^{\prime} > \rad$, $0 < T^{\prime} <T$, and $m\geq 0$, there are constants $N_m$ such that
\begin{equation}\label{eq:rfderivativequaddecay}
 |\nabla^{(m)}\Rm|(x, \tau) \leq N_m r^{-2}(x)
\end{equation}
on $\Cc_{\rad^{\prime}, T^{\prime}}$. In particular, using the standard evolution equations for the curvature tensor
and its covariant derivatives under the flow, the bounds \eqref{eq:rfderivativequaddecay} for $0 \leq m \leq 7$
imply that there is some constant $N$ such that
\[
   |\nabla^{(m)}\partial_{\tau}^{(l)}\Rc|(x, \tau) \leq Nr^{-2}(x)
\]
on $\Cc_{\rad^{\prime}, T^{\prime}}$ for $0\leq l \leq 2$ and $0 \leq m \leq 3$.
  Additionally,
\[
\Rc(x, 0)(\nabla r, \cdot) = \widehat{\Rc}(x, 0)(\delh r, \cdot) \equiv 0
\]
on $\Cc_{\rad}$
so both \eqref{eq:geninitial} and \eqref{eq:genquaddecay} are satisfied with $\beta = \Rc$. Thus $g(\tau)$ emanates
smoothly from $\hat{g}$ on $\Cc_{\rad^{\prime}, T^{\prime}}$ in the sense of Definition \ref{def:emsmooth}.
\end{proof}

For our purposes, it will always be enough to work on a slightly smaller domain $\Cc_{\rad^{\prime}, T^{\prime}}$ on which
the first several derivatives of $\Rm$ satisfy uniform quadratic decay bounds. We will say that a solution $g(\tau)$ to \eqref{eq:brf} has \emph{quadratic curvature decay with derivatives} on $\Cc_{\rad, T}$ if there is some constant $K$
such that
\begin{equation}\label{eq:quadwderiv}
        |\nabla^{(m)}\Rm|(x, t) \leq Kr^{-2}(x)
\end{equation}
for all $0\leq m \leq 7$.

Thus, a solution $g(\tau)$ to \eqref{eq:brf} which emanates from $\hat{g}$ and has quadratic decay with derivatives emanates smoothly from $\hat{g}$ in the sense we considered Section \ref{sec:acmetrics}.

\subsubsection{A remark on the shrinking self-similar case}\label{sec:shrinkercase} When $g(\tau)$ is a solution to \eqref{eq:brf} on $\Cc_{\rad, 1}$ arising from a normalized shrinking soliton $(M, f, g, 1)$
asymptotic to $\hat{g}$ on $\Cc_\rad$, and satisfying $g(0) = \hat{g}$ and $g(1) = g$ as described at the end of Section \ref{sec:intro}, we in fact have estimates of the form
\[
  |\nabla^{(m)}\Rm|(x, t) \leq K_mr^{-(m+2)}(x)
\]
on $\Cc_{\rad^{\prime}, 1}$ for all $m \geq 0$ and $\rad^{\prime} > \rad$. (Thus it will follow from Theorem \ref{thm:terminalcone2} that \emph{all} solutions to \eqref{eq:brf} on $\Cc_{\rad, 1}$ which emanate from $\hat{g}$ with quadratic curvature decay in fact satisfy such rates.) These rates are optimal as can be seen by considering the trivial shrinking soliton structure on a Ricci-flat cone. However, we will not need more than the quadratic decay of the derivatives of curvature to prove the result we are after.

\subsection{The PDE-ODE System}\label{sec:rfpdeode}
Our next task is recast the statement of Theorem \ref{thm:rfbu} as one for a system of inequalities
to which we may apply Theorem \ref{thm:xysysbu}. We do so with using
a slightly simpler version of the system from \cite{KotschwarRFBU}, based on the difference $\nabla\Rm - \delt\Rmt$ of the covariant derivatives of the curvature tensors.
In this discussion,
$g(\tau)$ and $\tilde{g}(\tau)$ will denote the two solutions to \eqref{eq:brf} from Theorem \ref{thm:rfbu}.
We will use $g(\tau)$ also as a family of  background metrics. All undecorated quantities
$|\cdot| = |\cdot|_{g(\tau)}$, $\nabla = \nabla_{g(\tau)}$, and $\Delta = \Delta_{g(\tau)}$ will represent the norms,
connections, and Laplace operators induced by $g(\tau)$.

Define
\[
  h \dfn g- \gt \quad \mbox{and} \quad  S = \nabla \Rm - \delt \Rmt.
\]
In Appendix \ref{app:pdeode}, by straightforward computations similar to those in \cite{KotschwarRFBU}, we show that $S$, and the first few covariant derivatives of $h$
satisfy the following schematic evolution equations on $\Cc_{\rad, T}$:
\begin{align}
 \label{eq:dth1}
    D_{\tau} h &= \nabla^{(2)} h + \gt^{-1}\ast (\nabla h)^2 + \big\{\gt^{-1}\ast\Rmt + \Rmt  +  \Rm\big\} \ast h,
\end{align}
\begin{align}
\begin{split}\label{eq:dtdh1}
    D_{\tau}\nabla h &= S
    + \big\{\gt^{-2}\ast\Rmt \ast h + \gt^{-1}\ast \Rmt + \Rm\big\}\ast \nabla h\\
    &\phantom{=}
    + \big\{\gt^{-1}\ast S + \gt^{-1}\ast\nabla \Rm + \nabla \Rm\big\}\ast h,
\end{split}
\end{align}
\begin{align}
\begin{split}\label{eq:dtddh1}
&D_{\tau} \nabla^{(2)} h = \nabla S
+ \bigg\{\gt^{-2}\ast\Rm + \gt^{-1}\ast \Rmt + \Rmt + \Rm\bigg\}\ast\nabla^{(2)} h\\
&\phantom{=} + \bigg\{\gt^{-1}\ast S + \gt^{-2}\ast h \ast S
+ \gt^{-2}\ast \Rmt\ast \nabla h+ \gt^{-3}\ast\Rmt \ast h\ast\nabla h
\\
&\phantom{= + \bigg\}}\quad + \gt^{-2}\ast \nabla \Rm \ast h
+ \gt^{-1}\ast \nabla \Rm + \nabla \Rm\bigg\}\ast \nabla h\\
&\phantom{=} + \bigg\{\gt^{-1}\ast\nabla S
+ \gt^{-1}\ast \nabla^{(2)} \Rm + \nabla^{(2)} \Rm\bigg\}\ast h,
\end{split}
\end{align}
and
\begin{align}
\begin{split}
 \label{eq:dts1}
 &(D_{\tau} + \Delta) S =  (\Rm + \Rmt) \ast S\\
 &\phantom{=}+ \bigg\{\gt^{-2} \ast \delt\Rmt + \gt^{-2}\ast\delt\Rmt\ast h + \nabla \Rm \bigg\}\ast \nabla\nabla h \\
&\phantom{=}
 + \bigg\{\gt^{-3}\ast\delt\Rmt \ast h \ast\nabla h   + \gt^{-3}  \ast \delt\Rmt \ast \nabla h
 + \gt^{-1}\ast\nabla\Rm \ast \nabla h\\
&\phantom{=\bigg\{\quad} + \gt^{-2}\ast\delt^{(2)}\Rmt\ast h  + \gt^{-2}\ast \delt^{(2)} \Rmt   \bigg\}\ast\nabla h\\
&\phantom{=}+\bigg\{\nabla \Rm \ast\Rmt + \gt^{-1}\ast \delt\Rmt \ast \Rmt + \gt^{-1}\ast\delt^{(3)}\Rmt
\bigg\}\ast h,
\end{split}
\end{align}
Here $A\ast B$ represents linear combinations of contractions of $A\otimes B$ with respect the metric $g$.
See Proposition \ref{prop:schemev} for the proofs of the above identities.

For our purposes, the important feature of the evolution equations \eqref{eq:dth1} - \eqref{eq:dts1}
is that each of the terms on the right contains a factor of $h$, $\nabla h$, $\nabla^2h$,
$S$, or $\nabla S$. Under our assumptions, the other factors will be bounded
or decay quadratically in space.

Now define the bundles
\[
 \Xc \dfn \pi^*T^5(T^*\Cc_{\rad}) \quad\mbox{and}\quad \Yc \dfn \pi^*T^2(T^*\Cc_{\rad}) \oplus \pi^*T^3(T^*\Cc_{\rad}) \oplus \pi^*T^4(T^*\Cc_{\rad}),
\]
over $\Cc_{\rad, T}$,
where, as before, $\pi:\Cc_{\rad, T}\longrightarrow \Cc_a$ is the projection. Then define the sections $\Xb\in \Gamma(\Xc)$, $\Yb\in \Gamma(\Yc)$ by
\begin{equation}\label{eq:rfxydef}
  \ve{X} = S, \quad\mbox{and}\quad \ve{Y} = (h, \nabla h, \nabla^2h).
\end{equation}

\begin{proposition}\label{prop:xysys} Assume $\rad \geq 1$ and $0 < T\leq 1$ and
suppose that $g(\tau)$ and $\gt(\tau)$ are smooth solutions to \eqref{eq:brf} on $\Cc_{\rad, T}$ which
emanate from $\hat{g}$ at $\tau =0$ and have quadratic decay with derivatives. Let $\Xb$ and $\Yb$ be the associated sections of $\mathcal{X}$ and $\mathcal{Y}$ defined as above. Then there is a constant $N > 0$ depending on the parameter $K$
from \eqref{eq:quadwderiv} such that
\begin{equation}\label{eq:xysysbounds}
 |\Xb| + |\nabla \Xb| + |\Yb| \leq \frac{N\tau}{r^2}
\end{equation}
and
\begin{align}\label{eq:rfxysys}
\begin{split}
    |D_{\tau} \Xb + \Delta \Xb| &\leq \frac{N}{r^2}(|\ve{X}| + |\ve{Y}|)\\
    |D_{\tau} \Yb| &\leq N(|\Xb| + |\nabla \Xb| + |\ve{Y}|)
\end{split}
\end{align}
on $\Cc_{\rad, T}$.
\end{proposition}
\begin{proof}
By Propositions \ref{prop:genmetricdecay} and \ref{prop:rfemsmooth} and our assumption of quadratic curvature decay with derivatives, the metrics $g(\tau)$ and $\gt(\tau)$ are uniformly
equivalent to each other for all $\tau$ (and in particular to the conical metric $\hat{g} = g(0) = \gt(0)$)
and
there is some $N > 0$ such that
\[
    |\delh^{(k)}(g-\hat{g})| + |\delh^{(k)}(\gt-\hat{g})| \leq \frac{N\tau}{r^2} \quad\mbox{and}\quad
     |\nabla^{(l)}\Rm| + |\delt^{(l)}\Rmt| \leq \frac{N}{r^{2}}
\]
on $\Cc_{\rad, T}$ for all $0 \leq k \leq 3$, and $0 \leq l \leq 7$. Thus also
\[
 |\delh^{(k)} h| \leq \frac{N\tau}{r^2}
\]
for $0 \leq k \leq 3$.
Using \eqref{eq:genhckest2}, we then have
\[
    |\nabla h| \leq  |\delh h| + C|\delh g||h|
    \leq \frac{N\tau}{r^{2}}, \quad\mbox{and}\quad  |\delt h| \leq  |\delh h| + C|\delh \gt||h|
    \leq \frac{N\tau}{r^{2}},
\]
and, arguing similarly for higher derivatives, that
\[
    |\nabla^{(k)}h| + |\delt^{(k)}h| \leq
    \frac{N\tau}{r^{2}},
\]
for $0 \leq k \leq 3$.
Hence we have in particular that
\begin{equation}\label{eq:ybounds}
    |\Yb| \leq \frac{N\tau}{r^{2}}
\end{equation}
on $\Cc_{\rad, T}$ for some $N$.

Next, using the evolution equations  for $\nabla \Rm$ and $\delt\Rmt$ (see \eqref{eq:drmev}) and the assumption \eqref{eq:quadwderiv} of quadratic curvature decay with derivatives, we see that $S = \nabla \Rm -\delt\Rmt$ satisfies
$|\partial_{\tau} S|_{\hat{g}} \leq Nr^{-2}$,
so that
\begin{equation}\label{eq:xbounds}
|\Xb| \leq \frac{N\tau}{r^{2}}.
\end{equation}
Similarly, using the evolution equations for $\nabla^{(2)}\Rm$ and $\delt^{(2)}\Rmt$
and estimating
\[
     |\nabla \Xb| \leq |\nabla^{(2)}\Rm -\delt^{(2)}\Rmt| + C|\gt^{-1}||\nabla h||\delt\Rmt|.
\]
we see that
\begin{equation}\label{eq:dxbounds}
  |\nabla \Xb| \leq \frac{N\tau}{r^2}
\end{equation}
Combining \eqref{eq:ybounds}, \eqref{eq:xbounds}, and \eqref{eq:dxbounds}
then gives \eqref{eq:xysysbounds}.

For the system \eqref{eq:rfxysys}, we then use the bounds just derived on the components of $\ve{X}$
and $\ve{Y}$ together with the evolution equations \eqref{eq:dth1} - \eqref{eq:dts1} to see that
\begin{align*}
    |D_{\tau} h| &\leq N(|h| + |\nabla h| + |\nabla^{(2)} h|),
\end{align*}
\begin{align*}
    |D_{\tau}\nabla h| &\leq N(|h| + |\nabla h| + |S|),
\end{align*}
\begin{align*}
|D_{\tau} \nabla^{(2)} h| &\leq N(|h| + |\nabla h| + |\nabla^2h| + |\nabla S|),
\end{align*}
and
\begin{align*}
|(D_{\tau} + \Delta) S| &\leq \frac{N\tau}{r^2}(|h| + |\nabla h| + |\nabla^{(2)}h| + |S|)
\end{align*}
on $\Cc_{\rad, T}$ for some constant $N$.
\end{proof}

\subsection{Proof of Theorem \ref{thm:rfbu}}

To prove Theorem \ref{thm:rfbu} --- which asserts that the solutions agree not only on some $\Cc_{\rad^{\prime}, T^{\prime}}\subset \Cc_{\rad, T}$ as Theorem \ref{thm:xysysbu} provides, but on all of $\Cc_{\rad, T}$ --- we will need to iteratively apply Theorem \ref{thm:xysysbu}. For this iteration we will need the following refinement of the usual statement of the real-analyticity
for solutions to the Ricci flow (stated here for the backward Ricci flow).
\begin{lemma}[Theorem 3.1, \cite{KotschwarWangIsometries}]\label{lem:analyticstructure} Suppose $g(\tau)$ is a smooth solution to \eqref{eq:brf} on $M\times [0, T]$. Then there exists a unique real-analytic structure $\Ac$
 relative to which $(M, g(\tau))$ is analytic for all $\tau\in [0, T)$. This structure is generated by the atlas of $g(\tau_0)$-normal coordinate charts for any $\tau_0\in [0, T)$.
 \end{lemma}
The real-analyticity of the time slices $(M, g(t))$ of a solution to Ricci flow for $t > 0$ was proven by Bando \cite{Bando} in the compact case, and his
proof extends verbatim to the case of complete solutions of bounded curvature. In
\cite{KotschwarRFAnalyticity}, we localized Bando's argument to apply to all smooth solutions. The point of Lemma \ref{lem:analyticstructure} is that the time slices $(M, g(\tau))$ are real-analytic relative
to a common (time-independent) atlas.

\begin{proof}[Proof of Theorem \ref{thm:rfbu}] As before, write $T_* = \min\{T, 1\}$.
By Proposition \ref{prop:rfemsmooth}, the solutions $g(\tau)$ and $\gt(\tau)$ have quadratic decay with derivatives
on $\Cc_{2\rad, T_*/2}$ and emanate smoothly from $\hat{g}$ at $\tau=0$. Let $h = g- \gt$ and form the sections
\[
  \ve{X}= S = \nabla \Rm - \delt\Rmt \quad \mbox{and} \quad \ve{Y} = (h, \nabla h, \nabla^{(2)}h),
\]
as in Section \ref{sec:rfpdeode}. By Proposition \ref{prop:xysys}, the pair $\ve{X}$ and $\ve{Y}$, together with the family of background metrics $g(\tau)$, meet the hypotheses of Theorem \ref{thm:xysysbu} on $\Cc_{2\rad, T_*/2}$. Thus, there are $\srad \geq 2\rad$ and $\lambda = \lambda(n) \in (0, 1)$ such that
\[
   \ve{X}(x, \tau) \equiv 0 \quad  \mbox{and} \quad \ve{Y}(x, \tau) \equiv 0,
\]
on $\Cc_{\srad, \lambda T_*/2}$.  In particular, $g(x, \tau) \equiv \gt(x, \tau)$ on $\Cc_{\srad, \lambda T_*/2}$.  According to Lemma \ref{lem:analyticstructure}, for each $\tau\in [0, T)$, $g(\tau)$ and $\gt(\tau)$ are both real-analytic with respect to a common atlas on $\Cc_{\rad}$
(since $g(0) = \hat{g} = \gt(0)$, they are real-analytic, e.g., with respect to the atlas generated by $\hat{g}$-normal coordinates). Thus, since  $g(x, \tau) \equiv \gt(x, \tau)$ on $\Cc_{\srad}\times\{\tau\}$ for each such $\tau$, we must actually have $g(x, \tau) \equiv \gt(x, \tau)$ --- and therefore also $\Xb \equiv 0$ and $\Yb \equiv 0$ --- on $\Cc_{\rad, \lambda T_*/2}$.

Define
\[
    \mathcal{A} = \{\,a \in [0, T]\,|\, g\equiv \gt \ \mbox{ on } \ \Cc_{\rad, a}\},
\]
and put $A = \sup \Ac$.  We know $A \geq \lambda T_*/2 > 0$. If $A = T$, then $g \equiv \gt$ on all of
$\Cc_{\rad, T}$ by continuity and the proof is complete.

So suppose $A < T$. By continuity, the sections $\Xb$ and $\Yb$ vanish identically on $\Cc_{\rad, A}$.
We cannot simply translate $\ve{X}(\tau)$ and $\ve{Y}(\tau)$ in time and apply Theorem \ref{thm:xysysbu} to the translates
as doing so would entail also translating the background metric $g(\tau)$, and the time-translate of $g(\tau)$ would not in general be guaranteed to
emanate smoothly from the cone at the new initial time. (This property of the metrics was used to derive the bounds \eqref{eq:eta1est} - \eqref{eq:eta4est} which were crucial in the error estimates for the weight functions in Propositions \ref{prop:carlemanest1}, \ref{prop:pdecarleman2}, and \ref{prop:odecarleman2}.) But we can more or less follow this strategy
if we first tinker a bit with the background metric.

Let
\[
T_0 \dfn \min\{T-A, 1\}/2, \quad 0 < \delta < \min\left\{\frac{A}{2}, \frac{\lambda T_0}{1-\lambda}\right\},
\]
and define $T_1 \dfn T_0 + \delta$. With these choices, we have
\begin{equation}\label{eq:deltaconseq}
   0 < T_1 < 1, \quad A+ T_1 - \delta \leq T, \quad\mbox{and}\quad \lambda T_1 > \delta.
\end{equation}
Then let $\gamma:[0, T_1]\longrightarrow [A-\delta, A+ T_1 -\delta]$ be a smooth, strictly increasing function satisfying
\[
   \gamma(\tau) = \left\{\begin{array}{rl}
                          \tau &  \tau \in [0, \delta/3],\\
                          \tau + A- \delta & \tau \in [2\delta/3, T_1].
                         \end{array}\right.
\]
Define
\[
  \bar{g}(\tau) \dfn g(\gamma(\tau)), \quad \ve{\bar{X}}(\tau) \dfn \Xb(\gamma(\tau)), \quad\mbox{and}\quad \ve{\bar{Y}}(\tau) \dfn\Yb(\gamma(\tau)).
\]
Then, since $\bar{\beta} \dfn \frac{1}{2}\partial_{\tau}\bar{g}$ satisfies
\[
 \partial_{\tau}^{(k)}\bar{\beta}(x, \tau) = \partial_{\tau}^{(k)}(\gamma^{\prime}(\tau)\beta(\gamma(\tau), x))
\]
and $\delb_{\bar{g}(\tau)}  = \nabla_{g(\gamma(\tau))}$, it follows that
$\gb(\tau)$ has quadratic decay with derivatives on $\Cc_{2\rad, T_1}$
and emanates smoothly from $\hat{g}$ at $\tau = 0$. On the intervals $[0, \delta/3]$ and $[2\delta/3, T_1]$, $\bar{g}(\tau)$ is an exact solution to the backward Ricci flow, and it coincides with $g(\tau + A-\delta)$
on the latter interval. Since
\[
 \ve{\bar{X}}(\tau) = \left\{\begin{array}{cl}
                          0 &  \tau \in [0, \delta]\\
                          \Xb(\tau + A- \delta) & \tau \in [2\delta/3, T_1]
                         \end{array}\right.,
\]
and
\[
\ve{\bar{Y}}(\tau) = \left\{\begin{array}{cl}
                          0 &  \tau \in [0, \delta]\\
                          \Yb(\tau + A- \delta) & \tau \in [2\delta/3, T_1]
                         \end{array}\right.,
\]
we see that $\ve{\bar{X}}$ and $\ve{\bar{Y}}$ satisfy
\begin{equation*}
 |\ve{\bar{X}}| + |\delb\ve{\bar{X}}| + |\ve{\bar{Y}}| \leq \frac{N}{r^2},
\end{equation*}
and
\begin{align*}\label{eq:rfxysys}
\begin{split}
    |\bar{D}_{\tau} \ve{\bar{X}} + \bar{\Delta}\ve{\bar{X}}| &\leq \frac{N}{r^2}(|\ve{\bar{X}}| + |\ve{\bar{Y}}|)\\
    |\bar{D}_{\tau} \ve{\bar{Y}}| &\leq N(|\ve{\bar{X}}| + |\delb\ve{\bar{X}}| + |\ve{\bar{Y}}|),
\end{split}
\end{align*}
on $\Cc_{2\rad, T_1}$, where in the above the norms are taken relative to $\bar{g}$.
Thus $\ve{X}(\tau)$ and $\ve{Y}(\tau)$ meet the hypotheses of Theorem \ref{thm:xysysbu}
on $\Cc_{2\rad, T_1}$ relative to the background metric $\bar{g}(\tau)$.

Applying that Theorem, we obtain that
$\ve{\bar{X}}\equiv 0$ and $\ve{\bar{Y}}\equiv 0$ on $\Cc_{\srad, \lambda T_1}$ for some $\srad \geq 2\rad$.
In particular,
\[
   g(x, \tau+A-\delta) \equiv \gt(x, \tau + A - \delta)
\]
on $\Cc_{\srad}\times [\delta, \lambda T_1]$, so $g(x, \tau) \equiv \gt(x, \tau)$ on $\Cc_{\srad, T_2}$ where $T_2 = \lambda T_1 + A -\delta$. Appealing to Lemma \ref{lem:analyticstructure} again as above, it follows that we  in fact have $g(x, \tau) \equiv \gt(x, \tau)$ on $\Cc_{\rad, T_2}$. So $T_2\in \Ac$. But, by \eqref{eq:deltaconseq}, our choice of $\delta$ implies that
$T_2 = \lambda_1 T + A - \delta > A$, so this contradicts that $A = \sup \Ac$.
\end{proof}

\section{An asymptotically conical shrinking soliton is a gradient soliton}\label{sec:gradsol}

Before we can use Theorem \ref{thm:rfbu}  to prove Theorem \ref{thm:terminalcone2}, we need to first establish one more fact.
Recall that a \emph{Ricci soliton structure} on a smooth manifold $M$ is a $4$-tuple $(M, g, X, \lambda)$
consisting of a Riemannian metric $g$, a smooth vector field $X$, and a constant $\lambda$ such that
\begin{equation}\label{eq:rcsol}
 \Rc(g) + \frac{1}{2}\mathcal{L}_X g= \frac{\lambda}{2} g.
\end{equation}
The structure is \emph{shrinking} if $\lambda  > 0$ and \emph{gradient} if $X= \nabla f$ for some smooth function $f$ on $M$.

Let us consider the case that $\lambda = 1$, and let $\Phi_{\tau}:\Cc_{\rad}\longrightarrow \Cc_{\rad}$ be the map $\Phi_{\tau}(r, \sigma) = (r/\sqrt{\tau}, \sigma)$. Then
\[
   \pd{}{\tau} \Phi_{\tau} = -\frac{r}{2\tau}\pd{}{r}, \quad \Phi_{1} = \operatorname{Id}.
\]
If the soliton structure $(\Cc_{\rad}, g, X, 1)$ is normalized so that $X = \frac{r}{2}\pd{}{r}$,
then $g(\tau) = \tau \Phi_\tau^*g$ solves the backward Ricci flow on $\Cc_{\rad}\times (0, 1]$ with $g(1) =g$.  If in addition $g(\tau)$ has uniform quadratic curvature decay and  converges to $\hat{g}$ locally smoothly on $\Cc_{\rad}$ as $\tau \searrow 0$, we say that $(\Cc_{\rad}, g, X, 1)$ is \emph{dynamically asymptotic} to $\hat{g}$.

It follows from Proposition 2.1 in \cite{KotschwarWangConical} (see Proposition 2.1 in \cite{KotschwarWangIsometries}) that if a gradient shrinker $(\Cc_{\rad}, g, f, 1)$ is $C^2$-asymptotic to $\hat{g}$ in the sense described in at the end of Section \ref{sec:intro},
then there is an injective local diffeomorphism $\varphi:\Cc_{\rad^{\prime}}\longrightarrow\Cc_{\rad}$ such that
$(\Cc_{\rad^{\prime}}, \varphi^*g, \varphi^*f, 1)$ is dynamically asymptotic to $\hat{g}$.

Using Theorem \ref{thm:rfbu}, we will be able to show that a solution to the backward Ricci flow on $\Cc_{\rad, T}$ which emanates smoothly from $\hat{g}$ and has quadratic curvature decay must coincide
on $\Cc_{\rad}\times (0, T]$ with the self-similar solution associated to a shrinking soliton structure $(\Cc_{\rad}, g, X, T^{-1})$ which is dynamically asymptotic to $\hat{g}$. We show next that such a structure
must be gradient.

\begin{theorem}\label{thm:xexact}
  Suppose that $(\Cc_{\rad}, \bar{g}, \bar{X}, 1)$ is a shrinking Ricci soliton which is dynamically asymptotic to $\hat{g}$. Then
there is $f\in C^{\infty}(\Cc_{\rad})$ such that $\operatorname{grad}_{\bar{g}} f = \bar{X}$. In fact, we may take
  \begin{equation}\label{eq:potentialexp}
    f(r, \sigma) = \frac{r^2}{4}\left(1 + 8\int_r^{\infty}s^{-3}\bar{R}(s, \sigma)\,ds\right).
  \end{equation}
\end{theorem}
\begin{proof} Let $\Phi_{\tau}:\Cc_{\rad}\longrightarrow\Cc_{\rad}$ be the map $\Phi_{\tau}(r/\sqrt{\tau}, \sigma)$ as above. The assumption that $(\Cc_{\rad}, \bar{g}, \bar{X}, 1)$ is dynamically asymptotic to $\hat{g}$
entails that $\bar{X} = \frac{r}{2}\pd{}{ r}$ on $\Cc_{\rad}$ and that the
associated solution to \eqref{eq:brf} extends to a smooth solution
on all $\Cc_{\rad, 1}$ defined by
\[
   g(\tau) = \left\{\begin{array}{rl}
                        \tau \Phi^*_{\tau}\bar{g} & \tau \in (0, 1],\\
                        \hat{g} & \tau = 0.
                    \end{array}\right.
\]
On $\Cc_{\rad}\times (0, 1]$, we have
\begin{equation}\label{eq:soltau}
    \Rc(g(\tau)) + \frac{1}{2}\mathcal{L}_{X_{\tau}}g(\tau) = \frac{1}{2\tau}g(\tau),
\end{equation}
where
\begin{equation}\label{eq:xtau}
   X = X_{\tau} =  \frac{r}{2\tau}\frac{\partial}{\partial r},
\end{equation}
on $\Cc_{\rad}\times (0, 1]$.

Notice that the family of one-forms
\[
   W \dfn W_{\tau} \dfn \tau g(\tau)(X_{\tau}, \cdot) = g(\tau)(\bar{X}, \cdot)
\]
is smooth on all of $\Cc_{\rad, 1}$ and satisfies
\[
 (D_{\tau} + \Delta)W = 0, \quad  W_0 = \hat{g}(X_1, \cdot) =  d(r^2/4).
\]
Indeed, from \eqref{eq:soltau} - \eqref{eq:xtau}, we have
\[
   \pd{}{\tau}(\tau X) = 0, \quad \mbox{and} \quad \Delta X + \Rc(X) = 0,
\]
and hence
\[
  \pd{}{\tau} W = 2\Rc(W) \quad \mbox{and} \quad \Delta W + \Rc(W) = 0,
\]
so
\[
  D_{\tau} W = \partial_{\tau} W - \Rc(W) = - \Delta W.
\]

Then a further computation using that $(D_{\tau} + \Delta)W = 0$ shows that the family of two-forms $A \dfn dW$ satisfies
\begin{equation}\label{eq:wev}
  (D_{\tau} + \Delta)A_{ij} = -2R_{ipqj}A_{pq}, \quad A_0 = 0.
\end{equation}
Since $g(\tau)$ has quadratic curvature decay, we thus have
\[
  |(D_{\tau} + \Delta) A| \leq \frac{N}{r^2}|A|
\]
on $\Cc_{\rad, 1}$.

Applying Theorem \ref{thm:xysysbu} to $A$, we obtain that
$A \equiv 0$ on $\Cc_{\rad_0, \tau_0}$ for some $0 < \tau_0 \leq 1$ and $\rad_0 \geq \rad$.
Now $X$ and $g$ are real-analytic with respect to a common atlas, and therefore so too are $W$ and $A$. So we must actually have $A \equiv 0$ on $\Cc_{\rad, \tau_0}$ But, for $0 <\tau\leq 1$, $A_{\tau}$ agrees with $A_{\tau_0}$ up to a scaling factor and pull-back by a diffeomorphism
on a neighborhood of infinity.  So,  for all $\tau\in [\tau_0, 1]$, we have that $A_{\tau}\equiv 0$ on some neighborhood of infinity and hence (again by analyticity) on all of $\Cc_a\times \{\tau\}$.  Thus $dW = 0$ on $\Cc_{\rad, 1}$.

Now let $V$ be any open neighborhood of $\Sigma$ diffeomorphic to an open euclidean ball, and let $U = (\rad, \infty)\times V$. Then $W_1$ is exact on $U$ and there is a smooth function $f_1$ on $U$ such that $W_1 = df_1$.
Writing $f= f(\tau) = f_1\circ\Phi_{\tau}$, we will have (after potentially adding a constant to $f_1$) that
\begin{equation}\label{eq:fsol}
   \Rc(g(\tau)) + \nabla\nabla f = \frac{g(\tau)}{2}, \quad |\nabla f(\tau)|^2 +  R = \frac{f}{\tau},
\end{equation}
on $U\times (0, 1]$. Since $\Phi_{\tau}$ preserves the neighborhood $U$, it follows as in \cite{KotschwarWangConical} from \eqref{eq:fsol} that $k \dfn \tau f$ will converge locally smoothly to $r^2/4$ on $U$ as $\tau \longrightarrow 0$.

From the latter equation in \eqref{eq:fsol}, we have in particular that $\partial_{\tau} k = \tau R$, so
\[
   k(r, \sigma, 1) = k(r, \sigma, \varepsilon) +\int_{\varepsilon}^{1} s R(r, \sigma, s)\,ds,
\]
for all $0 < \varepsilon < 1$. But since
\[
R(r, \sigma, s) = \frac{R(\Phi_{s}(r, \sigma), 1)}{s} =  \frac{\bar{R}(s^{-1/2}r, \sigma)}{s},
\]
for $s\in (0, 1]$,
this says
\begin{align*}
  f(r, \sigma, 1) &= k(r, \sigma, 1) = \lim_{\varepsilon\to 0}\left(k(r, \sigma, \varepsilon) +\int_{\varepsilon}^{1} s R(r, \sigma, s)\,ds\right)\\
  &= \frac{r^2}{4} +\lim_{\varepsilon\to 0} \int_{\varepsilon}^1 \bar{R}(s^{-1/2}r, \sigma)\,ds\\
  &=  \frac{r^2}{4}\left(1 + 8\lim_{\varepsilon\to 0} \int_{r}^{r/\sqrt{\epsilon}} s^{-3}\bar{R}(s, \sigma)\,ds\right)\\
  &= \frac{r^2}{4}\left(1 + 8\int_r^{\infty}s^{-3}\bar{R}(s, \sigma)\,ds\right)
\end{align*}
on $U$.  This representation shows that there is a globally defined $f\in C^{\infty}(\Cc_{\rad})$ with $df  = W_1 = \bar{X}^{\flat}$, and the proof is complete.
\end{proof}

\subsection{Proof of Theorem \ref{thm:terminalcone2}}
Now we have all that we need to prove the local version of our main result.
\begin{proof}[Proof of Theorem \ref{thm:terminalcone2}]
Note that  the radial dilation $\rho_{\lambda}(r, \sigma) = (\lambda r, \sigma)$ maps the end $\Cc_\rad$ into  $\Cc_{\lambda \rad}\subset \Cc_{\rad}$ for each $\lambda \geq 1$. Thus, for each $\lambda \geq 1$, the family of metrics
\[
   g_{\lambda}(\tau) \dfn \lambda^{-2}\rho_{\lambda}^*g(\lambda^2 \tau)
\]
is well-defined on $\Cc_{\rad, 1/\lambda^{2}}$, solves the backward Ricci flow \eqref{eq:brf} with $g_{\lambda}(0) = \hat{g}$, and has quadratic curvature decay.

Applying Theorem \ref{thm:rfbu} to $g_{\lambda}(\tau)$ and $g(\tau)$, we conclude that
\[
g(\tau) \equiv g_{\lambda}(\tau) \equiv \lambda^{-2}\rho_{\lambda}^*g(\lambda^2\tau)
\]
on $\Cc_{\rad, 1/\lambda^2}$ for each $\lambda \geq 1$. Since
\[
    \Phi_{\tau}(r, \sigma) = (r/\sqrt{\tau}, \sigma) = \rho_{1/\sqrt{\tau}}(r, \sigma),
\]
taking $\lambda = 1/\sqrt{\tau}$ for $\tau\in (0, 1]$ and writing $\bar{g} \dfn g(1)$, we have, in other words, that
\[
  g(\tau) = \tau\Phi_{\tau}^*\bar{g}
\]
on $\Cc_{\rad}\times (0, 1]$. This implies that
\[
  \Rc(g(\tau)) + \frac{1}{2} \mathcal{L}_{X_{\tau}}g(\tau) = \frac{g(\tau)}{2\tau}
\]
on $\Cc_{\rad}\times (0, 1]$, where $X_{\tau} = \frac{r}{2\tau}\pd{}{r}$.

By Theorem \ref{thm:xexact}, $X_1 = \grad_{\bar{g}} \bar{f}$ for some smooth $\bar{f}\in C^{\infty}(\Cc_{\rad})$ so that, taking $f = \bar{f}\circ \Phi_{\tau}$, we have
\[
  \Rc(g(\tau)) + \nabla\nabla f = \frac{g(\tau)}{2\tau}
\]
on $\Cc_{\rad}\times (0, 1]$. In fact, $\bar{f}$ is given by \eqref{eq:potentialexp}.
\end{proof}

\section{A unique continuation property for the soliton condition}\label{sec:solext}
Under the assumptions of Theorem \ref{thm:terminalcone}, we can use Theorem \ref{thm:terminalcone2}
to see that there is a smooth potential function on the end $V$ relative to which $g(1)$ satisfies the gradient shrinking soliton equation.
In order to complete the proof of Theorem \ref{thm:terminalcone},
we will need to know that the gradient shrinking soliton structure
defined on the end $V$ extends to all of $M$. We will prove this in part by way of the following more general extension property for soliton structures.

\begin{theorem}\label{thm:solext}
 Suppose $(M, g)$ is a connected and simply-connected real-analytic manifold and $U \subset M$ is a nonempty connected open set. If $X$ is a smooth vector field on $U$ for which $(U, g, X, \lambda)$ is a Ricci soliton, then $X$ extends uniquely to a vector field on $M$ such that $(M, g, X, \lambda)$
 is a Ricci soliton.  If $(U, g, X, \lambda)$ is gradient, so is $(M, g, X, \lambda)$.
\end{theorem}

Before we begin the proof --- which is modelled closely on the classical extension argument for Killing vector fields --- let us note that
Theorem \ref{thm:solext} admits a very simple proof for gradient soliton structures when the Ricci endomorphism $\Rc:TM\longrightarrow TM$ is nonsingular.
Indeed, if $X = \nabla f$ for some smooth $f$, it follows from the contracted second Bianchi identity that $\nabla R = 2\Rc(\nabla f)$ on $U$. When $\Rc$ is nonsingular, the one-form
\[
   W = \frac{1}{2}\left(\Rc^{-1}(\nabla R)\right)^{\flat}
\]
is well-defined and real-analytic on $M$ and satisfies  $W = df$ on $U$. Thus the tensors
\begin{equation}\label{eq:weq}
   A \dfn dW \quad \mbox{and} \quad T \dfn 2(\Rc(g) -\lambda g) + \mathcal{L}_{W^{\sharp}} g
\end{equation}
vanish on $U$. But $A$ and $T$ are real-analytic and so then must vanish everywhere on $M$. Since $M$ is simply-connected, $W$ is globally exact, and therefore, as $U$ is connected, $f$ will extend to a smooth function on $M$
such that $W = df$ and $(M, g, f, \lambda)$ is a Ricci soliton.

For the general case, we adapt the argument of Nomizu \cite{NomizuKilling} for the extension of Killing vector fields on a real-analytic manifold. The key observation is that, just as for Killing vector fields, a soliton vector field $X$ on $M$ can be recovered
from the value of $X$ and $\nabla X$ at a single point of the manifold.

\subsection{Extending soliton structures along paths}
Given a vector field $X$ on $M$, define $A_X\in \Gamma(\End(TM))$ by
\[
A_X(V) = \nabla_V X,
\]
and $E\in \Gamma(T^*M\otimes \operatorname{End}(TM))$ by
\begin{equation}\label{eq:e2def}
 E_{ij}^k = \nabla^k R_{ij} - \nabla_i R_{j}^k - \nabla_jR_i^k.
\end{equation}
The tensor $E$ is the first-variation $E = \delta(\nabla_g)[\Lc_Xg]$ of the Levi-Civita connection in the direction of $\Lc_Xg = -(2\Rc - \lambda g)$. It also arises in the following analog of Kostant's identity \cite{Kostant} for Killing vector fields.
\begin{lemma}\label{lem:kident}
Suppose $(M, g, X, \lambda)$ is a Ricci soliton. Then $A = A_X$ satisfies
\begin{equation}\label{eq:kident}
  \nabla_V A = R(V, X) + E(V).
\end{equation}
\end{lemma}
 \begin{proof} This is a standard identity; we give the proof here for completeness.
 From the soliton equation we have
\[
  \nabla_i X_j = -\nabla_j X_i + 2\lambda g_{ij} - 2R_{ij},
\]
so
\begin{align*}
 \nabla_i \nabla_j X_k &= \nabla_j \nabla_i X_k - R_{ijkp}X^p = -\nabla_j \nabla_k X_i - 2\nabla_j R_{ki} - R_{ijkp}X^p.
\end{align*}
Applying this identity iteratively, we find that
\begin{align*}
 \nabla_i \nabla_j X_k &= \nabla_k \nabla_i X_j + 2(\nabla_k R_{ij} - \nabla_j R_{ki}) + (R_{jkip} -R_{ijkp})X^p\\
 &= -\nabla_i\nabla_j X_k + 2(\nabla_k R_{ij} - \nabla_j R_{ki} -\nabla_{i}R_{jk})
 + (R_{jkip} - R_{ijkp} - R_{kijp})X^p.
\end{align*}
Then the Bianchi identity gives
\begin{equation}
 \nabla_i \nabla_j X_k = R_{ipjk}X^p + \nabla_k R_{ij} - \nabla_j R_{ki} -\nabla_{i}R_{jk},
\end{equation}
that is,
\begin{equation}\label{eq:ader}
 \nabla_i A_j^k = R_{ipj}^kX^p + E^k_{ij},
\end{equation}
and \eqref{eq:kident} follows.
\end{proof}

It follows from this identity that restrictions of $X$ and $A_X$ to any smooth path
satisfy a closed first-order system of ODE.
\begin{lemma}\label{lem:solsys} For any smooth path $\gamma: [a, b]\longrightarrow M$,
$\Xc(t) \dfn X|_{\gamma(t)}$ and $\Ac(t) \dfn A_{X}|_{\gamma(t)}$ satisfy the system
\begin{align}\label{eq:xasys}
 \begin{split}
    \frac{D}{\partial t}\Xc &=  \Ac(\Tc)\\
    \frac{D}{\partial t}\Ac &= \Rmc(\Tc, \Xc) + \Ec(\Tc)
 \end{split}
\end{align}
along $\gamma$, where $\Tc(t) = \dot{\gamma}(t)$, $\frac{D}{\partial t}$ denotes the covariant derivative along $\gamma$,
and $\Ec$ and $\Rmc$ are the restrictions of $E$ and the curvature endomorphism
$\Rm$ to $\gamma$.

In particular, if $X$ and $X^{\prime}$ are soliton vector fields with respect to $g$ and $\lambda$ on a connected open set $U$
and
\[
    X|_p= X^{\prime}|_p, \quad \nabla X|_p = \nabla X^{\prime}|_p
\]
at some $p\in U$, then $X\equiv X^{\prime}$ on $U$.
\end{lemma}
Of course, the last claim also follows directly from the corresponding statement for Killing vector fields, since the difference of two soliton vector fields
is Killing.

\subsection{Analytic continuations of soliton structures}

We will use Lemma \ref{lem:solsys} to prove that a Ricci soliton structure  admits an analytic continuation along any smooth path. We make this precise with the following definition, modelled on the corresponding definition for continuations of isometries in \cite{LeeRM}.
\begin{definition} Suppose that $(M, g)$ is a Riemannian manifold and $(U, g, X, \lambda)$ is a Ricci soliton structure on an open set $U\subset M$. Let
 $\gamma:[a, b]\longrightarrow M$ be a piecewise smooth path with $\gamma(a) = p_0\in U$.
 We will say that $\{(U_t, X_t)\}_{t\in [a, b]}$ is a \emph{continuation} of $(U, g, X, \lambda)$ along $\gamma$ if
  $X_a \equiv X$ on a neighborhood of $p_0$ and the following conditions are satisfied for each $t\in [a, b]$:
\begin{enumerate}
 \item[(i)] $U_t$ is an open neighborhood of $\gamma(t)$ and $X_t$ is a smooth vector field on $U_t$ such that
 $(U_t, g|_{U_t}, X_t, \lambda)$ is a Ricci soliton structure.
 \item[(ii)] There is a neighborhood $J\subset [a, b]$ of $t$ such that for all $s\in J$, $\gamma(s)\in U_s$ and $X_s \equiv X_t$ on $U_s\cap U_t$.
\end{enumerate}
\end{definition}

First we prove that that a soliton structure on connected open set may be continued along any continuous
path originating in that set.
\begin{proposition}[Existence of continuations]\label{prop:acexistence}
 Suppose $(M, g)$ is a real-analytic manifold, $U\subset M$ is a connected open set and $(U, g, X, \lambda)$
 a Ricci soliton structure on $U$. If $\gamma:[0, 1]\longrightarrow M$
 is any piecewise smooth path with $\gamma(0)= p_0\in U$, there exists a continuation $\{(U_t, X_t)\}_{t\in [0, 1]}$ of $(U, g, X, \lambda)$ along $\gamma$.
\end{proposition}
\begin{proof}
 We begin by describing how, for any $p\in M$, $V\in T_pM$, and $W\in \End(T_pM)$, we will define a particular analytic vector field $Y$
 on a normal neighborhood $B$ centered at $p$. Fix geodesic normal coordinates $(x^i)$ on $B$,
 and for each $x \in B$, let $\sigma_x:[0, 1]\longrightarrow V$ be the radial geodesic represented by $\sigma_{x}(t) = tx$ from $0$ to $x$. Then let
\[
       Y(x) = \Xc_{\sigma_x}(1)
\]
where
where $\Xc_{\sigma_{x}}(t)$ and $\Ac_{\sigma_{x}}(t)$ solve the system \eqref{eq:xasys} along $\sigma_{x}$ with
\[
    \Xc_{\sigma_{x}}(0)  = V, \quad \Ac_{\gamma_{x}}(0) = W.
\]
Since the metric $g$ is real-analytic (and hence real-analytic with respect to geodesic normal coordinates),
the pair $(\Xc_{\gamma_x}, \Ac_{\gamma_x})$ solve an inhomogeneous linear system of ODE with real-analytic coefficients
depending on a real-analytic parameter. It follows that $Y$ depends real-analytically on $x$.

Now we return to our main claim and cover $\gamma([0, 1])$ with convex normal neighborhoods $B_t$ centered at $\gamma(t)$ for $t\in [0, 1]$. Shrinking $B_{t_0}= B_0$ if necessary, we may assume that $B_{0}\subset U$. Then we may choose a partition
$\{t_0, t_1, \ldots, t_k\}$ of $[0, 1]$ such that the balls $\{B_{t_0}, \ldots B_{t_k}\}$ cover $\gamma([0, 1])$ and $\gamma([t_{i-1}, t_i])\subset B_{i-1}$ for each $i=1, \ldots, k$.

Since $B_0\subset U$, we have that $(B_0, g, X|_{B_0}, \lambda)$ is a soliton structure on $B_0$. Let $Y_{t_1}$ be the analytic vector field
on $B_{t_1}$ constructed as above using the data $(X|_{\gamma(t_1)}, \nabla X|_{\gamma(t_1)})$ at the center $\gamma(t_1)$.
Now $\gamma(t_1)\in B_0\cap B_{t_1}$. As $X$ is a soliton vector field, the restriction of $X$ and $\nabla X$ to the radial geodesics emanating from
$\gamma(t_1)$ also solve the system \eqref{eq:xasys} for as long as those geodesics remain in $B_{0}$. By the uniqueness assertion in Lemma \ref{lem:solsys}, it follows that $X \equiv Y_{t_1}$ on a neighborhood of $\gamma(t_1)$. Since both vector fields are real-analytic, we in fact have $X \equiv Y_{t_1}$ on all of the connected open set $B_0\cap B_{t_1}$. But the tensor
\[
\mathcal{S}(Y_{t_1}) \dfn 2\Rc(g) + \mathcal{L}_{Y_{t_1}}g - \lambda g
\]
is also real-analytic, and
\[
\mathcal{S}(Y_{t_1}) \equiv 2\Rc(g) + \mathcal{L}_{X}g - \lambda g \equiv 0
\]
on $B_0\cap B_{t_1}$. So $\mathcal{S}(Y_{t_1}) \equiv 0$ on all of $B_{t_1}$, i.e., $(B_{t_1}, g, Y_{t_1}, \lambda)$
is a soliton structure.

Continuing inductively in this fashion, we obtain vector fields $Y_{t_i}$ for $i=1, 2, \ldots, k$ such that $(B_{t_i}, g,  Y_{t_i}, \lambda)$
is a soliton structure and $Y_{t_{i-1}} \equiv Y_{t_{i}}$ on $B_{i-1}\cap B_i$.  Defining
\[
  (U_t, X_t) =\left\{\begin{array}{rl} (B_0, X|_{B_0}) & t\in [0, t_1],\\
                                        (B_{t_i}, X_{t_i}) & t\in (t_{i-1}, t_i], \ i =2, \dots, k,
                     \end{array}\right.
\]
we obtain a continuation $\{(U_t, X_t)\}_{t\in [0, 1]}$ of $(U, g, X, \lambda)$ along $\gamma$.
\end{proof}

Next we prove that the continuation of a soliton structure along a path is unique when it exists.
This argument is analogous to the standard one for the unique continuation of isometries (see, e.g., Lemma 12.1 in \cite{LeeRM}). We include the details here for completeness.

\begin{proposition}[Uniqueness of continuations]\label{prop:acuniqueness}
 Suppose $(M, g)$ is a real-analytic manifold, $(U, g, X, \lambda)$ is
 a Ricci soliton structure on a connected open subset $U\subset M$, and $\gamma:[0, 1]\longrightarrow M$
 is a continuous path with $\gamma(0)= p_0\in U$. Then, if  $\{(U_t, X_t)\}_{t\in [0, 1]}$ and
 $\{(\tilde{U}_t, \tilde{X}_t)\}_{t\in [0, 1]}$ are any two continuations of $(U, g, X, \lambda)$ along $\gamma$, we have
 $\tilde{X}_1 \equiv X_1$ in a neighborhood of $\gamma(1)$.
\end{proposition}

\begin{proof}
Suppose $S\subset [0, 1]$ is the set of $s$ such that $X_s = \tilde{X}_s$ on a neighborhood of $\gamma(s)$. We will show that $S = [0, 1]$. Note first that $0\in S$ since $X_0 \equiv X \equiv \tilde{X}_0$ on a neighborhood of $\gamma(0)$. So $S$ is nonempty.

Suppose that $s\in S$. Then, by the definition of continuation, for all $t$ sufficiently close to $s$, we know that $\gamma(t)\in U_s$ and that $X_t$
agrees with $X_s$ on $U_t\cap U_s$. Likewise, for $t$ sufficiently close to $s$, we have $\gamma(t)\in \tilde{U}_s$ and that $\tilde{X}_t$ agrees with $\tilde{X}_s$ on $\tilde{U}_t\cap \tilde{U}_s$. Since $X_s = \tilde{X_s}$, it follows that for all $t$ sufficiently close to $s$, we have $X_t \equiv \tilde{X}_t$ on a an open neighborhood of $\gamma(s)$. Thus $S$ is open.

Finally, to see that $S$ is closed, suppose that $\{s_i\}_{i=1}^{\infty}$ is a sequence of points in $S$ converging to $s\in [0, 1]$. Then, for each $i$, we have $X_{s_i}\equiv \tilde{X}_{s_i}$ on a neighborhood of $\gamma(s_i)$. For $i$ sufficiently large, we have $\gamma(s_i)\in U_{s}\cap \tilde{U}_{s}$, and $X_{s} \equiv X_{s_i} \equiv \tilde{X}_{s_i}\equiv \tilde{X}_{s}$
on a neighborhood of $\gamma(s_i)$. So
\[
X_s|_{\gamma(s_i)} = \tilde{X}_{s}|_{\gamma(s_i)} \quad\mbox{and}\quad  \nabla X_s|_{\gamma(s_i)} = \nabla \tilde{X}_s|_{\gamma(s_i)}
\]
for such $i$. Taking the limit, it follows by continuity that
\[
X_s|_{\gamma(s)} = \tilde{X}_{s}|_{\gamma(s)} \quad\mbox{and}\quad \nabla X_s|_{\gamma(s)} = \nabla \tilde{X}_s|_{\gamma(s)},
\]
and therefore, by Lemma \ref{lem:solsys}, that
$X_s \equiv \tilde{X}_s$ on a neighborhood of $\gamma(s)$. So $s\in S$ and $S$ is also closed. Thus $S= [0, 1]$.
\end{proof}

Given the existence and uniqueness results for the soliton continuations established  in Proposition \ref{prop:acexistence} and \ref{prop:acuniqueness}, we obtain the following monodromy
theorem with a few simple modifications to the standard argument for isometries. See. e.g., Theorem 12.2 in \cite{LeeRM}.
\begin{proposition}[A monodromy theorem for soliton structures]\label{prop:monodromy}
Suppose $(M, g)$ is a connected real-analytic Riemannian manifold, $U\subset M$ a connected open set,
and $(U, g, X, \lambda)$ a soliton structure on $U$. Let $p\in U$, $q\in M$, and $\gamma_0, \gamma_1:[0, 1]\longrightarrow M$ be any two path-homotopic paths from $p$ to $q$. If $\{(U^0_t, X^0_t)\}_{t\in [0, 1]}$
and $\{(U^1_t, X^1_t)\}_{t\in [0, 1]}$ are any continuations of $(U, X, g, \lambda)$ along $\gamma_0$ and $\gamma_1$, respectively, then $X^0_1$ and $X_1^1$ agree in a neighborhood of $q$.
\end{proposition}

\begin{proof}
 Let $H:[0, 1]\times [0, 1]\longrightarrow M$ be a homotopy between $\gamma_0$ and $\gamma_1$, and write $H_s(t) = H(t, s)$.  According to Proposition \ref{prop:acexistence}, there are continuations $\{(U^s_t, X^s_t)\}_{t\in [0, 1]}$ along $H_s(t)$ for each $s\in [0, 1]$. (We may take the continuations for $s=0$ and $s=1$ to be the given continuations along $\gamma_0 = H_0$ and $\gamma_1 = H_1$.)

 Let $F:[0, 1]\longrightarrow T_qM\oplus \End(T_qM)$ be the map $F(s) = (X_1^s|_q, \nabla X_1^s|_q)$. We claim that that $F$ is locally constant in $s$.
 Indeed, fixing $s\in [0, 1]$, for $s^{\prime}$ sufficiently close to $s$, the image of the curve
 $H_{s^{\prime}}([0, 1])$ will be compactly contained in the union of the sets $\cup_{t\in [0, 1]}U_s^{t}$,
 and the continuation $\{(U^s_t, X_t^s)\}_{t\in [0, 1]}$ will also be a continuation of $(U, g, X, \lambda)$ along $H_{s^{\prime}}$.

 By Proposition \ref{prop:acuniqueness}, $X^s_{1}$ and $X^{s^{\prime}}_1$
 then must agree in a neighborhood of $q$. Then $\nabla X^s_1$ and $\nabla X^{s^{\prime}}_1$ also agree
 in the same neighborhood, so
 \[
    F(s) = (X_1^s|_q, \nabla X_1^s|_q) = (X_1^{s^{\prime}}|_q,
    \nabla X_1^{s^{\prime}}|_q )= F(s^{\prime})
 \]
for all $s^{\prime}\in [0, 1]$ sufficiently close to $s$. So $F$ is locally constant, and therefore, by continuity, constant on $[0, 1]$. In particular, we have $X^0_1|_q = X^1_1|_q$ and $\nabla X^0_1|_q = \nabla X^1_1|_q$. By Lemma \ref{lem:solsys},  $X^0_1$ and $X_1^1$ agree on a neighborhood of $q$.
\end{proof}

\subsection{Extending the soliton structure} We have now all but proven Theorem \ref{thm:solext}.
\begin{proof}[Proof of Theorem \ref{thm:solext}] Denote the soliton structure on $U$ by $(U, g, \bar{X}, \lambda$) and fix $p_0\in U$.
For any $p\in M$, we define $X|_p$ to be $\bar{X}_1$ where $\{(U_t, \bar{X}_{t})\}_{t\in [0, 1]}$ is the continuation of $(U, g, \bar{X}, \lambda)$ along any continuous path $\gamma:[0, 1]\longrightarrow M$ connecting $p_0$
and $p_1$. By Propositions \ref{prop:acexistence} and \ref{prop:monodromy} and the simple-connectedness of of $M$, $X$ is well-defined and smooth on $M$. It follows that $(M, g, X, \lambda)$ is a Ricci soliton structure. Moreover, by Lemma \ref{lem:solsys} and Proposition \ref{prop:acuniqueness}, $X$ agrees with $\bar{X}$ on $U$.

For the second part of the claim, suppose now that $X=\nabla \bar{f}$ on $U$ for some $\bar{f}\in C^{\infty}(U)$. The one-form $X^{\flat}$ associated to $X$ is real-analytic, and thus so is
$dX^{\flat}$. Therefore as $dX^{\flat}\equiv 0$ on $U$, we have $dX^{\flat}\equiv 0$ on $M$. Since $M$
is simply-connected, there is $f\in C^{\infty}(M)$ such that $df = X^{\flat}$. Since $U$ is connected,
we may add a constant to $f$ to achieve that $f|_U \equiv \bar{f}$.
\end{proof}

\begin{remark}
We note that, just as for the continuation of isometries,
the requirement in Theorem \ref{thm:solext} that $M$ be simply-connected can be replaced by the weaker assumption that the fundamental group of the set $U$ surjects onto that of $M$, that is, that the inclusion map $\iota: U\longrightarrow M$
 induces a surjective map $\iota_{*}:\pi_1(U, p_0)\longrightarrow\pi_1(M, p_0)$ for some $p_0\in U$.
 \end{remark}

\section{Proof of Theorem \ref{thm:terminalcone}}\label{sec:mainproofs}

We now have all of the ingredients we need to prove the main result of this paper.

\begin{proof}[Proof of Theorem \ref{thm:terminalcone}] Writing $\tau = -t$, we may regard $g$ as a solution to the backward Ricci flow \eqref{eq:brf} on $M\times (0, 1]$. Our assumptions on the solution $g(\tau)$ imply that the pulled-back metrics $\bar{g}(\tau) = F^*g(\tau)$
on $\Cc_{\rad}\times (0, 1]$
extend to a smooth solution to \eqref{eq:brf} on $\Cc_{\rad, 1}$ with $\bar{g}(0) = \hat{g}$ and quadratic curvature decay.
We may then apply Theorem \ref{thm:terminalcone2} to conclude that
is a function $\bar{f}_1\in C^{\infty}(\Cc_{\rad})$ such that
\[
    \Rc(\bar{g}(1)) + \nabla\nabla_{\bar{g}(1)} \bar{f}_1 = \frac{\gb(1)}{2}, \quad \mbox{and}\quad   \grad_{\bar{g}(1)}\bar{f}_1 = \frac{r}{2}\pd{}{r},
\]
on $\Cc_{\rad}$. Moreover, writing $\bar{\Phi}_{\tau}(r, \sigma) = (r/\sqrt{\tau}, \sigma)$, we have  $\bar{g}(\tau) = \tau \bar{\Phi}_{\tau}^*\bar{g}(1)$ on $\Cc_{\rad}\times (0, 1]$, and that $\bar{f} = \bar{f_1}\circ\bar{\Phi}_{\tau}$ satisfies $\tau \bar{f}\longrightarrow r^2/4$ as $\tau\longrightarrow 0$ on $\Cc_{\rad}$.

Returning our attention to $M$, and defining $\Phi_{\tau} = \bar{\Phi}_{\tau}\circ F$ and $f = \bar{f}\circ F$,
we thus have that $g= g(\tau)$ satisfies
\[
  g(\tau) = \tau \Phi_{\tau}^*g(1) \quad\mbox{and}\quad \Rc(g(\tau)) + \nabla\nabla f = \frac{g}{2\tau}
\]
on $V\times (0, 1]$. Adding a constant to $\bar{f}$, if necessary, we may assume that
\[
   |\nabla f|^2 + R = \frac{f}{\tau}
\]
on $V\times (0, 1]$.

According to \cite{Bando, KotschwarRFAnalyticity}, $(M, g(\tau))$ is real-analytic for $\tau \in (0, 1)$.
Therefore the same is true of its lift $(\tilde{M}, \gt(\tau))$ to the universal cover $\tilde{M}$ of $M$.
Let $\tilde{V}$ be any one of the connected components of the preimage $\pi^{-1}(V)$ of $V$ in $\tilde{M}$ under the covering map $\pi:\tilde{M}\longrightarrow M$. For each $\tau \in (0, 1)$,
$(\tilde{V}, \tilde{g}(\tau), f\circ\pi|_{\tilde{V}}, 1)$ extends to a unique normalized gradient soliton structure $(\tilde{M}, \tilde{g}(\tau), \tilde{f}_{\tau}, \tau^{-1})$  on all of $\tilde{M}$ by Theorem \ref{thm:solext}.

We cannot yet say the same is true for $\tau =1$ as, a priori, we do not know that $(M, g(1))$ is real-analytic (except for on $V$, where we already know it to be a soliton). As a step in this direction, we first show that the soliton structures constructed
above on $\tilde{M}$ for $\tau\in (0, 1)$ are all time-slices of some \emph{common} self-similar Ricci flow.
To do this, we show that each slice can be represented as a rescaled pull-back of the $\tau=1/2$ slice.

Now, $(\tilde{M}, \tilde{g}(1/2), \tilde{f}_{1/2}, 2)$ is a complete gradient shrinking soliton structure. By \cite{ZhangCompletenessGRS}, $\operatorname{grad}_{\gt(1/2)}\ft_{1/2}$ is a complete vector field on $M$,
so if $\Psi_{\tau}$ is the family of diffeomorphisms satisfying the ODE
\[
     \pdtau \Psi_{\tau}  = -\frac{1}{2\tau}\delt_{\gt(1/2)}\ft_{1/2}\circ\Psi_{\tau}, \quad \Psi_{1/2} = \operatorname{Id},
\]
on $\tilde{M}\times (0, \infty)$, then
then $\ch{g}(\tau) = 2\tau\Psi_{\tau}^*\gt(1/2)$ solves the backward Ricci flow \eqref{eq:brf} on $\tilde{M}\times (0, \infty)$, and if $\ch{f}_{\tau} = \tilde{f}_{1/2}\circ \Psi_{\tau}$, then $(\tilde{M}, \ch{g}(\tau), \ch{f}_{\tau}, 1/\tau)$ is a shrinking soliton structure
on $\tilde{M}$ for each $\tau\in (0, \infty)$.

However, $\tilde{f}_{1/2}|_{\tilde{V}} = f_{1/2} \circ \pi|_{\tilde{V}} = f_1 \circ\Phi_{1/2} \circ \pi|_{\tilde{V}}$ on $\tilde{V}$.
Since $\Phi_{\tau}$ satisfies
\[
   \pd{}{\tau}\Phi_{\tau} = -\frac{1}{\tau}\operatorname{grad}_{g(1)}f_1\circ\Phi_{\tau}, \quad \Phi_1 = \operatorname{Id},
\]
and $\Phi_{2\tau} = \Phi_{1/2}^{-1}\circ \Phi_{\tau}$, we have
\[
 \pd{}{\tau}\Phi_{2\tau} = -\frac{1}{\tau}(\Phi_{1/2}^{-1})_{*}\left(\operatorname{grad}_{g(1)}f_1\right)\circ\Phi_{2\tau}= -\frac{1}{2\tau}\left(\operatorname{grad}_{g(1/2)}f_{1/2}\right)\circ\Phi_{2\tau}.
\]
It follows that $\Psi_{\tau} = \pi|_{\tilde{V}}^{-1}\circ \Phi_{2\tau}\circ \pi|_{\tilde{V}}$
for $\tau\in (0, 1)$ and hence that
\[
  \ch{g}(\tau) = 2\tau\Psi_{\tau}^*\gt(1/2) = 2\tau\pi^*(\Phi_{2\tau}^*g(1/2)) =
  \tau\pi^*\Phi_{\tau}^*g(1) = \gt(\tau),
\]
and
\[
 \ch{f} = \ft_{1/2}\circ \Psi_{\tau} = f_1\circ \Phi_{1/2}\circ \pi|_{\tilde{V}}\circ \pi|_{\tilde{V}}^{-1} \circ\Phi_{2\tau} \circ \pi|_{\tilde{V}} = f_1 \circ \Phi_{\tau}\circ\pi|_{\tilde{V}} = \ft,
\]
on $\tilde{V}\times (0, 1)$.

The upshot is that the soliton structures $(\tilde{M}, \ch{g}(\tau), \ch{f}_{\tau}, 1/\tau)$
are also extensions of the soliton structures $(\tilde{V}, \gt(\tau), \ft_{\tau}, 1/\tau)$ for each $\tau\in (0, 1)$, and so by Theorem \ref{thm:solext} must agree with the extensions $(\tilde{M}, \gt(\tau), \ft_{\tau}, 1/\tau)$ we have already
constructed on each time-slice via Theorem \ref{thm:solext} above. In particular, we have
$\gt(\tau) = \ch{g}(\tau)$ on $\tilde{M}\times (0, 1)$. Noting that $\ch{g}(\tau)$ is defined for $\tau\in (0, \infty)$, we may conclude by continuity that $\gt(1) = \ch{g}(1)$
on $\tilde{M}$. Thus, defining $\ft_1 = \ch{f}_1$, we find that
\[
   \Rc(\gt)+ \delt\delt \ft_{\tau} = \frac{\gt}{2\tau}
\]
on $\tilde{M}\times (0, 1]$.

It remains to show that the gradient soliton structure descends to $M$.  Let us write  $\gt= \gt_1$ and $\ft = \ft_1$ for simplicity, so that
\[
   \Rc(\gt) + \delt\delt \ft = \frac{\gt}{2}
\]
on $\tilde{M}$. Note that by our normalization on $(V, g|_V)$, we have $\tilde{R} + |\delt \ft|^2 = \tilde{f}$
on $\tilde{V}$, hence all on all of $\tilde{M}$ by analyticity.

If $(\tilde{M}, \gt)$ is flat, then $\gt(\tau)$ (being self-similar) is flat for $\tau \in (0, 1]$.  This implies that $(V, g(\tau))$ is flat for each $\tau\in (0, 1]$ and therefore so too is $(M, g(\tau))$ by analyticity. In this case we must have that $(M, g(1), |x|^2/4, 1)$ is the Gaussian soliton
on $M = \RR^n$.

So let us assume that $(\tilde{M}, \gt)$ is not flat. Then it follows from \cite{ZhangCompletenessGRS} and the strong maximum principle that $\tilde{R} > 0$ everywhere on $\tilde{M}$ and hence also that $\tilde{f} > 0$.
Let $\phi:\tilde{M}\longrightarrow\tilde{M}$ be any deck transformation, and let
$h = \tilde{f}\circ\phi - \tilde{f}$. We claim that $h\equiv 0$.

To see this, note that --- as both $\tilde{f}$ and $\ft\circ \phi$ are soliton potentials for $\gt$ --- we must at least have
$\delt\delt h \equiv 0$ on $\tilde{M}$. This implies that either $h$ is constant or
$(\tilde{M}, \gt)$ splits as $(\tilde{N}, g_{\tilde{N}})\times (\RR, ds^2)$. Let us assume first that $h$ is not constant. Then $|\delt h|_{\gt}\equiv a > 0$, and $h(x, s) = as + b$ for all $(x, s) \in \tilde{N}\times \RR$ for some $b\in \RR$.  Now, as $\delt \delt h \equiv 0$ implies that $\tilde{\Rc}(\delt h, \cdot) \equiv 0 $, it follows that
$l \dfn \langle \delt \ft,  \delt h\rangle_{\gt}$ satisfies
\[
   \delt l = \delt\delt \ft(\delt h, \cdot) = \frac{1}{2}\delt h
\]
on $M$
and so
\[
  \langle \delt \ft,  \delt h\rangle_{\gt} = \frac{h}{2} + \mathrm{const}
\]
This implies that, along any geodesic $\gamma_x(s) = (x, s)$ parallel to the $\RR$ factor, we have
\begin{equation}\label{eq:fsest}
     \ft(x, s) = \ft(x, 0) + \frac{s^2}{4} + cs + d.
\end{equation}
for some constants $c$ and $d$.

However, as $(\tilde{V}, \gt|_{\tilde{V}})$ is isometric to the asymptotically conical shrinking end $(V, g|_V) \approx (\Cc_{\rad}, \bar{g})$ and $\ft > 0$, we may assume that $\tilde{R} \leq C/\ft$ and $R\leq C/f$ for some constant $C$ on $\tilde{V}$ and $V$, respectively. Since $g$ is not flat, for $c_0$ sufficiently small, the set  $S = \{R= c_0\}\cap V$ is compactly contained in $V$. Let $\tilde{S}$ be any connected component of $\pi^{-1}(S)$ contained in $\tilde{V}$. Then let $p_0 = (x_0, s_0)\in \tilde{S}$ and define $\gamma_{p_0}(x, s) = (x_0, s_0+s)$. Since $\tilde{R}$ is constant along $\gamma_{p_0}$ (as $\delt h$ is Killing), we have in particular that $\gamma_{p_0}(s) \in \tilde{S}\subset \tilde{V}$
for all $s\in \RR$.  On the other hand, we have $\ft(x, s) > s^2/8$ for $s$ sufficiently large by \eqref{eq:fsest}. Thus
\[
    c_0 = \tilde{R}(\gamma_{p_0}(s)) \leq \frac{C}{\ft(x_0, s_0+s)} \leq \frac{8C}{s^2}
\]
for all $s$ sufficiently large, contradicting that $c_0 > 0$.

Thus it must be that $h$ is constant on $\tilde{M}$. Write $h = c$. We claim that $c = 0$.
For, observe that, by the definition of $h$,
\begin{equation}\label{eq:flower}
   \ft(\phi^{(n)}(p)) = \ft(p) + nc
\end{equation}
for all $n\in \ZZ$ and $p\in \tilde{M}$. In particular, if $c\neq 0$, $\phi$ must have infinite order.  But then \eqref{eq:flower} implies that $\ft$ is unbounded below on $\tilde{M}$, contradicting that $\ft > 0$.

Thus we must have $c=0$, and we conclude that $\ft$ is preserved
by every deck transformation. It follows that $\ft$ descends to a smooth function $f\in C^{\infty}$ on $M$ which
coincides with our  original function $f = f(1)$ on $V$ and satisfies
\[
  \Rc(g(1)) + \nabla\nabla_{g(1)} f = \frac{g(1)}{2}
\]
on all of $M$. Since $(M, g(1))$ is complete, $\nabla f$ is complete, and if $\Phi_{\tau}$ is the family of diffeomorphisms of $M$ satisfying
\[
  \pd{}{\tau}\Phi_{\tau} = -\frac{1}{\tau}\operatorname{grad}_{g(1)}f_1 \circ \Phi_{\tau}, \quad \Phi_1 = \operatorname{Id},
\]
then, as above, $g(\tau) = \tau \Phi^*_{\tau}g(1)$ for all $\tau\in (0, 1]$.
\end{proof}

\appendix
\section{Evolution equations for elements of the PDE-ODE system}\label{app:pdeode}

In this appendix, we compute the evolution equations for the components of the PDE-ODE system from Section \ref{sec:rfpdeode}.  Below, $g(\tau)$ and $\tilde{g}(\tau)$ will denote two solutions to the backward Ricci flow on $M\times [0, T]$ on some manifold $M$. As in Section \ref{sec:rfpdeode}, we will use $g(\tau)$ also as a family of  background metrics and write $|\cdot| = |\cdot|_{g(\tau)}$, $\nabla = \nabla_{g(\tau)}$, and $\Delta = \Delta_{g(\tau)}$ for the norms,
connections, and Laplace operators induced by $g(\tau)$.

In these computations, we will make use of the standard ``asterisk'' notation to denote various contractions of tensor products whose internal structure is unimportant for our purposes. However, since we have two metrics $g$ and $\gt$ lurking in the background of this discussion, let clarify here that $A\ast B$ will represent some linear combination of
 contractions of terms of the form $A^*\otimes B^*$ where $A^*$ and $B^*$ are obtained from $A$ and $B$
 by potentially raising and lowering indices with the metric $g$. The coefficients of these terms should depend at most on the dimension $n$.
 In particular, we will not conceal
 any factor of $\gt$ or $\gt^{-1}$ with the asterisk notation.  We will also use the shorthand notations
 \[
        A^k = \underbrace{A\ast A \ast \cdots \ast A}_{k}, \quad \gt^{-k} = (\gt^{-1})^k,
 \]
and
\[
     (A_1 + A_2 + \cdots + A_k) \ast B = A_1\ast B + A_2\ast B +
        \cdots + A_k \ast B
\]
to reduce some of the clutter in the expressions. With these conventions, we have in particular that
\[
        |(A_1 + A_2)\ast B| \leq C(|A_1||B| + |A_2||B|)
\]
for some constant $C =C(n)$.

First we recall some standard identities for the difference of quantities associated to $g$ and $\gt$.
 \begin{lemma}\label{lem:metdiff} Let $g$, $\gt$ be any two metrics and $h = g- \gt$.
Then
\begin{align}\label{eq:invdiff}
  g^{ij} - \gt^{ij} &= -\gt^{ia}g^{jb}h_{ab} = \gt^{-1}\ast h,\\
\label{eq:delgt}
  \nabla_k \gt^{ij} &= -\gt^{ia}\gt^{jb} \nabla_k h_{ab} = \gt^{-2} \ast \nabla h,\\
  \label{eq:rmdiff}
  \Rmt - \Rm & = \nabla^{(2)} h + \gt^{-1}\ast (\nabla h)^2 + \Rmt \ast h,
\end{align}
where $\Rm$ and $\Rmt$ denote the $(4, 0)$ curvature tensors of $g$ and $\gt$.

In addition,
\begin{align}
 \label{eq:d1diff}
 \delt V - \nabla V &= \gt^{-1}\ast V\ast \nabla h
\end{align}
and
\begin{align}
\begin{split}
\label{eq:laplacediff}
\tilde{\Delta} V - \Delta V &= \gt^{-2}\ast V\ast \nabla^{(2)}h
+\bigg\{\gt^{-3}\ast V\ast \nabla h + \gt^{-2}\ast\nabla V \bigg\}\ast\nabla h
\\
&\phantom{=} +\gt^{-1}\ast \nabla^{(2)}V \ast h
\end{split}
 \end{align}
 for any tensor $V$ of rank at least $1$.
\end{lemma}

Using these  identities together with the commutation formula
\begin{equation}\label{eq:dtnablacomm}
 [D_{\tau}, \nabla] V = \nabla \Rc \ast V +  \Rc \ast \nabla V,
\end{equation}
and the evolution equation
\begin{equation}\label{eq:drmev}
    \left(\pdtau + \Delta\right)\nabla \Rm = g^{-2}\ast\nabla \Rm \ast \Rm,
\end{equation}
we now verify Equations \ref{eq:dth1} - \ref{eq:dts1} from Section \ref{sec:rfpdeode}.
 \begin{proposition}\label{prop:schemev}
 The tensors $h = g- \gt$ and $S = \nabla \Rm - \delt\Rmt$ satisfy the schematic equations
\begin{align}
 \label{eq:dth}
    D_{\tau} h &= \nabla^{(2)} h + \gt^{-1}\ast (\nabla h)^2 + \big\{\gt^{-1}\ast\Rmt + \Rmt  +  \Rm\big\} \ast h,
\end{align}
\begin{align}
\begin{split}\label{eq:dtdh}
    D_{\tau}\nabla h &= S
    + \big\{\gt^{-2}\ast\Rmt \ast h + \gt^{-1}\ast \Rmt + \Rm\big\}\ast \nabla h\\
    &\phantom{=}
    + \big\{\gt^{-1}\ast S + \gt^{-1}\ast\nabla \Rm + \nabla \Rm\big\}\ast h,
\end{split}
\end{align}
\begin{align}
\begin{split}\label{eq:dtddh}
&D_{\tau} \nabla^{(2)} h = \nabla S
+ \bigg\{\gt^{-2}\ast\Rm + \gt^{-1}\ast \Rmt + \Rmt + \Rm\bigg\}\ast\nabla^{(2)} h\\
&\phantom{=} + \bigg\{\gt^{-1}\ast S + \gt^{-2}\ast h \ast S
+ \gt^{-2}\ast \Rmt\ast \nabla h+ \gt^{-3}\ast\Rmt \ast h\ast\nabla h
\\
&\phantom{= + \bigg\}}\quad + \gt^{-2}\ast \nabla \Rm \ast h
+ \gt^{-1}\ast \nabla \Rm + \nabla \Rm\bigg\}\ast \nabla h\\
&\phantom{=} + \bigg\{\gt^{-1}\ast\nabla S
+ \gt^{-1}\ast \nabla^{(2)} \Rm + \nabla^{(2)} \Rm\bigg\}\ast h,
\end{split}
\end{align}
and
\begin{align}
\begin{split}
 \label{eq:dts}
 &(D_{\tau} + \Delta) S =  (\Rm + \Rmt) \ast S\\
 &\phantom{=}+ \bigg\{\gt^{-2} \ast \delt\Rmt + \gt^{-2}\ast\delt\Rmt\ast h + \nabla \Rm \bigg\}\ast \nabla\nabla h \\
&\phantom{=}
 + \bigg\{\gt^{-3}\ast\delt\Rmt \ast h \ast\nabla h   + \gt^{-3}  \ast \delt\Rmt \ast \nabla h
 + \gt^{-1}\ast\nabla\Rm \ast \nabla h\\
&\phantom{=\bigg\{\quad} + \gt^{-2}\ast\delt^{(2)}\Rmt\ast h  + \gt^{-2}\ast \delt^{(2)} \Rmt   \bigg\}\ast\nabla h\\
&\phantom{=}+\bigg\{\nabla \Rm \ast\Rmt + \gt^{-1}\ast \delt\Rmt \ast \Rmt + \gt^{-1}\ast\delt^{(3)}\Rmt
\bigg\}\ast h,
\end{split}
\end{align}
\end{proposition}

\begin{proof} Let us temporarily write $U =\Rm -\Rmt$.
For \eqref{eq:dth}, we compute directly that
\begin{align}
\nonumber
    D_{\tau} h_{ij} &= 2(R_{ij}-\Rt_{ij}) + R^{p}_{i}h_{pj} + R^{p}_{j} h_{ip}\\
\nonumber
    &= 2g^{ab}(R_{iabj} -\Rt_{iabj}) + 2(g^{ab} - \gt^{ab})\Rt_{iabj} + R^{p}_{i}h_{pj} + R^{p}_{j} h_{ip}\\
\label{eq:dthvar}
    &= 2g^{ab}U_{iabj} + (\gt^{-1}\ast \Rmt + \Rm)\ast h.
\end{align}
Replacing $U = \Rm -\Rmt$ with the expression in \eqref{eq:rmdiff}, we obtain \eqref{eq:dtdh}.

For \eqref{eq:dtdh}, we use \eqref{eq:dtnablacomm} with \eqref{eq:dthvar} to compute that
\begin{align*}
  &D_{\tau}\nabla h = [D_{\tau}, \nabla h] + \nabla D_{\tau} h\\
  &\qquad = \nabla \Rm \ast h + \Rm \ast \nabla h + \nabla U  + \nabla(\gt^{-1}\ast \Rmt \ast h
  + \Rm \ast h).
\end{align*}
Now, on one hand, we have
\[
 \nabla \Rmt = -S + \nabla\Rm + \gt^{-1}\ast \Rmt \ast \nabla h,
\]
so
\begin{align*}
     \nabla U &= \nabla \Rm - \nabla\Rmt = S + \gt^{-1}\ast \Rmt \ast \nabla h.
\end{align*}
On the other hand, again using the above expression for $\nabla \Rmt$, we have
\begin{align*}
 \nabla(\gt^{-1}\ast \Rmt \ast h + \Rm \ast h) &= (\gt^{-2}\ast \Rmt \ast h + \gt^{-1}\ast \Rmt + \Rm)\ast \nabla h\\
 &\phantom{=}+ (\gt^{-1}\ast S + \gt^{-1}\ast \nabla \Rm + \nabla \Rm) \ast h.
\end{align*}
 Combining these equations yields \eqref{eq:dtdh}.  Equation \ref{eq:dtddh} follows similarly, writing
\[
 D_{\tau}\nabla\nabla h = \nabla (D_{\tau} \nabla h) + \Rm \ast \nabla h + \nabla \Rm \ast h
\]
and using \eqref{eq:dtdh}.

For \eqref{eq:dts}, we first use \eqref{eq:drmev}
to obtain that
\begin{align*}
    &(\partial_{\tau} + \Delta) S =
    (\Deltat - \Delta) \delt\Rmt + \gt^{-1}\ast \delt\Rmt \ast \Rmt \ast h
   + \nabla\Rm \ast U + \Rmt \ast S .
\end{align*}
By \eqref{eq:laplacediff}, the first term on the right is
\begin{align*}
(\Deltat - \Delta)\delt\Rmt &=  \gt^{-2} \ast \delt\Rmt \ast \nabla^{(2)} h +\gt^{-1} \ast \nabla^{(2)} \delt\Rmt \ast h\\
&\phantom{=} + \bigg\{\gt^{-2}\ast \nabla \delt \Rmt  + \gt^{-3}  \ast \delt\Rmt \ast \nabla h\bigg\}\ast\nabla h.
\end{align*}
Using \eqref{eq:d1diff}, we compute that
\[
 \nabla\delt\Rmt = \delt\delt\Rmt + \gt^{-1}\ast \delt\Rmt \ast \nabla h
\]
and
\begin{align*}
 \nabla^{(2)}\delt\Rmt &= \delt^{(3)}\Rmt + \gt^{-1}\ast \delt\Rmt\ast \nabla\nabla h\\
 &\phantom{=}
 + \bigg\{\gt^{-2}\ast\delt\Rmt \ast\nabla h + \gt^{-1} \ast \delt^{(2)}\Rmt \bigg\}\ast \nabla h,
\end{align*}
so
\begin{align*}
&(\Deltat - \Delta)\delt\Rmt =
\bigg\{\gt^{-2} \ast \delt\Rmt + \gt^{-2}\ast\delt\Rmt\ast h \bigg\}\ast \nabla\nabla h \\
&\phantom{=}
 + \bigg\{\gt^{-3}\ast\delt\Rmt \ast h \ast\nabla h   + \gt^{-3}  \ast \delt\Rmt \ast \nabla h
+ \gt^{-2}\ast\delt^{(2)}\Rmt\ast h \\
&\phantom{=\bigg\{\quad}+ \gt^{-2}\ast \delt^{(2)} \Rmt   \bigg\}\ast\nabla h
+ \gt^{-1}\ast\delt^{(3)}\Rmt \ast h.
\end{align*}
Also, by \eqref{eq:rmdiff},
\begin{align*}
\nabla\Rm \ast U &= \nabla \Rm \ast \nabla^{(2)}h + \gt^{-1}\ast\nabla\Rm\ast(\nabla h)^2 +
\nabla \Rm \ast\Rmt \ast h.
\end{align*}
So, putting things together, we have
\begin{align*}
    &(\partial_{\tau} + \Delta) S = \Rmt \ast S + \bigg\{\gt^{-2} \ast \delt\Rmt + \gt^{-2}\ast\delt\Rmt\ast h + \nabla \Rm \bigg\}\ast \nabla\nabla h \\
&\phantom{=}
 + \bigg\{\gt^{-3}\ast\delt\Rmt \ast h \ast\nabla h   + \gt^{-3}  \ast \delt\Rmt \ast \nabla h
 + \gt^{-1}\ast\nabla\Rm \ast \nabla h\\
&\phantom{=\bigg\{\quad} + \gt^{-2}\ast\delt^{(2)}\Rmt\ast h  + \gt^{-2}\ast \delt^{(2)} \Rmt   \bigg\}\ast\nabla h\\
&\phantom{=}+\bigg\{\nabla \Rm \ast\Rmt + \gt^{-1}\ast \delt\Rmt \ast \Rmt + \gt^{-1}\ast\delt^{(3)}\Rmt
\bigg\}\ast h,
\end{align*}
which, since
\[
 D_{\tau} S = \partial_{\tau} S + \Rm \ast S,
\]
gives \eqref{eq:dts}.
\end{proof}

\end{document}